\documentclass[11pt]{amsart}
\usepackage{amsmath}
\usepackage{amscd}
\usepackage{amssymb,latexsym}

\newtheorem{theorem}{Theorem}[section]

\newtheorem{lemma}[theorem]{Lemma}
\newtheorem{proposition}[theorem]{Proposition}

\theoremstyle{definition}

\newtheorem{example}[theorem]{Example}

\newcommand{\Mtwo}{{\sl Macaulay~2\,\,}}
\font\tteight=cmtt8
\def\A{\mathbb{A}}
\def\k{\mathbb{C}}
\def\N{\mathbb{N}}
\def\R{\mathbb{R}}
\def\P{\mathbb{P}}
\def\ZZ{\mathbb{Z}} 

\begin{document}
\bibliographystyle{plain}

\title{Algorithms for the toric Hilbert scheme}
\author{Michael Stillman$^1$}
\address{$^1$Cornell University, Department of Mathematics, 
Ithaca, NY 14853}
\author{Bernd Sturmfels$^2$} 
\address{$^2$UC Berkeley, Department of Mathematics, Berkeley, CA 94720}
\author{Rekha R. Thomas$^3$}
\address{$^3$University of Washington, Department of Mathematics, Seattle,
WA 98195}
\maketitle

\begin{abstract}
The toric Hilbert scheme parametrizes all algebras isomorphic to a
given semigroup algebra as a multigraded vectorspace. All components
of the scheme are toric varieties, and among them, there is a fairly
well understood coherent component. However, it is unknown whether
toric Hilbert schemes are always connected. In this chapter we
illustrate the use of \Mtwo for exploring the structure of toric
Hilbert schemes. In the process we will encounter algorithms from
commutative algebra, algebraic geometry, polyhedral theory and
geometric combinatorics.
\end{abstract}

\section*{Introduction}
Consider the multigrading of the polynomial ring $R =
\k[x_1,\ldots,x_n]$ specified by a non-negative integer $d \times
n$-matrix $A = (a_1,\ldots,a_n)$ such that $degree (x_i) = a_i \in
\N^d$. This defines a decomposition $\, R = \bigoplus_{b \in \N A} R_b
$, where $\N A$ is the subsemigroup of $\N^d$ spanned by
$a_1,\ldots,a_n$, and $R_b$ is the $\k$-span of all monomials $\, x^u
= x_1^{u_1}\cdots x_n^{u_n}$ with degree $Au = a_1 u_1 +\cdots + a_n
u_n = b$.  The {\it toric Hilbert scheme} $\,Hilb_A 
\,$ parametrizes all $A$-homogeneous ideals $I \subset R$ with the
property that $(R/I)_b$ is a $1$-dimensional $\k$-vectorspace, for all
$b \in \N A$. We call such an ideal $I$ an $A$-{\em graded} ideal.
Equivalently, $I$ is $A$-graded if it is $A$-homogeneous and $R/I$ is
isomorphic as a multigraded vectorspace to the semigroup algebra $\,
\k [ \N A ] = R/I_A$, where $$I_A := \,\langle x^u - x^v \, : \, Au =
Av \rangle \subset R$$ is the {\it toric ideal} of $A$. It follows 
from the definition that $A$-graded ideals are generated by monomials
and binomials in $R$. 

Our running example throughout 
this chapter is the following $2 \times 5$-matrix:
\begin{equation}
\label{OurMatrix}
A \quad = \quad 
\left( \begin{matrix}
           1 & 1 & 1 & 1 & 1  \\ 
           0 & 1 & 2 & 7 & 8 
\end{matrix} \right)
\end{equation}
which we input to \Mtwo as a list of lists of 
integers.
\par
\vskip 5 pt
\begingroup
\tteight
\baselineskip=7.4pt
\lineskip=0pt
\obeyspaces
i1\ :\ A\ =\ {\char`\{}{\char`\{}1,1,1,1,1{\char`\}},{\char`\{}0,1,2,7,8{\char`\}}{\char`\}};\ \leavevmode\hss\endgraf
\endgroup
\penalty-1000
\par
\vskip 1 pt
\noindent
The toric ideal of $A$ lives in the multigraded ring $R = \k [a,b,c,d,e]$.
\par
\vskip 5 pt
\begingroup
\tteight
\baselineskip=7.4pt
\lineskip=0pt
\obeyspaces
i2\ :\ \ R\ =\ QQ[a..e,Degrees=>transpose\ A];\ \leavevmode\hss\endgraf
\endgroup
\penalty-1000
\par
\vskip 1 pt
\noindent
\par
\vskip 5 pt
\begingroup
\tteight
\baselineskip=7.4pt
\lineskip=0pt
\obeyspaces
i3\ :\ \ describe\ R\ \leavevmode\hss\endgraf
\penalty-500\leavevmode\hss\endgraf
o3\ =\ QQ\ [a,\ b,\ c,\ d,\ e,\ Degrees\ =>\ {\char`\{}{\char`\{}1,\ 0{\char`\}},\ {\char`\{}1,\ 1{\char`\}},\ {\char`\{}1,\ 2{\char`\}},\ {\char`\{}1,\ 7{\char`\}},\ {\char`\{}1 $\cdot\cdot\cdot$\leavevmode\hss\endgraf
\penalty-500\leavevmode\hss\endgraf
o3\ :\ Adjacent\leavevmode\hss\endgraf
\endgroup
\penalty-1000
\par
\vskip 1 pt
\noindent

We use Algorithm 12.3 in \cite{St2} to compute $I_A$. The first step is
to find a matrix $B$ whose rows generate the lattice $ker_{\ZZ}(A)
:= \{x \in \ZZ^n : Ax = 0 \}$. 

\par
\vskip 5 pt
\begingroup
\tteight
\baselineskip=7.4pt
\lineskip=0pt
\obeyspaces
i4\ :\ \ B\ =\ transpose\ syz\ matrix\ A\ \leavevmode\hss\endgraf
\penalty-500\leavevmode\hss\endgraf
o4\ =\ |\ 1\ -2\ 1\ \ 0\ 0\ |\leavevmode\hss\endgraf
\ \ \ \ \ |\ 0\ 5\ \ -6\ 1\ 0\ |\leavevmode\hss\endgraf
\ \ \ \ \ |\ 0\ 6\ \ -7\ 0\ 1\ |\leavevmode\hss\endgraf
\penalty-500\leavevmode\hss\endgraf
\ \ \ \ \ \ \ \ \ \ \ \ \ \ 3\ \ \ \ \ \ \ \ 5\leavevmode\hss\endgraf
o4\ :\ Matrix\ ZZ\ \ <---\ ZZ\leavevmode\hss\endgraf
\endgroup
\penalty-1000
\par
\vskip 1 pt
\noindent

Although in theory any basis of $ker_{\ZZ}(A)$ will suffice, in
practice it is more efficient to use a {\em reduced} basis
\cite[\S 6.2]{Sch} which can be computed using the {\em basis
reduction} package {\tt LLL.m2}, a built-in feature of \Mtwo. 
The advantage of a reduced basis may not be apparent in small
examples. However, as the size of $A$ increases, it becomes
increasingly important for the termination of Algorithm 12.3. (To
appreciate this, consider the matrix (\ref{non-normal}) from 
Section~4.)

\par
\vskip 5 pt
\begingroup
\tteight
\baselineskip=7.4pt
\lineskip=0pt
\obeyspaces
i5\ :\ \ load\ "LLL.m2";\ \leavevmode\hss\endgraf
\endgroup
\penalty-1000
\par
\vskip 1 pt
\noindent
\par
\vskip 5 pt
\begingroup
\tteight
\baselineskip=7.4pt
\lineskip=0pt
\obeyspaces
i6\ :\ \ B\ =\ transpose\ matrix\ LLL\ syz\ matrix\ A\ \leavevmode\hss\endgraf
\penalty-500\leavevmode\hss\endgraf
o6\ =\ |\ 0\ 1\ \ -1\ -1\ 1\ \ |\leavevmode\hss\endgraf
\ \ \ \ \ |\ 1\ -1\ 0\ \ -1\ 1\ \ |\leavevmode\hss\endgraf
\ \ \ \ \ |\ 2\ 0\ \ -3\ 2\ \ -1\ |\leavevmode\hss\endgraf
\penalty-500\leavevmode\hss\endgraf
\ \ \ \ \ \ \ \ \ \ \ \ \ \ 3\ \ \ \ \ \ \ \ 5\leavevmode\hss\endgraf
o6\ :\ Matrix\ ZZ\ \ <---\ ZZ\leavevmode\hss\endgraf
\endgroup
\penalty-1000
\par
\vskip 1 pt
\noindent

A row $b = b^+ - b^-$ of $B$ is then coded as the binomial
$x^{b^+}-x^{b^-} \in R$, and we let $J$ be the ideal generated by all 
such binomials. 

\par
\vskip 5 pt
\begingroup
\tteight
\baselineskip=7.4pt
\lineskip=0pt
\obeyspaces
i7\ :\ \ toBinomial\ =\ (b,R)\ ->\ (\leavevmode\hss\endgraf
\ \ \ \ \ \ \ \ \ \ top\ :=\ 1{\char`\_}R;\ bottom\ :=\ 1{\char`\_}R;\leavevmode\hss\endgraf
\ \ \ \ \ \ \ \ \ \ scan({\char`\#}b,\ i\ ->\ if\ b{\char`\_}i\ >\ 0\ then\ top\ =\ top\ *\ R{\char`\_}i{\char`\^}(b{\char`\_}i)\leavevmode\hss\endgraf
\ \ \ \ \ \ \ \ \ \ \ \ \ \ \ else\ if\ b{\char`\_}i\ <\ 0\ then\ bottom\ =\ bottom\ *\ R{\char`\_}i{\char`\^}(-b{\char`\_}i));\leavevmode\hss\endgraf
\ \ \ \ \ \ \ \ \ \ top\ -\ bottom);\ \leavevmode\hss\endgraf
\endgroup
\penalty-1000
\par
\vskip 1 pt
\noindent

\par
\vskip 5 pt
\begingroup
\tteight
\baselineskip=7.4pt
\lineskip=0pt
\obeyspaces
i8\ :\ \ J\ =\ ideal\ apply(entries\ B,\ b\ ->\ toBinomial(b,R))\ \leavevmode\hss\endgraf
\penalty-500\leavevmode\hss\endgraf
\ \ \ \ \ \ \ \ \ \ \ \ \ \ \ \ \ \ \ \ \ \ \ \ \ \ \ \ \ \ \ \ \ \ \ \ \ \ \ 2\ 2\ \ \ \ 3\leavevmode\hss\endgraf
o8\ =\ ideal\ (-\ c*d\ +\ b*e,\ -\ b*d\ +\ a*e,\ a\ d\ \ -\ c\ e)\leavevmode\hss\endgraf
\penalty-500\leavevmode\hss\endgraf
o8\ :\ Ideal\ of\ R\leavevmode\hss\endgraf
\endgroup
\penalty-1000
\par
\vskip 1 pt
\noindent
The toric ideal equals $(J : (x_1 \cdots x_n)^\infty)$ which is 
computed via $n$ successive saturations as follows.
\par
\vskip 5 pt
\begingroup
\tteight
\baselineskip=7.4pt
\lineskip=0pt
\obeyspaces
i9\ :\ \ scan(gens\ ring\ J,\ f\ ->\ (J\ =\ saturate(J,f);))\leavevmode\hss\endgraf
\endgroup
\penalty-1000
\par
\vskip 1 pt
\noindent

Putting the above pieces of code together, we get the following
procedure for computing the toric ideal of a matrix $A$.

\par
\vskip 5 pt
\begingroup
\tteight
\baselineskip=7.4pt
\lineskip=0pt
\obeyspaces
i10\ :\ toricIdeal\ =\ (A)\ ->\ (\leavevmode\hss\endgraf
\ \ \ \ \ \ \ \ \ \ n\ :=\ {\char`\#}(A{\char`\_}0);\ \ \leavevmode\hss\endgraf
\ \ \ \ \ \ \ \ \ \ R\ =\ QQ[vars(0..n-1),Degrees=>transpose\ A,MonomialSize=>16];\ \leavevmode\hss\endgraf
\ \ \ \ \ \ \ \ \ \ B\ :=\ transpose\ matrix\ LLL\ syz\ matrix\ A;\leavevmode\hss\endgraf
\ \ \ \ \ \ \ \ \ \ J\ :=\ ideal\ apply(entries\ B,\ b\ ->\ toBinomial(b,R));\leavevmode\hss\endgraf
\ \ \ \ \ \ \ \ \ \ scan(gens\ ring\ J,\ f\ ->\ (J\ =\ saturate(J,f);));\leavevmode\hss\endgraf
\ \ \ \ \ \ \ \ \ \ J\leavevmode\hss\endgraf
\ \ \ \ \ \ \ \ \ \ );\ \leavevmode\hss\endgraf
\endgroup
\penalty-1000
\par
\vskip 1 pt
\noindent

In our example, 
$I_A = \langle cd-be,bd-ae,b^2-ac,a^2d^2-c^3e,c^4-a^3e,bc^3-a^3d,
ad^4-c^2e^3,d^6-ce^5 \rangle$ which we now compute using this procedure.
\par
\vskip 5 pt
\begingroup
\tteight
\baselineskip=7.4pt
\lineskip=0pt
\obeyspaces
i11\ :\ \ I\ =\ toricIdeal\ A\ \leavevmode\hss\endgraf
\penalty-500\leavevmode\hss\endgraf
\ \ \ \ \ \ \ \ \ \ \ \ \ \ \ \ \ \ \ \ \ \ \ \ \ \ \ \ \ \ \ \ \ \ \ \ 2\ \ \ \ \ \ \ \ \ 2\ 2\ \ \ \ 3\ \ \ \ 4\ \ \ \ 3\ \ \ \ \ \  $\cdot\cdot\cdot$\leavevmode\hss\endgraf
o11\ =\ ideal\ (c*d\ -\ b*e,\ b*d\ -\ a*e,\ b\ \ -\ a*c,\ a\ d\ \ -\ c\ e,\ c\ \ -\ a\ e,\ b*c $\cdot\cdot\cdot$\leavevmode\hss\endgraf
\penalty-500\leavevmode\hss\endgraf
o11\ :\ Ideal\ of\ R\leavevmode\hss\endgraf
\endgroup
\penalty-1000
\par
\vskip 1 pt
\noindent
This ideal defines an embedding of $\P^1$ as a degree $8$
curve into $\P^4$. We will see in Section 3 that its toric 
Hilbert scheme $Hilb_A$ has a non-reduced component. 

We recommend \cite{St2} as an introductory reference for the 
topics in this chapter.
The study of toric Hilbert schemes for $d=1$ goes back to
Arnold \cite{Arn} and Korkina et.al.\cite{KPR}, and it was
further developed by Sturmfels  (\cite{St1} and \cite[\S 10]{St2}). 
Peeva and Stillman \cite{PS1} introduced the scheme structure 
which gives the toric Hilbert scheme its universal property,
and from this they derive a formula for the tangent space
of a point on  $\, Hilb_A $. Maclagan recently showed that the 
quadratic binomials in \cite[\S 5]{St1} define the same scheme as the
determinantal equations in \cite{PS1}.
Both of these systems of global equations are 
generally much too big for 
practical computations. Instead, most of our algorithms are based on
the local equations given by Peeva and Stillman in \cite{PS2}
and the combinatorial approach of Maclagan and Thomas in \cite{MT}.

This chapter is organized into four sections and two appendices as
follows. The main goal in Section~1 is to describe an algorithm for
generating all monomial $A$-graded ideals for a given $A$. These
monomial ideals are the vertices of the {\em flip graph} of $A$ whose
connectivity is equivalent to the connectivity of $Hilb_A$. We
describe how all neighbors of a given vertex of this graph can be
calculated. In Section~2, we explain the role of polyhedral geometry
in the study of $Hilb_A$. Our first algorithm tests for {\em
  coherence} in a monomial $A$-graded ideal. We then show how to
compute the polyhedral complexes supporting $A$-graded ideals which in
turn relates the flip graph of $A$ to the {\em Baues graph} of $A$.
For unimodular matrices, these two graphs coincide and hence our
method of computing the flip graph can be used to compute the Baues
graph. Section~3 explores the components of $Hilb_A$ via local
equations around the torus fixed points of the scheme. We include a
combinatorial interpretation of these local equations from the point
of view of integer programming.  The scheme $Hilb_A$ has a {\em
  coherent} component which is examined in detail in Section~4. We
prove that this component is, in general, not normal and that its
normalization is the toric variety of the Gr\"obner fan of $I_A$. We
conclude the chapter with two appendices, each containing one large
piece of \Mtwo code that we use in this chapter. In Appendix 1 we
describe {\tt polarCone.m2} which is a procedure to convert a
generator representation of a polyhedron to an inequality
representation and vice versa. In Appendix 2 we explain a simplified
version of the procedure {\tt minPres.m2}, available in \Mtwo, for
computing minimal presentations of polynomial quotient rings. The main
ingredient of this package is the subroutine {\tt
  removeRedundantVariables} which is what we use in this chapter.

\section{Generating monomial ideals}
We start out by computing the {\it Graver basis} $Gr_A$, which is the
set of binomials in $I_A$ which are minimal with respect to the
partial order $$\, x^u - x^v \,\leq\, x^{u'} - x^{v'} \quad :\iff
\quad \hbox{ $x^u$ divides $x^{u'}$ \ and \ $x^v$ divides $x^{v'}$.}
$$ The set $Gr_A$ is a {\em universal Gr\"obner basis} of $I_A$ and
has its origins in the theory of integer programming \cite{Gra}. It
can be computed using \cite[Algorithm 7.2]{St2}, a \Mtwo version of
which is given below.

\par
\vskip 5 pt
\begingroup
\tteight
\baselineskip=7.4pt
\lineskip=0pt
\obeyspaces
i12\ :\ \ graver\ =\ (I)\ ->\ (\leavevmode\hss\endgraf
\ \ \ \ \ \ \ \ \ \ R\ :=\ ring\ I;\leavevmode\hss\endgraf
\ \ \ \ \ \ \ \ \ \ kk\ :=\ coefficientRing\ R;\leavevmode\hss\endgraf
\ \ \ \ \ \ \ \ \ \ n\ :=\ numgens\ R;\leavevmode\hss\endgraf
\ \ \ \ \ \ \ \ \ \ --\ construct\ new\ ring\ S\ with\ 2n\ variables\leavevmode\hss\endgraf
\ \ \ \ \ \ \ \ \ \ S\ :=\ kk[Variables=>2*n,MonomialSize=>16];\leavevmode\hss\endgraf
\ \ \ \ \ \ \ \ \ \ toS\ :=\ map(S,R,(vars\ S){\char`\_}{\char`\{}0..n-1{\char`\}});\leavevmode\hss\endgraf
\ \ \ \ \ \ \ \ \ \ toR\ :=\ map(R,S,vars\ R\ |\ matrix(R,\ {\char`\{}toList(n:1){\char`\}}));\leavevmode\hss\endgraf
\ \ \ \ \ \ \ \ \ \ --\ embed\ I\ in\ S\leavevmode\hss\endgraf
\ \ \ \ \ \ \ \ \ \ J\ :=\ gens\ toS\ I;\leavevmode\hss\endgraf
\ \ \ \ \ \ \ \ \ \ --\ construct\ the\ toric\ ideal\ of\ the\ Lawrence\ \leavevmode\hss\endgraf
\ \ \ \ \ \ \ \ \ \ --\ lifting\ of\ A\leavevmode\hss\endgraf
\ \ \ \ \ \ \ \ \ \ i\ :=\ 0;\leavevmode\hss\endgraf
\ \ \ \ \ \ \ \ \ \ while\ i\ <\ n\ do\ (\leavevmode\hss\endgraf
\ \ \ \ \ \ \ \ \ \ \ \ \ \ wts\ :=\ join(toList(i:0),{\char`\{}1{\char`\}},toList(n-i-1:0));\leavevmode\hss\endgraf
\ \ \ \ \ \ \ \ \ \ \ \ \ \ wts\ =\ join(wts,wts);\leavevmode\hss\endgraf
\ \ \ \ \ \ \ \ \ \ \ \ \ \ J\ =\ homogenize(J,S{\char`\_}(n+i),wts);\leavevmode\hss\endgraf
\ \ \ \ \ \ \ \ \ \ \ \ \ \ i=i+1;\leavevmode\hss\endgraf
\ \ \ \ \ \ \ \ \ \ \ \ \ \ );\leavevmode\hss\endgraf
\ \ \ \ \ \ \ \ \ J\ =\ ideal\ J;\leavevmode\hss\endgraf
\ \ \ \ \ \ \ \ \ scan(gens\ ring\ J,\ f\ ->\ (J\ =\ saturate(J,f);));\leavevmode\hss\endgraf
\ \ \ \ \ \ \ \ \ --\ apply\ the\ map\ toR\ to\ the\ minimal\ generators\ of\ J\ \leavevmode\hss\endgraf
\ \ \ \ \ \ \ \ \ J\ =\ matrix\ entries\ toR\ mingens\ J;\leavevmode\hss\endgraf
\ \ \ \ \ \ \ \ \ p\ :=\ sortColumns\ J;\leavevmode\hss\endgraf
\ \ \ \ \ \ \ \ \ J{\char`\_}p)\ ;\ \ \leavevmode\hss\endgraf
\endgroup
\penalty-1000
\par
\vskip 1 pt
\noindent
   
   The above piece of code first constructs a new polynomial ring $S$
   in $n$ more variables than $R$. Assume $S = \k [x_1, \ldots, x_n,
   y_1, \ldots, y_n]$. The inclusion map {\tt toS} $: R \rightarrow
   S$ embeds the toric ideal $I$ in $S$ and collects its generators in
   the matrix $J$. A binomial $x^a - x^b$ lies in $Gr_A$ if and only
   if $x^ay^b-x^by^a$ is a minimal generator of the toric ideal in $S$
   of the $(d+n) \times 2n$ matrix $\Lambda(A) := \left (
     \begin{array}{cc} A & 0 \\ I_n & I_n \end{array} \right)$ which
   is called the {\em Lawrence lifting} of $A$. Since $u \in
   ker_{\ZZ}(A) \Leftrightarrow (u,-u) \in ker_{\ZZ} (\Lambda(A))$, we
   use the {\tt while loop} to homogenize the binomials in $J$ with
   respect to $\Lambda(A)$, using the $n$ new variables in $S$. This
   converts a binomial $x^a-x^b \in J$ to the binomial
   $x^ay^b-x^by^a$.  The ideal generated by these new binomials is
   also labeled $J$. As before, we can now successively saturate $J$
   to get the toric ideal of $\Lambda(A)$ in $S$. The image of the
   minimal generators of this toric ideal under the map {\tt toR}
   $: S \rightarrow R$ such that $x_i \mapsto x_i$ and $y_i \mapsto 1$
   is precisely the Graver basis $Gr_A$. This list is then sorted and 
   output.

In our example $Gr_A$ consists of $42$ binomials.
\par
\vskip 5 pt
\begingroup
\tteight
\baselineskip=7.4pt
\lineskip=0pt
\obeyspaces
i13\ :\ Graver\ =\ graver\ I\ \leavevmode\hss\endgraf
\penalty-500\leavevmode\hss\endgraf
o13\ =\ {\char`\{}0,\ 0{\char`\}}\ |\ -cd+be\ -bd+ae\ -b2+ac\ -cd2+ae2\ -a2d2+c3e\ -c4+a2bd\ -c4+a3 $\cdot\cdot\cdot$\leavevmode\hss\endgraf
\penalty-500\leavevmode\hss\endgraf
\ \ \ \ \ \ \ \ \ \ \ \ \ \ 1\ \ \ \ \ \ \ 42\leavevmode\hss\endgraf
o13\ :\ Matrix\ R\ \ <---\ R\leavevmode\hss\endgraf
\endgroup
\penalty-1000
\par
\vskip 1 pt
\noindent

Returning to the general case, an element $b $ of $\N A$ is called a
{\it Graver degree} if there exists a binomial $x^u - x^v$ in the
Graver basis $Gr_A$ such that $Au = Av = b$. If $b$ is a Graver degree
then the set of monomials in $R_b$ is the corresponding {\it Graver
  fiber}.  In our running example there are $37$ Graver fibers each
corresponding to a different Graver degree. To compute the Graver
fibers of $A$, we first construct the {\tt DoubleIdeal} $DI := 
\langle x^ax^b : x^a-x^b  \in Gr_A \rangle$ which is contained in
every monomial ideal of $Hilb_A$. Since our purpose in constructing
Graver fibers is to choose standard monomials from each fiber, it
suffices to list the monomials in each Graver degree that do not lie in
$DI$. Since $R$ is multigraded by $A$, we can obtain such a 
presentation of a Graver fiber by simply asking for the basis of $R$ 
in degree $b$ modulo $DI$. 
\par
\vskip 5 pt
\begingroup
\tteight
\baselineskip=7.4pt
\lineskip=0pt
\obeyspaces
i14\ :\ graverFibers\ =\ (Graver)\ ->\ (\leavevmode\hss\endgraf
\ \ \ \ \ \ \ \ \ \ \ DoubleIdeal\ :=\ (I)\ ->\ (\ trim\ ideal(\leavevmode\hss\endgraf
\ \ \ \ \ \ \ \ \ \ \ \ \ apply(numgens\ I,\ a\ ->\ (\leavevmode\hss\endgraf
\ \ \ \ \ \ \ \ \ \ \ \ \ \ \ \ \ f\ :=\ I{\char`\_}a;\ leadTerm\ f\ *\ (f\ -\ leadTerm\ f)))));\leavevmode\hss\endgraf
\ \ \ \ \ \ \ \ \ \ \ DI\ :=\ DoubleIdeal\ ideal\ Graver;\leavevmode\hss\endgraf
\ \ \ \ \ \ \ \ \ \ \ R\ :=\ ring\ Graver;\leavevmode\hss\endgraf
\ \ \ \ \ \ \ \ \ \ \ result\ :=\ new\ MutableHashTable;\leavevmode\hss\endgraf
\ \ \ \ \ \ \ \ \ \ \ scan(degrees\ source\ Graver,\ d\ ->\ (\leavevmode\hss\endgraf
\ \ \ \ \ \ \ \ \ \ \ \ \ \ \ \ if\ not\ result{\char`\#}?d\ then\ \leavevmode\hss\endgraf
\ \ \ \ \ \ \ \ \ \ \ \ \ \ \ \ \ \ \ \ \ result{\char`\#}d\ =\ compress\ (basis(d,R)\ {\char`\%}\ DI)));\leavevmode\hss\endgraf
\ \ \ \ \ \ \ \ \ \ \ result);\ \leavevmode\hss\endgraf
\endgroup
\penalty-1000
\par
\vskip 1 pt
\noindent
\par
\vskip 5 pt
\begingroup
\tteight
\baselineskip=7.4pt
\lineskip=0pt
\obeyspaces
i15\ :\ \ fibers\ =\ graverFibers\ Graver;\ \leavevmode\hss\endgraf
\endgroup
\penalty-1000
\par
\vskip 1 pt
\noindent
\par
\vskip 5 pt
\begingroup
\tteight
\baselineskip=7.4pt
\lineskip=0pt
\obeyspaces
i16\ :\ \ peek\ fibers\ \leavevmode\hss\endgraf
\penalty-500\leavevmode\hss\endgraf
o16\ =\ MutableHashTable{\char`\{}{\char`\{}2,\ 2{\char`\}}\ =>\ {\char`\{}0,\ 0{\char`\}}\ |\ ac\ b2\ |\ \ \ \ \ \ \ \ \ \ \ \ \ \ \ \ \ \ \ \ \  $\cdot\cdot\cdot$\leavevmode\hss\endgraf
\ \ \ \ \ \ \ \ \ \ \ \ \ \ \ \ \ \ \ \ \ \ \ {\char`\{}2,\ 8{\char`\}}\ =>\ {\char`\{}0,\ 0{\char`\}}\ |\ ae\ bd\ |\leavevmode\hss\endgraf
\ \ \ \ \ \ \ \ \ \ \ \ \ \ \ \ \ \ \ \ \ \ \ {\char`\{}2,\ 9{\char`\}}\ =>\ {\char`\{}0,\ 0{\char`\}}\ |\ be\ cd\ |\leavevmode\hss\endgraf
\ \ \ \ \ \ \ \ \ \ \ \ \ \ \ \ \ \ \ \ \ \ \ {\char`\{}3,\ 16{\char`\}}\ =>\ {\char`\{}0,\ 0{\char`\}}\ |\ ae2\ bde\ cd2\ |\leavevmode\hss\endgraf
\ \ \ \ \ \ \ \ \ \ \ \ \ \ \ \ \ \ \ \ \ \ \ {\char`\{}4,\ 14{\char`\}}\ =>\ {\char`\{}0,\ 0{\char`\}}\ |\ a2d2\ c3e\ |\leavevmode\hss\endgraf
\ \ \ \ \ \ \ \ \ \ \ \ \ \ \ \ \ \ \ \ \ \ \ {\char`\{}4,\ 7{\char`\}}\ =>\ {\char`\{}0,\ 0{\char`\}}\ |\ a3d\ bc3\ |\leavevmode\hss\endgraf
\ \ \ \ \ \ \ \ \ \ \ \ \ \ \ \ \ \ \ \ \ \ \ {\char`\{}4,\ 8{\char`\}}\ =>\ {\char`\{}0,\ 0{\char`\}}\ |\ a3e\ a2bd\ c4\ |\leavevmode\hss\endgraf
\ \ \ \ \ \ \ \ \ \ \ \ \ \ \ \ \ \ \ \ \ \ \ {\char`\{}5,\ 10{\char`\}}\ =>\ {\char`\{}0,\ 0{\char`\}}\ |\ a3ce\ a2b2e\ a2bcd\ ab3d\ c5\ |\leavevmode\hss\endgraf
\ \ \ \ \ \ \ \ \ \ \ \ \ \ \ \ \ \ \ \ \ \ \ {\char`\{}5,\ 14{\char`\}}\ =>\ {\char`\{}0,\ 0{\char`\}}\ |\ a3d2\ ac3e\ b2c2e\ bc3d\ |\leavevmode\hss\endgraf
\ \ \ \ \ \ \ \ \ \ \ \ \ \ \ \ \ \ \ \ \ \ \ {\char`\{}5,\ 16{\char`\}}\ =>\ {\char`\{}0,\ 0{\char`\}}\ |\ a3e2\ a2cd2\ ab2d2\ c4e\ |\leavevmode\hss\endgraf
\ \ \ \ \ \ \ \ \ \ \ \ \ \ \ \ \ \ \ \ \ \ \ {\char`\{}5,\ 21{\char`\}}\ =>\ {\char`\{}0,\ 0{\char`\}}\ |\ a2d3\ bc2e2\ c3de\ |\leavevmode\hss\endgraf
\ \ \ \ \ \ \ \ \ \ \ \ \ \ \ \ \ \ \ \ \ \ \ {\char`\{}5,\ 22{\char`\}}\ =>\ {\char`\{}0,\ 0{\char`\}}\ |\ a2d2e\ abd3\ c3e2\ |\leavevmode\hss\endgraf
\ \ \ \ \ \ \ \ \ \ \ \ \ \ \ \ \ \ \ \ \ \ \ {\char`\{}5,\ 28{\char`\}}\ =>\ {\char`\{}0,\ 0{\char`\}}\ |\ ad4\ c2e3\ |\leavevmode\hss\endgraf
\ \ \ \ \ \ \ \ \ \ \ \ \ \ \ \ \ \ \ \ \ \ \ {\char`\{}5,\ 7{\char`\}}\ =>\ {\char`\{}0,\ 0{\char`\}}\ |\ a4d\ abc3\ b3c2\ |\leavevmode\hss\endgraf
\ \ \ \ \ \ \ \ \ \ \ \ \ \ \ \ \ \ \ \ \ \ \ {\char`\{}5,\ 8{\char`\}}\ =>\ {\char`\{}0,\ 0{\char`\}}\ |\ a4e\ a3bd\ ac4\ b2c3\ |\leavevmode\hss\endgraf
\ \ \ \ \ \ \ \ \ \ \ \ \ \ \ \ \ \ \ \ \ \ \ {\char`\{}6,\ 12{\char`\}}\ =>\ {\char`\{}0,\ 0{\char`\}}\ |\ a3c2e\ a2bc2d\ ab4e\ b5d\ c6\ |\leavevmode\hss\endgraf
\ \ \ \ \ \ \ \ \ \ \ \ \ \ \ \ \ \ \ \ \ \ \ {\char`\{}6,\ 14{\char`\}}\ =>\ {\char`\{}0,\ 0{\char`\}}\ |\ a4d2\ a2c3e\ abc3d\ b4ce\ b3c2d $\cdot\cdot\cdot$\leavevmode\hss\endgraf
\ \ \ \ \ \ \ \ \ \ \ \ \ \ \ \ \ \ \ \ \ \ \ {\char`\{}6,\ 18{\char`\}}\ =>\ {\char`\{}0,\ 0{\char`\}}\ |\ a3ce2\ a2b2e2\ a2c2d2\ b4d2\ c5 $\cdot\cdot\cdot$\leavevmode\hss\endgraf
\ \ \ \ \ \ \ \ \ \ \ \ \ \ \ \ \ \ \ \ \ \ \ {\char`\{}6,\ 21{\char`\}}\ =>\ {\char`\{}0,\ 0{\char`\}}\ |\ a3d3\ abc2e2\ ac3de\ b3ce2\ bc3 $\cdot\cdot\cdot$\leavevmode\hss\endgraf
\ \ \ \ \ \ \ \ \ \ \ \ \ \ \ \ \ \ \ \ \ \ \ {\char`\{}6,\ 24{\char`\}}\ =>\ {\char`\{}0,\ 0{\char`\}}\ |\ a3e3\ a2cd2e\ abcd3\ b3d3\ c4e2 $\cdot\cdot\cdot$\leavevmode\hss\endgraf
\ \ \ \ \ \ \ \ \ \ \ \ \ \ \ \ \ \ \ \ \ \ \ {\char`\{}6,\ 28{\char`\}}\ =>\ {\char`\{}0,\ 0{\char`\}}\ |\ a2d4\ ac2e3\ b2ce3\ c3d2e\ |\leavevmode\hss\endgraf
\ \ \ \ \ \ \ \ \ \ \ \ \ \ \ \ \ \ \ \ \ \ \ {\char`\{}6,\ 30{\char`\}}\ =>\ {\char`\{}0,\ 0{\char`\}}\ |\ a2d2e2\ acd4\ b2d4\ c3e3\ |\leavevmode\hss\endgraf
\ \ \ \ \ \ \ \ \ \ \ \ \ \ \ \ \ \ \ \ \ \ \ {\char`\{}6,\ 35{\char`\}}\ =>\ {\char`\{}0,\ 0{\char`\}}\ |\ ad5\ bce4\ c2de3\ |\leavevmode\hss\endgraf
\ \ \ \ \ \ \ \ \ \ \ \ \ \ \ \ \ \ \ \ \ \ \ {\char`\{}6,\ 36{\char`\}}\ =>\ {\char`\{}0,\ 0{\char`\}}\ |\ ad4e\ bd5\ c2e4\ |\leavevmode\hss\endgraf
\ \ \ \ \ \ \ \ \ \ \ \ \ \ \ \ \ \ \ \ \ \ \ {\char`\{}6,\ 42{\char`\}}\ =>\ {\char`\{}0,\ 0{\char`\}}\ |\ ce5\ d6\ |\leavevmode\hss\endgraf
\ \ \ \ \ \ \ \ \ \ \ \ \ \ \ \ \ \ \ \ \ \ \ {\char`\{}6,\ 7{\char`\}}\ =>\ {\char`\{}0,\ 0{\char`\}}\ |\ a5d\ a2bc3\ b5c\ |\leavevmode\hss\endgraf
\ \ \ \ \ \ \ \ \ \ \ \ \ \ \ \ \ \ \ \ \ \ \ {\char`\{}6,\ 8{\char`\}}\ =>\ {\char`\{}0,\ 0{\char`\}}\ |\ a5e\ a4bd\ a2c4\ b4c2\ |\leavevmode\hss\endgraf
\ \ \ \ \ \ \ \ \ \ \ \ \ \ \ \ \ \ \ \ \ \ \ {\char`\{}7,\ 14{\char`\}}\ =>\ {\char`\{}0,\ 0{\char`\}}\ |\ a5d2\ a3c3e\ a2bc3d\ b6e\ b5cd\  $\cdot\cdot\cdot$\leavevmode\hss\endgraf
\ \ \ \ \ \ \ \ \ \ \ \ \ \ \ \ \ \ \ \ \ \ \ {\char`\{}7,\ 21{\char`\}}\ =>\ {\char`\{}0,\ 0{\char`\}}\ |\ a4d3\ a2bc2e2\ a2c3de\ abc3d2\  $\cdot\cdot\cdot$\leavevmode\hss\endgraf
\ \ \ \ \ \ \ \ \ \ \ \ \ \ \ \ \ \ \ \ \ \ \ {\char`\{}7,\ 28{\char`\}}\ =>\ {\char`\{}0,\ 0{\char`\}}\ |\ a3d4\ a2c2e3\ ac3d2e\ b4e3\ bc3 $\cdot\cdot\cdot$\leavevmode\hss\endgraf
\ \ \ \ \ \ \ \ \ \ \ \ \ \ \ \ \ \ \ \ \ \ \ {\char`\{}7,\ 35{\char`\}}\ =>\ {\char`\{}0,\ 0{\char`\}}\ |\ a2d5\ abce4\ ac2de3\ b3e4\ c3d3 $\cdot\cdot\cdot$\leavevmode\hss\endgraf
\ \ \ \ \ \ \ \ \ \ \ \ \ \ \ \ \ \ \ \ \ \ \ {\char`\{}7,\ 42{\char`\}}\ =>\ {\char`\{}0,\ 0{\char`\}}\ |\ ace5\ ad6\ b2e5\ c2d2e3\ |\leavevmode\hss\endgraf
\ \ \ \ \ \ \ \ \ \ \ \ \ \ \ \ \ \ \ \ \ \ \ {\char`\{}7,\ 49{\char`\}}\ =>\ {\char`\{}0,\ 0{\char`\}}\ |\ be6\ cde5\ d7\ |\leavevmode\hss\endgraf
\ \ \ \ \ \ \ \ \ \ \ \ \ \ \ \ \ \ \ \ \ \ \ {\char`\{}7,\ 7{\char`\}}\ =>\ {\char`\{}0,\ 0{\char`\}}\ |\ a6d\ a3bc3\ b7\ |\leavevmode\hss\endgraf
\ \ \ \ \ \ \ \ \ \ \ \ \ \ \ \ \ \ \ \ \ \ \ {\char`\{}7,\ 8{\char`\}}\ =>\ {\char`\{}0,\ 0{\char`\}}\ |\ a6e\ a5bd\ a3c4\ b6c\ |\leavevmode\hss\endgraf
\ \ \ \ \ \ \ \ \ \ \ \ \ \ \ \ \ \ \ \ \ \ \ {\char`\{}8,\ 56{\char`\}}\ =>\ {\char`\{}0,\ 0{\char`\}}\ |\ ae7\ bde6\ cd2e5\ d8\ |\leavevmode\hss\endgraf
\ \ \ \ \ \ \ \ \ \ \ \ \ \ \ \ \ \ \ \ \ \ \ {\char`\{}8,\ 8{\char`\}}\ =>\ {\char`\{}0,\ 0{\char`\}}\ |\ a7e\ a6bd\ a4c4\ b8\ |\leavevmode\hss\endgraf
\penalty-500\leavevmode\hss\endgraf
o16\ :\ Net\leavevmode\hss\endgraf
\endgroup
\penalty-1000
\par
\vskip 1 pt
\noindent

For example, the Graver degree $(8,8)$ corresponds to the Graver fiber
$$ \bigl\{\,
\underline{a^7 e}, \, \underline{a^6 b d},\,  \underline{a^4 c^4}, \,
a^3 b^2 c^3,\, a^2 b^4 c^2,\,  a b^6 c, \, \underline{b^8} \,\bigr\}.$$
Our \Mtwo code only outputs the four underlined monomials,
in the format {\tt  | a7e a6bd a4c4 b8 |}. The three non-underlined 
monomials lie in the {\tt DoubleIdeal}. Graver degrees are
important because of the following result.

\begin{lemma} \cite[Lemma 10.5]{St2}.
The multidegree of any minimal generator of any ideal 
$I$ in $Hilb_A$ is a Graver degree.
\end{lemma}

The next step in constructing the toric Hilbert scheme is to compute
all its fixed points with respect to the scaling action of the
$n$-dimensional algebraic torus $(\k^*)^n$.  These fixed
points are the monomial ideals $M$ lying on $Hilb_A$.  Every term
order $\prec$ on the polynomial ring $R$ gives such a monomial ideal:
$M = in_\prec(I_A )$, the initial ideal of the toric ideal $I_A$ with 
respect to $\prec$. We call
these the {\it coherent} monomial ideals on $Hilb_A$. They can be
computed by \cite[Algorithm 3.6]{St2} applied to
$I_A$. A refinement and fast implementation can be found in the 
software package {\tt TiGERS} by Huber and Thomas \cite{HT}. 

Now we wish to compute all monomial ideals $M$ on $Hilb_A$ regardless
of whether $M$ is coherent or not. For this we use the procedure
{\tt generateAmonos} given below. This procedure takes in the Graver
basis $Gr_A$ and records the numerator of the Hilbert series of $I_A$
in {\tt trueHS}. It then computes the Graver fibers of $A$, sorts them
and calls the subroutine {\tt selectStandard} to generate a
candidate for a monomial ideal on $Hilb_A$.\\

\par
\vskip 5 pt
\begingroup
\tteight
\baselineskip=7.4pt
\lineskip=0pt
\obeyspaces
i17\ :\ generateAmonos\ =\ (Graver)\ ->\ (\leavevmode\hss\endgraf
\ \ \ \ \ \ \ \ \ \ \ trueHS\ =\ poincare\ coker\ Graver;\leavevmode\hss\endgraf
\ \ \ \ \ \ \ \ \ \ \ fibers\ =\ graverFibers\ Graver;\leavevmode\hss\endgraf
\ \ \ \ \ \ \ \ \ \ \ fibers\ =\ apply(sort\ pairs\ fibers,\ last);\leavevmode\hss\endgraf
\ \ \ \ \ \ \ \ \ \ \ monos\ =\ {\char`\{}{\char`\}};\leavevmode\hss\endgraf
\ \ \ \ \ \ \ \ \ \ \ selectStandard(fibers,\ ideal(0{\char`\_}(ring\ Graver)));\leavevmode\hss\endgraf
\ \ \ \ \ \ \ \ \ \ \ )\ ;\ \leavevmode\hss\endgraf
\endgroup
\penalty-1000
\par
\vskip 1 pt
\noindent

The input to the subroutine {\tt selectStandard}
are the Graver fibers given as a list of matrices and a monomial 
ideal $J$ which should be included in every $A$-graded ideal 
that we generate. The subroutine then loops through each Graver fiber, 
and at each step selects a standard monomial from that fiber and 
updates the ideal $J$ by adding the other monomials in this fiber 
to $J$. The final $J$ output by the subroutine is the candidate ideal
that is sent back to {\tt generateAmonos}. It is stored by the program 
if its Hilbert series agrees with that of $I_A$. 

\par
\vskip 5 pt
\begingroup
\tteight
\baselineskip=7.4pt
\lineskip=0pt
\obeyspaces
i18\ :\ \ selectStandard\ =\ (fibers,\ J)\ ->\ (\leavevmode\hss\endgraf
\ \ \ \ \ \ \ \ \ \ \ if\ {\char`\#}fibers\ ==\ 0\ then\ (\leavevmode\hss\endgraf
\ \ \ \ \ \ \ \ \ \ \ \ \ \ if\ trueHS\ ==\ poincare\ coker\ gens\ J\leavevmode\hss\endgraf
\ \ \ \ \ \ \ \ \ \ \ \ \ \ then\ (monos\ =\ append(monos,flatten\ entries\ mingens\ J));\leavevmode\hss\endgraf
\ \ \ \ \ \ \ \ \ \ \ )\ else\ (\leavevmode\hss\endgraf
\ \ \ \ \ \ \ \ \ \ \ \ \ \ P\ :=\ fibers{\char`\_}0;\leavevmode\hss\endgraf
\ \ \ \ \ \ \ \ \ \ \ \ \ \ fibers\ =\ drop(fibers,1);\leavevmode\hss\endgraf
\ \ \ \ \ \ \ \ \ \ \ \ \ \ P\ =\ compress(P\ {\char`\%}\ J);\leavevmode\hss\endgraf
\ \ \ \ \ \ \ \ \ \ \ \ \ \ nP\ :=\ numgens\ source\ P;\ \leavevmode\hss\endgraf
\ \ \ \ \ \ \ \ \ \ \ \ \ \ --\ nP\ is\ the\ number\ of\ monomials\ not\ in\ J.\leavevmode\hss\endgraf
\ \ \ \ \ \ \ \ \ \ \ \ \ \ if\ nP\ >\ 0\ then\ (\leavevmode\hss\endgraf
\ \ \ \ \ \ \ \ \ \ \ \ \ \ \ \ \ if\ nP\ ==\ 1\ then\ selectStandard(fibers,J)\leavevmode\hss\endgraf
\ \ \ \ \ \ \ \ \ \ \ \ \ \ \ \ \ else\ (--remove\ one\ monomial\ from\ P,take\ the\ rest.\leavevmode\hss\endgraf
\ \ \ \ \ \ \ \ \ \ \ \ \ \ \ \ \ \ \ \ \ \ \ P\ =\ flatten\ entries\ P;\leavevmode\hss\endgraf
\ \ \ \ \ \ \ \ \ \ \ \ \ \ \ \ \ \ \ \ \ \ \ scan({\char`\#}P,\ i\ ->\ (\leavevmode\hss\endgraf
\ \ \ \ \ \ \ \ \ \ \ \ \ \ \ \ \ \ \ \ \ \ \ \ \ \ \ \ J1\ :=\ J\ +\ ideal\ drop(P,{\char`\{}i,i{\char`\}});\leavevmode\hss\endgraf
\ \ \ \ \ \ \ \ \ \ \ \ \ \ \ \ \ \ \ \ \ \ \ \ \ \ \ \ selectStandard(fibers,\ J1)))));\leavevmode\hss\endgraf
\ \ \ \ \ \ \ \ \ \ \ ));\ \leavevmode\hss\endgraf
\endgroup
\penalty-1000
\par
\vskip 1 pt
\noindent

All the monomial $A$-graded ideals are stored in the list {\tt monos}.
Below, we ask \Mtwo for the cardinality of {\tt monos} and its 
first ten elements.

\par
\vskip 5 pt
\begingroup
\tteight
\baselineskip=7.4pt
\lineskip=0pt
\obeyspaces
i19\ :\ generateAmonos(Graver);\leavevmode\hss\endgraf
\endgroup
\penalty-1000
\par
\vskip 1 pt
\noindent
\par
\vskip 5 pt
\begingroup
\tteight
\baselineskip=7.4pt
\lineskip=0pt
\obeyspaces
i20\ :\ {\char`\#}monos\ \leavevmode\hss\endgraf
\penalty-500\leavevmode\hss\endgraf
o20\ =\ 281\leavevmode\hss\endgraf
\endgroup
\penalty-1000
\par
\vskip 1 pt
\noindent
\par
\vskip 5 pt
\begingroup
\tteight
\baselineskip=7.4pt
\lineskip=0pt
\obeyspaces
i21\ :\ scan(0..9,\ i\ ->\ print\ toString\ monos{\char`\#}i)\leavevmode\hss\endgraf
{\char`\{}c*d,\ b*d,\ b{\char`\^}2,\ c{\char`\^}3*e,\ c{\char`\^}4,\ b*c{\char`\^}3,\ c{\char`\^}2*e{\char`\^}3,\ b*c{\char`\^}2*e{\char`\^}2,\ b*c*e{\char`\^}4,\ d{\char`\^}6{\char`\}}\leavevmode\hss\endgraf
{\char`\{}c*d,\ b*d,\ b{\char`\^}2,\ c{\char`\^}3*e,\ c{\char`\^}4,\ b*c{\char`\^}3,\ c{\char`\^}2*e{\char`\^}3,\ b*c{\char`\^}2*e{\char`\^}2,\ c*e{\char`\^}5,\ b*c*e{\char`\^}4, $\cdot\cdot\cdot$\leavevmode\hss\endgraf
{\char`\{}c*d,\ b*d,\ b{\char`\^}2,\ c{\char`\^}3*e,\ c{\char`\^}4,\ b*c{\char`\^}3,\ c{\char`\^}2*e{\char`\^}3,\ b*c{\char`\^}2*e{\char`\^}2,\ c*e{\char`\^}5,\ b*c*e{\char`\^}4, $\cdot\cdot\cdot$\leavevmode\hss\endgraf
{\char`\{}c*d,\ b*d,\ b{\char`\^}2,\ c{\char`\^}3*e,\ c{\char`\^}4,\ b*c{\char`\^}3,\ c{\char`\^}2*e{\char`\^}3,\ b*c{\char`\^}2*e{\char`\^}2,\ c*e{\char`\^}5,\ b*c*e{\char`\^}4, $\cdot\cdot\cdot$\leavevmode\hss\endgraf
{\char`\{}c*d,\ b*d,\ b{\char`\^}2,\ c{\char`\^}3*e,\ c{\char`\^}4,\ b*c{\char`\^}3,\ c{\char`\^}2*e{\char`\^}3,\ b*c{\char`\^}2*e{\char`\^}2,\ d{\char`\^}6,\ a*d{\char`\^}5{\char`\}}\leavevmode\hss\endgraf
{\char`\{}c*d,\ b*d,\ b{\char`\^}2,\ c{\char`\^}3*e,\ c{\char`\^}4,\ b*c{\char`\^}3,\ b*c{\char`\^}2*e{\char`\^}2,\ a*d{\char`\^}4,\ d{\char`\^}6{\char`\}}\leavevmode\hss\endgraf
{\char`\{}c*d,\ b*d,\ b{\char`\^}2,\ c{\char`\^}3*e,\ c{\char`\^}4,\ b*c{\char`\^}3,\ a*d{\char`\^}4,\ a{\char`\^}2*d{\char`\^}3,\ d{\char`\^}6{\char`\}}\leavevmode\hss\endgraf
{\char`\{}c*d,\ b*d,\ b{\char`\^}2,\ a{\char`\^}2*d{\char`\^}2,\ c{\char`\^}4,\ b*c{\char`\^}3,\ a*d{\char`\^}4,\ d{\char`\^}6{\char`\}}\leavevmode\hss\endgraf
{\char`\{}c*d,\ b*d,\ b{\char`\^}2,\ a{\char`\^}2*d{\char`\^}2,\ a{\char`\^}3*d,\ c{\char`\^}4,\ a*d{\char`\^}4,\ d{\char`\^}6{\char`\}}\leavevmode\hss\endgraf
{\char`\{}c*d,\ b*d,\ b{\char`\^}2,\ a{\char`\^}3*e,\ a{\char`\^}2*d{\char`\^}2,\ a{\char`\^}3*d,\ a*d{\char`\^}4,\ d{\char`\^}6{\char`\}}\leavevmode\hss\endgraf
\endgroup
\penalty-1000
\par
\vskip 1 pt
\noindent

The monomial ideals (torus-fixed points) on $Hilb_A$ form the vertices
of the {\it flip graph} of $A$ whose edges correspond to the
torus-fixed curves on $Hilb_A$. This graph was introduced in \cite{MT}
and provides structural information about $Hilb_A$.  The edges
emanating from a monomial ideal $M$ can be constructed as follows.
For any minimal generator $x^u$ of $M$, let
$x^v$ be the unique monomial with $x^v \not\in M$ and $Au = Av$. Form
the {\it wall ideal} which is generated by $x^u - x^v$ and all minimal
generators of $M$ other than $x^u$, and let $M'$ be the unique initial
monomial ideal of the wall ideal which contains $x^v$. If $M'$ lies on
$Hilb_A$ then $\{M, M'\}$ is an edge of the flip graph. We now
illustrate the \Mtwo procedure for computing all flip neighbors of a
monomial $A$-graded ideal.
 
\par
\vskip 5 pt
\begingroup
\tteight
\baselineskip=7.4pt
\lineskip=0pt
\obeyspaces
i22\ :\ findPositiveVector\ =\ (m,s)\ ->\ (\leavevmode\hss\endgraf
\ \ \ \ \ \ \ \ \ \ \ expvector\ :=\ first\ exponents\ s\ -\ first\ exponents\ m;\leavevmode\hss\endgraf
\ \ \ \ \ \ \ \ \ \ \ n\ :=\ {\char`\#}expvector;\leavevmode\hss\endgraf
\ \ \ \ \ \ \ \ \ \ \ i\ :=\ first\ positions(0..n-1,\ j\ ->\ expvector{\char`\_}j\ >\ 0);\leavevmode\hss\endgraf
\ \ \ \ \ \ \ \ \ \ \ splice\ {\char`\{}i:0,\ 1,\ (n-i-1):0{\char`\}}\leavevmode\hss\endgraf
\ \ \ \ \ \ \ \ \ \ \ );\leavevmode\hss\endgraf
\endgroup
\penalty-1000
\par
\vskip 1 pt
\noindent

\par
\vskip 5 pt
\begingroup
\tteight
\baselineskip=7.4pt
\lineskip=0pt
\obeyspaces
i23\ :\ flips\ =\ (M)\ ->\ (\leavevmode\hss\endgraf
\ \ \ \ \ \ \ \ \ \ \ R\ :=\ ring\ M;\leavevmode\hss\endgraf
\ \ \ \ \ \ \ \ \ \ \ --\ store\ generators\ of\ M\ in\ monoms\leavevmode\hss\endgraf
\ \ \ \ \ \ \ \ \ \ \ monoms\ :=\ first\ entries\ generators\ M;\leavevmode\hss\endgraf
\ \ \ \ \ \ \ \ \ \ \ result\ :=\ {\char`\{}{\char`\}};\leavevmode\hss\endgraf
\ \ \ \ \ \ \ \ \ \ \ --\ test\ each\ generator\ of\ M\ to\ see\ if\ it\ leads\ to\ a\ neighbo $\cdot\cdot\cdot$\leavevmode\hss\endgraf
\ \ \ \ \ \ \ \ \ \ \ scan({\char`\#}monoms,\ i\ ->\ (\leavevmode\hss\endgraf
\ \ \ \ \ \ \ \ \ \ \ \ \ m\ :=\ monoms{\char`\_}i;\leavevmode\hss\endgraf
\ \ \ \ \ \ \ \ \ \ \ \ \ rest\ :=\ drop(monoms,{\char`\{}i,i{\char`\}});\leavevmode\hss\endgraf
\ \ \ \ \ \ \ \ \ \ \ \ \ b\ :=\ basis(degree\ m,\ R);\leavevmode\hss\endgraf
\ \ \ \ \ \ \ \ \ \ \ \ \ s\ :=\ (compress\ (b\ {\char`\%}\ M)){\char`\_}(0,0);\leavevmode\hss\endgraf
\ \ \ \ \ \ \ \ \ \ \ \ \ J\ :=\ ideal(m-s)\ +\ ideal\ rest;\leavevmode\hss\endgraf
\ \ \ \ \ \ \ \ \ \ \ \ \ if\ poincare\ coker\ gens\ J\ ==\ poincare\ coker\ gens\ M\ then\ (\leavevmode\hss\endgraf
\ \ \ \ \ \ \ \ \ \ \ \ \ \ \ w\ :=\ findPositiveVector(m,s);\leavevmode\hss\endgraf
\ \ \ \ \ \ \ \ \ \ \ \ \ \ \ R1\ :=\ (coefficientRing\ R)[generators\ R,\ Weights=>w];\leavevmode\hss\endgraf
\ \ \ \ \ \ \ \ \ \ \ \ \ \ \ J\ =\ substitute(J,R1);\leavevmode\hss\endgraf
\ \ \ \ \ \ \ \ \ \ \ \ \ \ \ J\ =\ trim\ ideal\ leadTerm\ J;\leavevmode\hss\endgraf
\ \ \ \ \ \ \ \ \ \ \ \ \ \ \ result\ =\ append(result,J);\leavevmode\hss\endgraf
\ \ \ \ \ \ \ \ \ \ \ \ \ \ \ \ )));\leavevmode\hss\endgraf
\ \ \ \ \ \ \ \ \ \ \ scan({\char`\#}result,\ i->(print\ gens\ result{\char`\_}i)));\leavevmode\hss\endgraf
\endgroup
\penalty-1000
\par
\vskip 1 pt
\noindent

The above code inputs a monomial $A$-graded ideal $M$ whose 
minimal generators are stored in the list {\tt monoms}. The flip 
neighbors of $M$ will be stored in {\tt result}. For each monomial 
$x^u$ in {\tt monoms} we need to test whether it yields a flip 
neighbor of $M$ or not. At the $i$-th step of this loop, we let {\tt m}
be the $i$-th monomial in {\tt monoms}. The list {\tt rest} contains 
all monomials in {\tt monoms} except {\tt m}. We compute 
the standard monomial {\tt s} of $M$ of the same degree as $m$.
The {\em wall ideal} of $m-s$ is the binomial ideal 
$J$ generated by $m-s$ and the monomials in {\tt rest}. We then check 
whether $J$ is $A$-graded by comparing its Hilbert series with that 
of $M$. If this is the case, we use the subroutine {\tt
findPositiveVector} to find a unit vector $w = (0,\ldots,1,\ldots,0)$
such that $w \cdot s > w \cdot m$. The flip neighbor is then the 
initial ideal of $J$ with respect to $w$ and it is stored in 
{\tt result}. The program outputs the minimal generators of
each flip neighbor. Here is an example.
 
\par
\vskip 5 pt
\begingroup
\tteight
\baselineskip=7.4pt
\lineskip=0pt
\obeyspaces
i24\ :\ R\ =\ QQ[a..e,Degrees=>transpose\ A];\leavevmode\hss\endgraf
\endgroup
\penalty-1000
\par
\vskip 1 pt
\noindent
\par
\vskip 5 pt
\begingroup
\tteight
\baselineskip=7.4pt
\lineskip=0pt
\obeyspaces
i25\ :\ M\ =\ ideal(a*e,c*d,a*c,a{\char`\^}2*d{\char`\^}2,a{\char`\^}2*b*d,a{\char`\^}3*d,c{\char`\^}2*e{\char`\^}3,\leavevmode\hss\endgraf
\ \ \ \ \ \ \ \ \ \ \ \ \ \ \ \ c{\char`\^}3*e{\char`\^}2,c{\char`\^}4*e,c{\char`\^}5,c*e{\char`\^}5,a*d{\char`\^}5,b*e{\char`\^}6);\leavevmode\hss\endgraf
\penalty-500\leavevmode\hss\endgraf
o25\ :\ Ideal\ of\ R\leavevmode\hss\endgraf
\endgroup
\penalty-1000
\par
\vskip 1 pt
\noindent
\par
\vskip 5 pt
\begingroup
\tteight
\baselineskip=7.4pt
\lineskip=0pt
\obeyspaces
i26\ :\ flips\ M\ \leavevmode\hss\endgraf
{\char`\{}0{\char`\}}\ |\ ae\ cd\ ac\ a2d2\ a3d\ c4\ c2e3\ c3e2\ ad5\ ce5\ be6\ |\leavevmode\hss\endgraf
{\char`\{}0{\char`\}}\ |\ cd\ ae\ ac\ a2d2\ a2bd\ a3d\ c3e2\ c4e\ c5\ ad4\ ce5\ c2e4\ be6\ |\leavevmode\hss\endgraf
{\char`\{}0{\char`\}}\ |\ ae\ cd\ ac\ a2d2\ a3d\ a2bd\ c2e3\ c3e2\ c4e\ c5\ ce5\ bce4\ ad6\ be6\ |\leavevmode\hss\endgraf
{\char`\{}0{\char`\}}\ |\ ae\ ac\ cd\ a2bd\ a3d\ a2d2\ c2e3\ c3e2\ c4e\ c5\ ce5\ ad5\ d7\ |\leavevmode\hss\endgraf
\endgroup
\penalty-1000
\par
\vskip 1 pt
\noindent

It is an open problem whether the toric Hilbert scheme $Hilb_A$ is
connected. Recent work in geometric combinatorics \cite{San} suggests
that this is probably false for some $A$. This result and its 
implications for $Hilb_A$ will be discussed further in Section 2.
However, we have the following result of Maclagan and Thomas
\cite{MT}.

\begin{theorem} 
The flip graph of $A$ is connected if and only if the toric Hilbert 
scheme $Hilb_A$ is connected.
\end{theorem}

We now have two algorithms for listing monomial ideals on $Hilb_A$.
First, there is the {\it backtracking algorithm} whose \Mtwo
implementation was described above.  Second, there is the {\it flip
  search algorithm} which starts with any coherent monomial ideal $M$
and then constructs the connected component of $M$ in the flip graph
of $A$ by carrying out local flips as above.  This procedure is also 
implemented in {\tt TiGERS} \cite{HT}. Clearly, the two algorithms
will produce the same answer if and only if $Hilb_A$ is connected. In
other words, finding an example where $Hilb_A$ is disconnected is
equivalent to finding a matrix $A$ for which the flip search algorithm
produces fewer monomial ideals than the backtracking algorithm.

\section{Polyhedral Geometry}

Algorithms from polyhedral geometry are essential in the study of the
toric Hilbert scheme. Consider the problem of deciding whether or not
a given monomial ideal $M$ in $Hilb_A$ is coherent.  This problem
gives rise to a system of linear inequalities as follows. Let
$x^{u_1}, \ldots, x^{u_r}$ be the minimal generators of $M$, and let
$x^{v_i}$ be the unique standard monomial with $A u_i = A v_i$. Then
$M$ is coherent if and only if there exists a vector $w \in \R^n$ such
that $\,w \cdot (u_i - v_i) > 0\,$ for $i =1,\ldots,r$.  Thus the test
for coherence amounts to solving a {\sl feasibility problem of linear
  programming}, and there are many highly efficient implementations
(based on the simplex algorithms or interior point methods) available
for this task. For our experimental purposes, it is convenient to use
the code {\tt polarCone.m2}, given in Appendix 1, which is based on
the (inefficient but easy-to-implement) {\em Fourier-Motzkin
  elimination} method \cite{Zie}.  This code converts the generator
representation of a polyhedron to its inequality representation and
vice versa. A simple example is given in Appendix 1. In
particular, given a Gr\"obner basis $\mathcal G$ of $I_A$, the
function {\tt polarCone} will compute all the extreme rays of the {\em
  Gr\"obner cone} $\,\{ w \in \R^n \,: \,w \cdot (u_i - v_i) \geq 0\,$
for each $x^{u_i}-x^{v_i} \in {\mathcal G}\}.$

We now show how to use \Mtwo to decide whether a 
monomial $A$-graded ideal $M$ is coherent. The first step in 
this calculation is to compute all the standard monomials of $M$ 
of the same degree as the minimal generators of $M$. We do this 
using the procedure {\tt stdMonomials}.

\par
\vskip 5 pt
\begingroup
\tteight
\baselineskip=7.4pt
\lineskip=0pt
\obeyspaces
i27\ :\ stdMonomials\ =\ (M)\ ->\ (\leavevmode\hss\endgraf
\ \ \ \ \ \ \ \ \ \ \ R\ :=\ ring\ M;\leavevmode\hss\endgraf
\ \ \ \ \ \ \ \ \ \ \ RM\ :=\ R/M;\leavevmode\hss\endgraf
\ \ \ \ \ \ \ \ \ \ \ apply(numgens\ M,\ i\ ->\ (\leavevmode\hss\endgraf
\ \ \ \ \ \ \ \ \ \ \ \ \ \ \ \ \ s\ :=\ basis(degree(M{\char`\_}i),RM);\ lift(s{\char`\_}(0,0),\ R)))\leavevmode\hss\endgraf
\ \ \ \ \ \ \ \ \ \ \ );\ \leavevmode\hss\endgraf
\endgroup
\penalty-1000
\par
\vskip 1 pt
\noindent

As an example, consider the following monomial $A$-graded ideal.

\par
\vskip 5 pt
\begingroup
\tteight
\baselineskip=7.4pt
\lineskip=0pt
\obeyspaces
i28\ :\ \ R\ =\ QQ[a..e,Degrees\ =>\ transpose\ A\ ];\ \leavevmode\hss\endgraf
\endgroup
\penalty-1000
\par
\vskip 1 pt
\noindent
\par
\vskip 5 pt
\begingroup
\tteight
\baselineskip=7.4pt
\lineskip=0pt
\obeyspaces
i29\ :\ M\ =\ ideal(a{\char`\^}3*d,\ a{\char`\^}2*b*d,\ a{\char`\^}2*d{\char`\^}2,\ a*b{\char`\^}3*d,\ a*b{\char`\^}2*d{\char`\^}2,\ a*b*d{\char`\^}3,\ \leavevmode\hss\endgraf
\ \ \ \ \ \ \ \ \ \ \ \ \ \ \ \ a*c,\ a*d{\char`\^}4,\ a*e,\ b{\char`\^}5*d,\ b{\char`\^}4*d{\char`\^}2,\ b{\char`\^}3*d{\char`\^}3,\ b{\char`\^}2*d{\char`\^}4,\ \leavevmode\hss\endgraf
\ \ \ \ \ \ \ \ \ \ \ \ \ \ \ \ b*d{\char`\^}5,\ b*e,\ c*e{\char`\^}5);\ \leavevmode\hss\endgraf
\penalty-500\leavevmode\hss\endgraf
o29\ :\ Ideal\ of\ R\leavevmode\hss\endgraf
\endgroup
\penalty-1000
\par
\vskip 1 pt
\noindent
\par
\vskip 5 pt
\begingroup
\tteight
\baselineskip=7.4pt
\lineskip=0pt
\obeyspaces
i30\ :\ \ toString\ stdMonomials\ M\ \leavevmode\hss\endgraf
\penalty-500\leavevmode\hss\endgraf
o30\ =\ {\char`\{}b*c{\char`\^}3,\ c{\char`\^}4,\ c{\char`\^}3*e,\ c{\char`\^}5,\ c{\char`\^}4*e,\ c{\char`\^}3*e{\char`\^}2,\ b{\char`\^}2,\ c{\char`\^}2*e{\char`\^}3,\ b*d,\ c{\char`\^}6, $\cdot\cdot\cdot$\leavevmode\hss\endgraf
\penalty-500\leavevmode\hss\endgraf
o30\ :\ String\leavevmode\hss\endgraf
\endgroup
\penalty-1000
\par
\vskip 1 pt
\noindent

From the pairs $(x^u,x^v)$ of minimal generators $x^u$ and
corresponding standard monomials $x^v$, the function {\tt inequalities}
creates a matrix whose columns are the vectors $u-v$. 

\par
\vskip 5 pt
\begingroup
\tteight
\baselineskip=7.4pt
\lineskip=0pt
\obeyspaces
i31\ :\ inequalities\ =\ (M)\ ->\ (\leavevmode\hss\endgraf
\ \ \ \ \ \ \ \ \ \ \ \ \ \ stds\ :=\ stdMonomials(M);\leavevmode\hss\endgraf
\ \ \ \ \ \ \ \ \ \ \ \ \ \ transpose\ matrix\ apply(numgens\ M,\ i\ ->\ (\leavevmode\hss\endgraf
\ \ \ \ \ \ \ \ \ \ \ \ \ \ \ \ \ \ flatten\ exponents(M{\char`\_}i)\ -\ \leavevmode\hss\endgraf
\ \ \ \ \ \ \ \ \ \ \ \ \ \ \ \ \ \ \ \ \ \ flatten\ exponents(stds{\char`\_}i))));\ \leavevmode\hss\endgraf
\endgroup
\penalty-1000
\par
\vskip 1 pt
\noindent
\par
\vskip 5 pt
\begingroup
\tteight
\baselineskip=7.4pt
\lineskip=0pt
\obeyspaces
i32\ :\ \ inequalities\ M\ \leavevmode\hss\endgraf
\penalty-500\leavevmode\hss\endgraf
o32\ =\ |\ 3\ \ 2\ \ 2\ \ 1\ \ 1\ \ 1\ \ 1\ \ 1\ \ 1\ \ 0\ \ 0\ \ 0\ \ 0\ \ 0\ \ 0\ \ 0\ \ |\leavevmode\hss\endgraf
\ \ \ \ \ \ |\ -1\ 1\ \ 0\ \ 3\ \ 2\ \ 1\ \ -2\ 0\ \ -1\ 5\ \ 4\ \ 3\ \ 2\ \ 1\ \ 1\ \ 0\ \ |\leavevmode\hss\endgraf
\ \ \ \ \ \ |\ -3\ -4\ -3\ -5\ -4\ -3\ 1\ \ -2\ 0\ \ -6\ -5\ -4\ -3\ -2\ -1\ 1\ \ |\leavevmode\hss\endgraf
\ \ \ \ \ \ |\ 1\ \ 1\ \ 2\ \ 1\ \ 2\ \ 3\ \ 0\ \ 4\ \ -1\ 1\ \ 2\ \ 3\ \ 4\ \ 5\ \ -1\ -6\ |\leavevmode\hss\endgraf
\ \ \ \ \ \ |\ 0\ \ 0\ \ -1\ 0\ \ -1\ -2\ 0\ \ -3\ 1\ \ 0\ \ -1\ -2\ -3\ -4\ 1\ \ 5\ \ |\leavevmode\hss\endgraf
\penalty-500\leavevmode\hss\endgraf
\ \ \ \ \ \ \ \ \ \ \ \ \ \ \ 5\ \ \ \ \ \ \ \ 16\leavevmode\hss\endgraf
o32\ :\ Matrix\ ZZ\ \ <---\ ZZ\leavevmode\hss\endgraf
\endgroup
\penalty-1000
\par
\vskip 1 pt
\noindent

It is convenient to simplify the output of the next procedure 
using the following program to divide an integer vector 
by the g.c.d. of its components. We also load {\tt polarCone.m2} 
which is needed in {\tt decideCoherence} below.

\par
\vskip 5 pt
\begingroup
\tteight
\baselineskip=7.4pt
\lineskip=0pt
\obeyspaces
i33\ :\ \ primitive\ :=\ (L)\ ->\ (\leavevmode\hss\endgraf
\ \ \ \ \ \ \ \ \ \ \ n\ :=\ {\char`\#}L-1;\ g\ :=\ L{\char`\#}n;\leavevmode\hss\endgraf
\ \ \ \ \ \ \ \ \ \ \ while\ n\ >\ 0\ do\ (n\ =\ n-1;\ g\ =\ gcd(g,\ L{\char`\#}n););\leavevmode\hss\endgraf
\ \ \ \ \ \ \ \ \ \ \ if\ g\ ===\ 1\ then\ L\ else\ apply(L,\ i\ ->\ i\ //\ g));\leavevmode\hss\endgraf
\endgroup
\penalty-1000
\par
\vskip 1 pt
\noindent

\par
\vskip 5 pt
\begingroup
\tteight
\baselineskip=7.4pt
\lineskip=0pt
\obeyspaces
i34\ :\ \ load\ "polarCone.m2"\ \leavevmode\hss\endgraf
\endgroup
\penalty-1000
\par
\vskip 1 pt
\noindent

\par
\vskip 5 pt
\begingroup
\tteight
\baselineskip=7.4pt
\lineskip=0pt
\obeyspaces
i35\ :\ decideCoherence\ =\ (M)\ ->\ (\leavevmode\hss\endgraf
\ \ \ \ \ \ \ \ \ \ \ ineqs\ :=\ inequalities\ M;\leavevmode\hss\endgraf
\ \ \ \ \ \ \ \ \ \ \ c\ :=\ first\ polarCone\ ineqs;\leavevmode\hss\endgraf
\ \ \ \ \ \ \ \ \ \ \ m\ :=\ -\ sum(numgens\ source\ c,\ i\ ->\ c{\char`\_}{\char`\{}i{\char`\}});\leavevmode\hss\endgraf
\ \ \ \ \ \ \ \ \ \ \ prods\ :=\ (transpose\ m)\ *\ ineqs;\leavevmode\hss\endgraf
\ \ \ \ \ \ \ \ \ \ \ if\ numgens\ source\ prods\ !=\ numgens\ source\ compress\ prods\leavevmode\hss\endgraf
\ \ \ \ \ \ \ \ \ \ \ then\ false\ else\ primitive\ (first\ entries\ transpose\ m));\ \leavevmode\hss\endgraf
\endgroup
\penalty-1000
\par
\vskip 1 pt
\noindent
 
Let $K$ be the cone $\{x \in {\mathbb R}^n : g \cdot x \leq 0$,
for all columns $g$ of {\tt ineqs} \}. The command {\tt
polarCone ineqs} computes a pair of matrices $P$ and $Q$ such
that $K$ is the sum of the cone generated by the columns of $P$
and the subspace generated by the columns of $Q$. Let {\tt m} be
the negative of the sum of the columns of $P$. Then {\tt m} lies
in the cone $-K$. The entries in the matrix {\tt prods} are the
dot products $g \cdot m$ for each column $g$ of {\tt ineqs}.
Since $M$ is a monomial $A$-graded ideal, it is coherent if and
only if $K$ is full dimensional which is the case if and only if
no dot product $g \cdot m$ is zero. This is the conditional in
the {\tt if .. then} statement of {\tt decideCoherence}. If $M$
is coherent, the program outputs the primitive representative of
{\tt m} and otherwise returns the boolean {\tt false}. Notice that 
if $M$ is coherent, the cone $-K$ is the Gr\"obner cone corresponding 
to $M$ and the vector {\tt m} is a weight vector $w$ such that
$in_w(I_A) = M$. We now test whether the ideal $M$ from 
line {\tt i29} is coherent.

\par
\vskip 5 pt
\begingroup
\tteight
\baselineskip=7.4pt
\lineskip=0pt
\obeyspaces
i36\ :\ decideCoherence(M)\ \leavevmode\hss\endgraf
\penalty-500\leavevmode\hss\endgraf
o36\ =\ {\char`\{}0,\ 0,\ 1,\ 15,\ 18{\char`\}}\leavevmode\hss\endgraf
\penalty-500\leavevmode\hss\endgraf
o36\ :\ List\leavevmode\hss\endgraf
\endgroup
\penalty-1000
\par
\vskip 1 pt
\noindent

Hence, $M$ is coherent: it is the initial ideal with respect to the 
weight vector $w = (0,0,1,15,18)$ of the toric ideal in our running
example. The matrix in (\ref{OurMatrix}) has 55 noncoherent
monomial $A$-graded ideals in total and here is one of them.

\par
\vskip 5 pt
\begingroup
\tteight
\baselineskip=7.4pt
\lineskip=0pt
\obeyspaces
i37\ :\ \ N\ =\ ideal(a*e,c*d,a*c,c{\char`\^}3*e,a{\char`\^}3*d,c{\char`\^}4,a*d{\char`\^}4,a{\char`\^}2*d{\char`\^}3,c*e{\char`\^}5,\leavevmode\hss\endgraf
\ \ \ \ \ \ \ \ \ \ \ \ \ \ \ \ \ c{\char`\^}2*e{\char`\^}4,d{\char`\^}7);\leavevmode\hss\endgraf
\penalty-500\leavevmode\hss\endgraf
o37\ :\ Ideal\ of\ R\leavevmode\hss\endgraf
\endgroup
\penalty-1000
\par
\vskip 1 pt
\noindent
\par
\vskip 5 pt
\begingroup
\tteight
\baselineskip=7.4pt
\lineskip=0pt
\obeyspaces
i38\ :\ \ decideCoherence(N)\ \leavevmode\hss\endgraf
\penalty-500\leavevmode\hss\endgraf
o38\ =\ false\leavevmode\hss\endgraf
\endgroup
\penalty-1000
\par
\vskip 1 pt
\noindent

In the rest of this section, we study the connection between
$A$-graded ideals and polyhedral complexes defined on $A$ which will
relate the flip graph to the Baues graph of $A$. (See \cite{Reiner} for a
survey of the Baues problem and its relatives).  Let $pos(A) := \{ Au
: u \in \R^n, u \geq 0 \}$ be the cone generated by the columns of $A$
in $\R^d$. A {\em polyhedral subdivision} $\Delta$ of $A$ is a
collection of full dimensional subcones $pos(A_{\sigma})$ of $pos(A)$
such that the union of these subcones is $pos(A)$ and the intersection
of any two subcones is a face of each.  Here $A_{\sigma} := \{a_j : j
\in \sigma \subseteq \{1,\ldots,n\} \}$.  It is customary to identify 
$\Delta$ with the set of sets $\{ \sigma : pos(A_{\sigma}) \in \Delta
\}$. If every cone in the 
subdivision $\Delta$ is simplicial (the number of extreme rays of the
cone equals the dimension of the cone), we say that $\Delta$ is a {\em
  triangulation} of $A$. The simplicial complex corresponding
to a triangulation $\Delta$ is uniquely obtained by including in
$\Delta$ all the subsets of every $\sigma \in \Delta$. We refer the
reader to \cite[\S 8]{St2} for more details.

For each $\sigma \in \Delta$, let $I_{\sigma}$ be the prime ideal 
that is the sum of the toric ideal $I_{A_{\sigma}}$ and the monomial 
ideal $\langle x_j :j \not \in \sigma \rangle$. Recall that 
$(\k^{\ast})^n$ acts on $R$ by scaling variables : $\lambda
\mapsto \lambda \cdot x := (\lambda_1 x_1, \ldots, \lambda_n x_n)$. Two
ideals $J$ and $J'$ are said to be 
{\em torus isomorphic} if $J = \lambda \cdot J'$ for some $\lambda \in 
(\k^{\ast})^n$. The following theorem shows that polyhedral
subdivisions of $A$ are related to $A$-graded ideals via their 
radicals.

\begin{theorem}(Theorem~10.10 \cite[\S 10]{St2}) \label{polysubdivisions}
  If $I$ is any $A$-graded ideal, then there exists a polyhedral
  subdivision $\Delta(I)$ of $A$ such that $\sqrt{I} = \cap_{\sigma
    \in \Delta(I)} J_{\sigma}$ where each component $J_{\sigma}$ is a
  prime ideal that is torus isomorphic to $I_{\sigma}$.
\end{theorem}

We say that $\Delta(I)$ supports the $A$-graded ideal $I$.
When $M$ is a monomial $A$-graded ideal, $\Delta(M)$ is a 
triangulation of $A$. In particular, if $M$ is coherent (i.e, $M =
in_w(I_A)$ for some weight vector $w$), then $\Delta(M)$ is the {\em
  regular} or {\em coherent} triangulation of $A$ induced by $w$
\cite[\S 8]{St2}. The coherent triangulations of $A$ are in bijection
with the vertices of the {\em secondary polytope} of $A$ \cite{BFS},
\cite{GKZ}.  

It is convenient to represent a triangulation $\Delta$ of $A$ by its 
{\em Stanley-Reisner} ideal $I_{\Delta} := \langle x_{i_1}x_{i_2}
\cdots x_{i_k} : \{ i_1, i_2, \ldots, i_k \}$ is a non-face of  
$\Delta \rangle$. If $M$ is a monomial $A$-graded ideal,
Theorem~\ref{polysubdivisions} implies that $I_{\Delta(M)}$ is the 
radical of $M$. Hence, we will represent triangulations 
of $A$ by their Stanley-Reisner ideals. As seen below, the matrix in
our running example has eight distinct triangulations 
corresponding to the eight distinct radicals of the 281 monomial 
$A$-graded ideals computed earlier. All eight are coherent.\\

\begin{tabular}{lll}
{$\{\{1,2\},\{2,3\},\{3,4\},\{4,5\}\}$}
&\,\,\,\,$\leftrightarrow$\,\,\,\,& $\langle ac, ad, ae, bd, be, ce
\rangle$ \\  
{$\{\{1,3\},\{3,4\},\{4,5\}\}$} &\,\,\,\,$\leftrightarrow$\,\,\,\,&
$\langle b, ad, ae, ce \rangle$ \\  
{$\{\{1,2\},\{2,4\},\{4,5\}\}$} &\,\,\,\,$\leftrightarrow$\,\,\,\,&
$\langle c, ad, ae, be \rangle$ \\ 
{$\{\{1,2\},\{2,3\},\{3,5\}\}$} &\,\,\,\,$\leftrightarrow$\,\,\,\,&
$\langle d, ac, ae, be \rangle$ \\  
{$\{\{1,3\},\{3,5\}\}$} &\,\,\,\,$\leftrightarrow$\,\,\,\,& $\langle
b, d, ae \rangle$ \\ 
{$\{\{1,4\},\{4,5\}\}$} &\,\,\,\,$\leftrightarrow$\,\,\,\,& $\langle
b, c, ae \rangle$ \\  
{$\{\{1,2\},\{2,5\}\}$} &\,\,\,\,$\leftrightarrow$\,\,\,\,& $\langle
c, d, ae \rangle$ \\ 
{$\{\{1,5\}\}$} &\,\,\,\,$\leftrightarrow$\,\,\,\,& $\langle b, c, d
\rangle$   
\end{tabular}\\

The Baues graph of $A$ is a graph on all the triangulations of
$A$ in which two triangulations are adjacent if they differ by a
single {\em bistellar flip} \cite{Reiner}. The {\em Baues problem} from
discrete geometry asked whether the Baues graph of a point
configuration can be disconnected for some $A$. Every edge of the
secondary polytope of $A$ corresponds to a bistellar flip and hence,
the subgraph of the Baues graph that is induced by the coherent
triangulations of $A$ is indeed connected: it is precisely the edge
graph of the secondary polytope of $A$.  The Baues problem was
recently settled by Santos \cite{San} who gave an example of a six
dimensional point configuration with $324$ points for which there is
an isolated (necessarily non-regular) triangulation.

Santos' configuration would also have a disconnected flip graph and hence
a disconnected toric Hilbert scheme if it were true that {\em every} 
triangulation of $A$ supports a monomial $A$-graded
ideal. However, Peeva has shown that this need not be the case
(Theorem~10.13 in \cite[\S 10]{St2}). Hence, the map from the set of
all monomial $A$-graded ideals to the set of all triangulations of
$A$ that sends $M \mapsto \Delta(M)$ is not always
surjective, and it is unknown whether Santos' $6 \times 324$ 
configuration has a disconnected toric Hilbert scheme.  

Thus, even though one cannot in general conclude that the existence of
a disconnected Baues graph implies the existence of a disconnected
flip graph, there is an important special situation in which such a
conclusion is possible. We call an integer matrix $A$ of full row rank
{\em unimodular} if the absolute value of each of its non-zero maximal
minors is the same constant. A matrix $A$ is unimodular if and only if
every monomial $A$-graded ideal is square-free. For a unimodular
matrix $A$, the Baues graph of $A$ coincides with the flip graph of
$A$, and as you might expect, Santos' configuration is not unimodular.

\begin{theorem} (Lemma~10.14 \cite[\S 10]{St2}) \label{unimodular}
If $A$ is unimodular, then each triangulation of $A$  
supports a unique (square-free) monomial $A$-graded ideal. Moreover, a
monomial $A$-graded ideal is coherent if and only if the 
triangulation supporting it is coherent.
\end{theorem}

Theorem~\ref{unimodular} provides a competitive algebraic algorithm
for computing all the triangulations of a unimodular matrix since 
they are precisely the polyhedral complexes supporting monomial
$A$-graded ideals. Then we could enumerate the connected component 
of a coherent monomial $A$-graded ideal in the flip graph of $A$ to
decide whether the Baues/flip graph is disconnected. 

Let $\Delta_r$ be the standard $r$-simplex which 
is the convex hull of the $r+1$ unit vectors in $\R^{r+1}$ and let 
$A(r,s)$ be the $(r+s+2) \times (r+1)(s+1)$ matrix whose columns 
are the products of the vertices of $\Delta_r$ and $\Delta_s$. All 
matrices of type $A(r,s)$ are unimodular. From the
product of two triangles we get $$A(2,2) := 
\left ( \begin{array}{ccccccccc}
1&1&1&0&0&0&0&0&0\\
0&0&0&1&1&1&0&0&0\\
0&0&0&0&0&0&1&1&1\\
1&0&0&1&0&0&1&0&0\\
0&1&0&0&1&0&0&1&0\\
0&0&1&0&0&1&0&0&1 \end{array} \right ).$$
We can now use our algebraic algorithms to compute all
the triangulations of $A(2,2)$. Since \Mtwo requires the first entry 
of the degree of every variable in a ring to be positive, we use 
the following matrix with the same row space as $A(2,2)$ for our 
computation:

\par
\vskip 5 pt
\begingroup
\tteight
\baselineskip=7.4pt
\lineskip=0pt
\obeyspaces
i39\ :\ \ A22\ =\leavevmode\hss\endgraf
\ \ \ \ \ \ \ \ {\char`\{}{\char`\{}1,1,1,1,1,1,1,1,1{\char`\}},{\char`\{}0,0,0,1,1,1,0,0,0{\char`\}},{\char`\{}0,0,0,0,0,0,1,1,1{\char`\}},\leavevmode\hss\endgraf
\ \ \ \ \ \ \ \ {\char`\{}1,0,0,1,0,0,1,0,0{\char`\}},{\char`\{}0,1,0,0,1,0,0,1,0{\char`\}},{\char`\{}0,0,1,0,0,1,0,0,1{\char`\}}{\char`\}};\ \leavevmode\hss\endgraf
\endgroup
\penalty-1000
\par
\vskip 1 pt
\noindent
\par
\vskip 5 pt
\begingroup
\tteight
\baselineskip=7.4pt
\lineskip=0pt
\obeyspaces
i40\ :\ \ I22\ =\ toricIdeal\ A22\leavevmode\hss\endgraf
\penalty-500\leavevmode\hss\endgraf
o40\ =\ ideal\ (f*h\ -\ e*i,\ c*h\ -\ b*i,\ f*g\ -\ d*i,\ e*g\ -\ d*h,\ c*g\ -\ a*i,\ b* $\cdot\cdot\cdot$\leavevmode\hss\endgraf
\penalty-500\leavevmode\hss\endgraf
o40\ :\ Ideal\ of\ R\leavevmode\hss\endgraf
\endgroup
\penalty-1000
\par
\vskip 1 pt
\noindent
The ideal {\tt I22} is generated by the 2 by 2 minors of a 3 by 3
matrix of indeterminates.  This is the ideal of $\P^2 \times \P^2$
embedded in $\P^8$ via the Segre embedding.
\par
\vskip 5 pt
\begingroup
\tteight
\baselineskip=7.4pt
\lineskip=0pt
\obeyspaces
i41\ :\ \ Graver22\ =\ graver\ I22;\ \leavevmode\hss\endgraf
\penalty-500\leavevmode\hss\endgraf
\ \ \ \ \ \ \ \ \ \ \ \ \ \ 1\ \ \ \ \ \ \ 15\leavevmode\hss\endgraf
o41\ :\ Matrix\ R\ \ <---\ R\leavevmode\hss\endgraf
\endgroup
\penalty-1000
\par
\vskip 1 pt
\noindent
\par
\vskip 5 pt
\begingroup
\tteight
\baselineskip=7.4pt
\lineskip=0pt
\obeyspaces
i42\ :\ \ generateAmonos(Graver22);\ \leavevmode\hss\endgraf
\endgroup
\penalty-1000
\par
\vskip 1 pt
\noindent
\par
\vskip 5 pt
\begingroup
\tteight
\baselineskip=7.4pt
\lineskip=0pt
\obeyspaces
i43\ :\ \ {\char`\#}monos\ \leavevmode\hss\endgraf
\penalty-500\leavevmode\hss\endgraf
o43\ =\ 108\leavevmode\hss\endgraf
\endgroup
\penalty-1000
\par
\vskip 1 pt
\noindent
\par
\vskip 5 pt
\begingroup
\tteight
\baselineskip=7.4pt
\lineskip=0pt
\obeyspaces
i44\ :\ \ scan(0..9,i->print\ toString\ monos{\char`\#}i)\ \leavevmode\hss\endgraf
{\char`\{}f*h,\ c*h,\ f*g,\ e*g,\ c*g,\ b*g,\ c*e,\ c*d,\ b*d{\char`\}}\leavevmode\hss\endgraf
{\char`\{}f*h,\ d*h,\ c*h,\ f*g,\ c*g,\ b*g,\ c*e,\ c*d,\ b*d{\char`\}}\leavevmode\hss\endgraf
{\char`\{}d*i,\ f*h,\ d*h,\ c*h,\ c*g,\ b*g,\ c*e,\ c*d,\ b*d{\char`\}}\leavevmode\hss\endgraf
{\char`\{}e*i,\ c*h,\ f*g,\ e*g,\ c*g,\ b*g,\ c*e,\ c*d,\ b*d{\char`\}}\leavevmode\hss\endgraf
{\char`\{}e*i,\ d*i,\ c*h,\ e*g,\ c*g,\ b*g,\ c*e,\ c*d,\ b*d{\char`\}}\leavevmode\hss\endgraf
{\char`\{}e*i,\ d*i,\ d*h,\ c*h,\ c*g,\ b*g,\ c*e,\ c*d,\ b*d{\char`\}}\leavevmode\hss\endgraf
{\char`\{}f*h,\ c*h,\ f*g,\ e*g,\ c*g,\ b*g,\ c*e,\ a*e,\ c*d{\char`\}}\leavevmode\hss\endgraf
{\char`\{}e*i,\ c*h,\ f*g,\ e*g,\ c*g,\ b*g,\ c*e,\ a*e,\ c*d,\ b*d*i{\char`\}}\leavevmode\hss\endgraf
{\char`\{}e*i,\ c*h,\ f*g,\ e*g,\ c*g,\ b*g,\ c*e,\ a*e,\ c*d,\ a*f*h{\char`\}}\leavevmode\hss\endgraf
{\char`\{}e*i,\ d*i,\ c*h,\ e*g,\ c*g,\ b*g,\ c*e,\ a*e,\ c*d{\char`\}}\leavevmode\hss\endgraf
\endgroup
\penalty-1000
\par
\vskip 1 pt
\noindent

Thus there are 108 monomial $A(2,2)$-graded ideals and 
{\tt decideCoherence} will check that all of them 
are coherent. Since $A(2,2)$ is unimodular, each monomial 
$A(2,2)$-graded ideal is square-free and is hence 
radical. These 108 ideals represent the 108 triangulations of 
$A(2,2)$ and we have listed ten of them above.
The flip graph (equivalently, Baues graph) of $A(2,2)$ is connected.
However, it is unknown whether the Baues graph of $A(r,s)$ is 
connected for all values of $(r,s)$.

\section{Local Equations}
Consider the reduced Gr\"obner basis of a toric ideal $I_A$ for a
term order $w$:
\begin{equation}
\label{GrobnerBasis} \bigl\{ \,
 x^{u_1} -  x^{v_1} \, , \,\, x^{u_2} -  x^{v_2} \,, \,\, \ldots \, , \,\,
x^{u_r} -  x^{v_r}\, \bigr\} .
\end{equation}
The initial ideal $\,M = in_w(I_A) = \langle x^{u_1}, x^{u_2}, \ldots,
x^{u_r} \rangle \,$ is a coherent monomial $A$-graded ideal. In
particular, it is a $(\k^*)^n$-fixed point on the toric Hilbert scheme
$Hilb_A$.  We shall explain a method, due to Peeva and Stillman
\cite{PS2}, for computing local equations of $Hilb_A$ around such a
fixed point.  A variant of this method also works for computing the
local equations around a non-coherent monomial ideal $M$, but that
variant involves local algebra, specifically Mora's tangent cone
algorithm, which is not yet fully implemented in \Mtwo. See \cite{PS2}
for details.

We saw how to compute the flip graph of $A$ in Section~1. The vertices
of this graph are the $(\k^*)^n$-fixed points $M$ and its edges
correspond to the $(\k^*)^n$-fixed curves.  By computing and
decomposing the local equations around each $M$, we get a complete
description of the scheme $Hilb_A$.

The first step is to introduce a new variable $ \, z_i \,$ for each 
binomial in our Gr\"obner basis (\ref{GrobnerBasis}) and to consider 
the following $r$ binomials:
\begin{equation}
\label{FlatFamily}
 x^{u_1} -  z_1 \cdot x^{v_1} \,, \,\,
 x^{u_2} - z_2  \cdot x^{v_2} \,,\, \, \ldots \, ,
\,\, x^{u_r} -   z_r \cdot x^{v_r} 
\end{equation}
in the polynomial ring $\k[x,z]$ in  $n+r$ indeterminates.
The term order $w$ can be extended to an elimination term order
in $\k[x,z]$ so that $x^{u_i}$ is the leading term of
$ x^{u_i} -  z_i \cdot x^{v_i} $ for all $i$. 
We compute the minimal first syzygies
of the monomial ideal $M$, and form the
corresponding $S$-pairs of binomials in (\ref{FlatFamily}).
For each $S$-pair
$$
\frac{lcm(x^{u_i},x^{u_j})}{x^{u_i}} \cdot (x^{u_i} - z_i \cdot
x^{v_i} ) \,\,\, - \,\,\, \frac{lcm(x^{u_i},x^{u_j})}{x^{u_j}} \cdot
(x^{u_j} - z_j \cdot x^{v_j}) $$
we compute a normal form with respect
to (\ref{FlatFamily}) using the extended term order $w$.  The result
is a binomial in $\k[x,z]$ which factors as 
$$  x^\alpha \cdot z^\beta \cdot  ( z^\gamma - z^\delta ) , $$
where $\alpha \in \N^n$ and $\beta,\gamma,\delta \in \N^r$.
Note that this normal form is not unique but depends on our
choice of a reduction path.
Let $J_M$ denote the ideal in $\k[z_1 , \ldots, z_r]$ generated by all 
binomials $\, z^\beta \cdot  ( z^\gamma - z^\delta ) \,$
gotten from normal forms of all the $S$-pairs considered above.

\begin{proposition} {\rm (Peeva-Stillman \cite{PS2}) } \label{localeqns}
The ideal $J_M$ is independent of the reduction paths chosen.
It defines a subscheme of $\k^r$ isomorphic to
an affine open neighborhood of the point $M$ on 
the toric Hilbert scheme $Hilb_A$.
\end{proposition}

We apply this technique to compute a particularly interesting affine
chart of $Hilb_A$ for our running example.
Consider the following set of $13$ binomials:
\begin{eqnarray*}
& \bigl\{ \,a e - z_1 b d ,  \,
 c d - z_2 b e , \,
 a c - z_3 b^2 , \,
 a^2 d^2 - z_4 c^3 e , \,
 a^2 b d - z_5 c^4 , \\ &
 a^3 d - z_6 b c^3 , \,
 c^2 e^3 - z_7 a d^4 , \, 
 c^3 e^2 - z_8 a b d^3 , \,
 c^4 e - z_9 a b^2 d^2 , \\ &
 c^5 - z_{10} a b^3 d , \,
 c e^5 - z_{11} d^6 , \,
 a d^5 - z_{12} b c e^4 , \,
 b e^6 - z_{13} d^7  \, \bigr\}.
\end{eqnarray*}
If we set $\, z_1 = z_2 = \cdots = z_{13} = 1\,$
then we get a generating set for the toric ideal $I_A$.
The $13$ monomials obtained by setting
$\, z_1 = z_2 = \cdots = z_{13} = 0 \,$
generate the initial monomial ideal $ M = in_w (I_A)$
with respect to the weight vector $w = (9, 3, 5, 0, 0)$.
Thus $M$ is one of the $226$ coherent monomial 
$A$-graded ideals of our running example. The above set of 
13 binomials in $\k[x,z]$ give the universal family 
for $Hilb_A$ around this $M$.

The local chart of $Hilb_A$ around the point $M$
is a subscheme of affine space $\k^{13}$ with coordinates 
$z_1, \ldots, z_{13}$, whose
defining equations are obtained as follows.
Extend the weight vector $w$ by assigning
weight zero to all variables $z_i$, so that
the first term in each of the above $13$ binomials
is the leading term. For each pair of binomials corresponding to a  
minimal syzygy of $M$, form their $S$-pair and then reduce it to
normal form with respect to the $13$ binomials above.
For instance,
$$
S \bigl(
 c^5 - z_{10} a b^3 d , 
 c e^5 - z_{11} d^6 \bigr)
\,\,\, = \,\,\,
 z_{11} c^4 d^6  - z_{10} a b^3 d e^5
\, \longrightarrow \,
b^4 d^2 e^4 \cdot (z_2^4 z_{11} - z_1 z_{10}).
$$
Each such normal form is a monomial in $a,b,c,d,e$
times a binomial in $z_1, \ldots, z_{13}$.
The set of all these binomials generates the ideal $J_M$ of
local equations of $Hilb_A$ around $M$.
In our example, $J_M$ is generated by $27$ nonzero binomials.
This computation can be done in \Mtwo using
the procedure {\tt localCoherentEquations}.

\vskip .2cm

\par
\vskip 5 pt
\begingroup
\tteight
\baselineskip=7.4pt
\lineskip=0pt
\obeyspaces
i45\ :\ \ localCoherentEquations\ =\ (IA)\ ->\ (\leavevmode\hss\endgraf
\ \ \ \ \ \ \ \ \ \ \ --\ IA\ is\ the\ toric\ ideal\ of\ A\ living\ in\ a\ ring\ equipped\leavevmode\hss\endgraf
\ \ \ \ \ \ \ \ \ \ \ --\ with\ weight\ order\ w,\ if\ we\ are\ computing\ the\ local\ \leavevmode\hss\endgraf
\ \ \ \ \ \ \ \ \ \ \ --\ equations\ about\ the\ initial\ ideal\ of\ IA\ w.r.t.\ w.\leavevmode\hss\endgraf
\ \ \ \ \ \ \ \ \ \ \ R\ :=\ ring\ IA;\leavevmode\hss\endgraf
\ \ \ \ \ \ \ \ \ \ \ w\ :=\ (monoid\ R).Options.Weights;\leavevmode\hss\endgraf
\ \ \ \ \ \ \ \ \ \ \ M\ :=\ ideal\ leadTerm\ IA;\leavevmode\hss\endgraf
\ \ \ \ \ \ \ \ \ \ \ S\ :=\ first\ entries\ ((gens\ M)\ {\char`\%}\ IA);\leavevmode\hss\endgraf
\ \ \ \ \ \ \ \ \ \ \ --\ Make\ the\ universal\ family\ J\ in\ a\ new\ ring.\leavevmode\hss\endgraf
\ \ \ \ \ \ \ \ \ \ \ nv\ :=\ numgens\ R;\ n\ :=\ numgens\ M;\leavevmode\hss\endgraf
\ \ \ \ \ \ \ \ \ \ \ T\ =\ (coefficientRing\ R)[generators\ R,\ z{\char`\_}1\ ..\ z{\char`\_}n,\ \leavevmode\hss\endgraf
\ \ \ \ \ \ \ \ \ \ \ \ \ \ \ \ \ \ \ \ \ \ \ \ \ \ \ \ \ \ \ \ \ \ \ \ Weights\ =>\ flatten\ splice{\char`\{}w,\ n:0{\char`\}},\leavevmode\hss\endgraf
\ \ \ \ \ \ \ \ \ \ \ \ \ \ \ \ \ \ \ \ \ \ \ \ \ \ \ \ \ \ \ \ \ \ \ \ MonomialSize\ =>\ 16];\leavevmode\hss\endgraf
\ \ \ \ \ \ \ \ \ \ \ M\ =\ substitute(generators\ M,T);\leavevmode\hss\endgraf
\ \ \ \ \ \ \ \ \ \ \ S\ =\ apply(S,\ s\ ->\ substitute(s,T));\leavevmode\hss\endgraf
\ \ \ \ \ \ \ \ \ \ \ J\ =\ ideal\ apply(n,\ i\ ->\ \leavevmode\hss\endgraf
\ \ \ \ \ \ \ \ \ \ \ \ \ \ \ \ \ \ \ \ \ M{\char`\_}(0,i)\ -\ T{\char`\_}(nv\ +\ i)\ *\ S{\char`\_}i);\leavevmode\hss\endgraf
\ \ \ \ \ \ \ \ \ \ \ --\ Find\ the\ ideal\ Ihilb\ of\ local\ equations\ about\ M:\leavevmode\hss\endgraf
\ \ \ \ \ \ \ \ \ \ \ spairs\ :=\ (gens\ J)\ *\ (syz\ M);\leavevmode\hss\endgraf
\ \ \ \ \ \ \ \ \ \ \ g\ :=\ forceGB\ gens\ J;\leavevmode\hss\endgraf
\ \ \ \ \ \ \ \ \ \ \ B\ =\ (coefficientRing\ R)[z{\char`\_}1\ ..\ z{\char`\_}n,MonomialSize=>16];\leavevmode\hss\endgraf
\ \ \ \ \ \ \ \ \ \ \ Fones\ :=\ map(B,T,\ matrix(B,{\char`\{}splice\ {\char`\{}nv:1{\char`\}}{\char`\}})\ |\ vars\ B);\leavevmode\hss\endgraf
\ \ \ \ \ \ \ \ \ \ \ Ihilb\ :=\ ideal\ Fones\ (spairs\ {\char`\%}\ g);\leavevmode\hss\endgraf
\ \ \ \ \ \ \ \ \ \ \ Ihilb\leavevmode\hss\endgraf
\ \ \ \ \ \ \ \ \ \ \ );\leavevmode\hss\endgraf
\endgroup
\penalty-1000
\par
\vskip 1 pt
\noindent

Suppose we wish to calculate the local equations about $M =
in_w(I_A)$.  The input to {\tt localCoherentEquations} is the
toric ideal $I_A$ living in a polynomial ring equipped with the 
weight order specified by $w$. This is done as follows.

\par
\vskip 5 pt
\begingroup
\tteight
\baselineskip=7.4pt
\lineskip=0pt
\obeyspaces
i46\ :\ IA\ =\ toricIdeal\ A;\ \leavevmode\hss\endgraf
\penalty-500\leavevmode\hss\endgraf
o46\ :\ Ideal\ of\ R\leavevmode\hss\endgraf
\endgroup
\penalty-1000
\par
\vskip 1 pt
\noindent
\par
\vskip 5 pt
\begingroup
\tteight
\baselineskip=7.4pt
\lineskip=0pt
\obeyspaces
i47\ :\ Y\ =\ QQ[a..e,\ MonomialSize\ =>\ 16,\ \leavevmode\hss\endgraf
\ \ \ \ \ \ \ \ \ \ \ \ \ \ \ \ \ \ Degrees\ =>\ transpose\ A,\ Weights\ =>\ {\char`\{}9,3,5,0,0{\char`\}}];\leavevmode\hss\endgraf
\endgroup
\penalty-1000
\par
\vskip 1 pt
\noindent
\par
\vskip 5 pt
\begingroup
\tteight
\baselineskip=7.4pt
\lineskip=0pt
\obeyspaces
i48\ :\ IA\ =\ substitute(IA,Y);\ \leavevmode\hss\endgraf
\penalty-500\leavevmode\hss\endgraf
o48\ :\ Ideal\ of\ Y\leavevmode\hss\endgraf
\endgroup
\penalty-1000
\par
\vskip 1 pt
\noindent

The initial ideal $M$ is calculated in the third line of the
algorithm, and {\tt S} stores the standard monomials of $M$ of the
same degrees as the minimal generators of $M$. We could have
calculated {\tt S} using our old procedure {\tt stdMonomials} but this
involves computing the monomials in $R_b$ for various values of $b$
which can be slow on large examples. As by-products, {\tt
  localCoherentEquations} also gets {\tt J}, the ideal of the
universal family for $Hilb_A$ about $M$, the ring {\tt T} of this
ideal, and the ring {\tt B} of {\tt Ihilb} which is the ideal of the
affine patch of $Hilb_A$ about $M$. The matrix {\tt spairs} contains
all the $S$-pairs between generators of {\tt J} corresponding to the
minimal first syzygies of $M$. The command {\tt forceGB} is used to
declare the generators of {\tt J} to be a Gr\"obner basis, and {\tt
  Fones} is the ring map from {\tt T} to {\tt B} that sends each of
$a,b,c,d,e$ to one and the $z$ variables to themselves.  The columns
of the matrix {\tt (spairs \% g)} are the normal forms of the
polynomials in {\tt spairs} with respect to the forced Gr\"obner basis
{\tt g} and the ideal {\tt Ihilb} of local equations is generated by
the image of these normal forms in the ring {\tt B} under the map {\tt
  Fones}.

\par
\vskip 5 pt
\begingroup
\tteight
\baselineskip=7.4pt
\lineskip=0pt
\obeyspaces
i49\ :\ JM\ =\ localCoherentEquations(IA)\leavevmode\hss\endgraf
\penalty-500\leavevmode\hss\endgraf
\ \ \ \ \ \ \ \ \ \ \ \ \ \ \ \ \ \ \ \ \ \ \ \ \ \ \ \ \ \ \ \ \ \ \ \ \ \ \ \ \ \ \ \ \ \ \ \ \ \ \ \ \ \ \ \ \ \ \ \ \ \ \ \ \ \ \ \ \ \  $\cdot\cdot\cdot$\leavevmode\hss\endgraf
o49\ =\ ideal\ (z\ z\ \ -\ z\ ,\ z\ z\ \ -\ z\ ,\ -\ z\ z\ \ +\ z\ ,\ -\ z\ z\ \ +\ z\ ,\ -\ z\ z\ \ +\  $\cdot\cdot\cdot$\leavevmode\hss\endgraf
\ \ \ \ \ \ \ \ \ \ \ \ \ \ 1\ 2\ \ \ \ 3\ \ \ 1\ 2\ \ \ \ 3\ \ \ \ \ 4\ 7\ \ \ \ 2\ \ \ \ \ 5\ 8\ \ \ \ 2\ \ \ \ \ 1\ 5\ \ \  $\cdot\cdot\cdot$\leavevmode\hss\endgraf
\penalty-500\leavevmode\hss\endgraf
o49\ :\ Ideal\ of\ B\leavevmode\hss\endgraf
\endgroup
\penalty-1000
\par
\vskip 1 pt
\noindent

\vskip .2cm

Removing duplications among the generators, \\

$J_M = \langle
z_1-z_{10}z_{11},
z_2-z_4z_7,
z_2-z_5z_8,
z_2-z_{11}z_{12},
z_2-z_1z_{11}z_{13},\\
z_3-z_1z_2,
z_3-z_5z_9,
z_4-z_1z_5,
z_6-z_3z_5,
z_6-z_1z_2z_5,
z_7-z_1z_{10},
z_8-z_1z_7,\\
z_9-z_1z_8,
z_{12}-z_1z_{13},
z_1z_2-z_5z_9,
z_1z_2-z_1z_5z_8,
z_1z_2-z_1^2z_4z_{10},
z_1z_2-z_1^2z_5z_7,\\
z_1z_2-z_1z_{11}z_{12},
z_1z_2-z_2z_{10}z_{11},
z_1^3z_4-z_3z_{11},
z_1z_5z_8-z_4z_8,
z_2z_{10}-z_1z_{12},\\
z_3z_4-z_1z_6,
z_3z_7-z_2z_8,
z_3z_8-z_2z_9,
z_3z_{10}-z_2z_7
\rangle$. \\

Notice that there are many generators of $J_M$ which have a single
variable as one of its terms. Using these generators we can remove
variables from other binomials. This is done in \Mtwo using the
subroutine {\tt removeRedundantVariables} which is the main ingredient
of the package {\tt minPres.m2} for computing the minimal
presentations of polynomial quotient rings. Both {\tt
  removeRedundantVariables} and {\tt minPres.m2} are explained in
Appendix 2. The command {\tt removeRedundantVariables} applied to an
ideal in a polynomial ring (not quotient ring) creates a ring map from
the ring to itself that sends the redundant variables to polynomials 
in the non-redundant variables and the non-redundant variables to 
themselves. Applying this to our ideal $J_M$ we obtain the following 
simplifications.

\par
\vskip 5 pt
\begingroup
\tteight
\baselineskip=7.4pt
\lineskip=0pt
\obeyspaces
i50\ :\ load\ "minPres.m2";\leavevmode\hss\endgraf
\endgroup
\penalty-1000
\par
\vskip 1 pt
\noindent
\par
\vskip 5 pt
\begingroup
\tteight
\baselineskip=7.4pt
\lineskip=0pt
\obeyspaces
i51\ :\ G\ =\ removeRedundantVariables\ JM\ \leavevmode\hss\endgraf
\penalty-500\leavevmode\hss\endgraf
\ \ \ \ \ \ \ \ \ \ \ \ \ \ \ \ \ \ \ \ \ \ \ \ \ \ 3\ \ 2\ \ \ \ \ \ 4\ \ 3\ \ \ \ \ \ \ \ \ \ \ \ \ \ \ \ \ \ 2\ 4\ \ 3\ \ \ \ 2\  $\cdot\cdot\cdot$\leavevmode\hss\endgraf
o51\ =\ map(B,B,{\char`\{}z\ \ z\ \ ,\ z\ z\ \ z\ \ ,\ z\ z\ \ z\ \ ,\ z\ z\ \ z\ \ ,\ z\ ,\ z\ z\ \ z\ \ ,\ z\ \  $\cdot\cdot\cdot$\leavevmode\hss\endgraf
\ \ \ \ \ \ \ \ \ \ \ \ \ \ \ \ 10\ 11\ \ \ 5\ 10\ 11\ \ \ 5\ 10\ 11\ \ \ 5\ 10\ 11\ \ \ 5\ \ \ 5\ 10\ 11\ \ \ 10 $\cdot\cdot\cdot$\leavevmode\hss\endgraf
\penalty-500\leavevmode\hss\endgraf
o51\ :\ RingMap\ B\ <---\ B\leavevmode\hss\endgraf
\endgroup
\penalty-1000
\par
\vskip 1 pt
\noindent
\par
\vskip 5 pt
\begingroup
\tteight
\baselineskip=7.4pt
\lineskip=0pt
\obeyspaces
i52\ :\ ideal\ gens\ gb(G\ JM)\leavevmode\hss\endgraf
\penalty-500\leavevmode\hss\endgraf
\ \ \ \ \ \ \ \ \ \ \ \ \ \ \ 3\ \ 2\ \ \ \ \ \ \ \ 2\leavevmode\hss\endgraf
o52\ =\ ideal(z\ z\ \ z\ \ \ -\ z\ \ z\ \ z\ \ )\leavevmode\hss\endgraf
\ \ \ \ \ \ \ \ \ \ \ \ \ 5\ 10\ 11\ \ \ \ 10\ 11\ 13\leavevmode\hss\endgraf
\penalty-500\leavevmode\hss\endgraf
o52\ :\ Ideal\ of\ B\leavevmode\hss\endgraf
\endgroup
\penalty-1000
\par
\vskip 1 pt
\noindent

Thus our affine patch of $Hilb_A$ has the coordinate ring 
$$\k[z_1,z_2,\ldots,z_{13}]/J_M \,\, \simeq \,\,
\frac{\k[z_5,z_{10},z_{11},z_{13}]}{ \langle z_5 z_{{10}}^3 z_{11}^2 -
  z_{10}z_{11}^2 z_{13} \rangle} = \frac{\k[z_5,z_{10},z_{11},z_{13}]}
{\langle (z_5 z_{10}^2 -z_{13}) z_{10}z_{11}^2 \rangle}.$$
Hence, we see immediately that there are three
components through the point $M$ on $Hilb_A$. The restriction of the
coherent component to the affine neighborhood of $M$ on $Hilb_A$ is
defined by the ideal quotient $\, (J_M : (z_1 z_2 \cdots
z_{13})^\infty) $ and hence the first of the above components 
is an affine patch of the coherent component. Locally near $M$ it is  
given by the single equation $z_5 z_{10}^2 - z_{13} = 0$ in $\A^4$. 
It is smooth and as expected, has dimension three. The second
component, $z_{10} = 0$ is also of dimension three and is smooth at $M$.
The third component, given by $z_{11}^2 = 0$ is more interesting.  It
has dimension three as well, but is not reduced.  Thus we have proved
the following result.

\begin{proposition}
The toric Hilbert scheme $Hilb_A$ of the matrix $A$ in
(\ref{OurMatrix}) is not reduced.
\end{proposition}

We can use the ring map {\tt G} from above to simplify {\tt J} so as
to involve only the four variables $z_5, z_{10},z_{11}$ and $z_{13}$.

\vskip .2cm 

\par
\vskip 5 pt
\begingroup
\tteight
\baselineskip=7.4pt
\lineskip=0pt
\obeyspaces
i53\ :\ CX\ =\ QQ[a..e,\ z{\char`\_}5,z{\char`\_}10,z{\char`\_}11,z{\char`\_}13,\ Weights\ =>\ splice\leavevmode\hss\endgraf
\ \ \ \ \ \ \ \ \ \ \ \ {\char`\{}9,3,5,0,0,0,0,0,0{\char`\}}];\leavevmode\hss\endgraf
\endgroup
\penalty-1000
\par
\vskip 1 pt
\noindent
 
\par
\vskip 5 pt
\begingroup
\tteight
\baselineskip=7.4pt
\lineskip=0pt
\obeyspaces
i54\ :\ F\ =\ map(CX,\ ring\ J,\ matrix{\char`\{}{\char`\{}a,b,c,d,e{\char`\}}{\char`\}}\ |\ \leavevmode\hss\endgraf
\ \ \ \ \ \ \ \ \ \ \ \ \ \ \ \ \ \ substitute(G.matrix,CX))\leavevmode\hss\endgraf
\penalty-500\leavevmode\hss\endgraf
\ \ \ \ \ \ \ \ \ \ \ \ \ \ \ \ \ \ \ \ \ \ \ \ \ \ \ \ \ \ \ \ \ \ \ \ \ \ \ \ \ \ 3\ \ 2\ \ \ \ \ \ 4\ \ 3\ \ \ \ \ \ \ \ \ \ \ \ \ \  $\cdot\cdot\cdot$\leavevmode\hss\endgraf
o54\ =\ map(CX,T,{\char`\{}a,\ b,\ c,\ d,\ e,\ z\ \ z\ \ ,\ z\ z\ \ z\ \ ,\ z\ z\ \ z\ \ ,\ z\ z\ \ z\ \ ,\ z $\cdot\cdot\cdot$\leavevmode\hss\endgraf
\ \ \ \ \ \ \ \ \ \ \ \ \ \ \ \ \ \ \ \ \ \ \ \ \ \ \ \ \ \ \ \ 10\ 11\ \ \ 5\ 10\ 11\ \ \ 5\ 10\ 11\ \ \ 5\ 10\ 11\ \ \  $\cdot\cdot\cdot$\leavevmode\hss\endgraf
\penalty-500\leavevmode\hss\endgraf
o54\ :\ RingMap\ CX\ <---\ T\leavevmode\hss\endgraf
\endgroup
\penalty-1000
\par
\vskip 1 pt
\noindent
Applying this map to {\tt J} we get the ideal $J1$, 
\par
\vskip 5 pt
\begingroup
\tteight
\baselineskip=7.4pt
\lineskip=0pt
\obeyspaces
i55\ :\ J1\ =\ F\ J\ \leavevmode\hss\endgraf
\penalty-500\leavevmode\hss\endgraf
\ \ \ \ \ \ \ \ \ \ \ \ \ \ \ \ \ \ \ \ \ \ \ \ \ \ \ \ \ \ \ \ \ \ \ \ \ \ \ \ \ \ \ \ 3\ \ 2\ \ \ \ \ \ \ \ \ \ 2\ \ \ 4\ \ 3\ \ \ \  $\cdot\cdot\cdot$\leavevmode\hss\endgraf
o55\ =\ ideal\ (c*d\ -\ b*e*z\ \ z\ \ ,\ a*e\ -\ b*d*z\ z\ \ z\ \ ,\ a*c\ -\ b\ z\ z\ \ z\ \ ,\ a $\cdot\cdot\cdot$\leavevmode\hss\endgraf
\ \ \ \ \ \ \ \ \ \ \ \ \ \ \ \ \ \ \ \ \ \ \ \ 10\ 11\ \ \ \ \ \ \ \ \ \ \ \ \ 5\ 10\ 11\ \ \ \ \ \ \ \ \ \ \ 5\ 10\ 11\ \ \  $\cdot\cdot\cdot$\leavevmode\hss\endgraf
\penalty-500\leavevmode\hss\endgraf
o55\ :\ Ideal\ of\ CX\leavevmode\hss\endgraf
\endgroup
\penalty-1000
\par
\vskip 1 pt
\noindent

\noindent and adding the ideal $\langle z_{11}^2 \rangle$ to $J1$ 
we obtain the universal family for the non-reduced component of
$Hilb_A$ about $M$. 

\par
\vskip 5 pt
\begingroup
\tteight
\baselineskip=7.4pt
\lineskip=0pt
\obeyspaces
i56\ :\ \ substitute(ideal(z{\char`\_}11{\char`\^}2),CX)\ +\ J1\ \leavevmode\hss\endgraf
\penalty-500\leavevmode\hss\endgraf
\ \ \ \ \ \ \ \ \ \ \ \ \ \ 2\ \ \ \ \ \ \ \ \ \ \ \ \ \ \ \ \ \ \ \ \ \ \ \ \ \ \ \ \ \ \ \ \ \ 3\ \ 2\ \ \ \ \ \ \ \ \ \ 2\ \ \ 4\ \  $\cdot\cdot\cdot$\leavevmode\hss\endgraf
o56\ =\ ideal\ (z\ \ ,\ c*d\ -\ b*e*z\ \ z\ \ ,\ a*e\ -\ b*d*z\ z\ \ z\ \ ,\ a*c\ -\ b\ z\ z\ \ z $\cdot\cdot\cdot$\leavevmode\hss\endgraf
\ \ \ \ \ \ \ \ \ \ \ \ \ \ 11\ \ \ \ \ \ \ \ \ \ \ \ \ 10\ 11\ \ \ \ \ \ \ \ \ \ \ \ \ 5\ 10\ 11\ \ \ \ \ \ \ \ \ \ \ 5\ 10\  $\cdot\cdot\cdot$\leavevmode\hss\endgraf
\penalty-500\leavevmode\hss\endgraf
o56\ :\ Ideal\ of\ CX\leavevmode\hss\endgraf
\endgroup
\penalty-1000
\par
\vskip 1 pt
\noindent

\vskip .2cm

In the rest of this section, we present an interpretation of
the ideal $J_M$ in terms of the combinatorial theory
of {\it integer programming}. See, for instance, 
\cite[\S 4]{St2} or \cite{Tho} for 
the relevant background. Our reduced Gr\"obner basis 
(\ref{GrobnerBasis}) is the {\it minimal test set} for
the family of integer programs
\begin{equation}
\label{IP}
{\rm Minimize} \quad
w \cdot u \,\,\quad
{\rm subject} \,\, {\rm to } \,\,\,
A \cdot u = b   \,\,\, {\rm and}
 \,\,\,u \in \N^n, 
\end{equation}
where $A \in \N^{d \times n}$ and $w
\in \ZZ^n$ are fixed and $b$ ranges over $\N^d$.
If $u' \in \N^n$ is any feasible solution
to (\ref{IP}), then the corresponding optimal solution
$u \in \N^n$ is computed as follows: the monomial
$x^u $ is the unique normal form of $x^{u'}$ modulo
the Gr\"obner basis (\ref{GrobnerBasis}).

Suppose we had reduced $x^{u'}$
modulo the binomials (\ref{FlatFamily}) instead of (\ref{GrobnerBasis}).
Then the output has a $z$-factor which depends on 
our choice of reduction path. To be precise, suppose the
reduction path has length $m$ and at the $j$-th step we had used the
reduction $\, x^{u_{\mu_j}} \rightarrow  z_{\mu_j} \cdot x^{v_{\mu_j}}
$. Then we would obtain the normal form
$$ \, z_{\mu_1} z_{\mu_2} z_{\mu_3} \cdots z_{\mu_m} \cdot x^u.$$
Reduction paths can have different lengths. If we take
another  path  which
has length $m'$ and  uses
$\, x^{u_{\nu_j}} \rightarrow  z_{\nu_j} \cdot x^{v_{\nu_j}} \,$
at the $j$-th step, then the output would be
$$ \, z_{\nu_1} z_{\nu_2} z_{\nu_3} \cdots z_{\nu_{m'}} \cdot x^u  .$$

\begin{theorem} \label{paths}
The ideal $J_M$ of local equations
on $Hilb_A$ is generated by the binomials
$$ \, z_{\mu_1} z_{\mu_2} z_{\mu_3} \cdots z_{\mu_m} \quad - \quad
 z_{\nu_1} z_{\nu_2} z_{\nu_3} \cdots z_{\nu_{m'}}  $$
each encoding a pair of distinct reduction sequences from a feasible
solution of (\ref{IP}) to the corresponding optimal solution
using the minimal test set in (\ref{GrobnerBasis}).
\end{theorem}

\noindent {\sl Proof: }
The given ideal is contained in $J_M$ because its generators
are differences of monomials arising from the possible
reduction paths of $\,{lcm(x^{u_i},x^{u_j})} $,
for $1 \leq i,j \leq r $. Conversely, any reduction
sequence can be transformed into an equivalent reduction sequence
using S-pair reductions. This follows from standard
arguments in the proof of Buchberger's criterion
\cite[\S 2.6, Theorem 6]{CLO}, and it implies that
the binomials  $ \, z_{\mu_1}  \cdots z_{\mu_m} -
 z_{\nu_1}  \cdots z_{\nu_{m'}}  \,$ are $\k[z]$-linear
combinations of the generators of $J_M$.
\qed

\vskip .2cm

A given feasible solution of the integer program (\ref{IP})
usually has many different reduction paths to the optimal solution
using  (\ref{GrobnerBasis}). 
For our matrix (\ref{OurMatrix}) and cost vector 
$w = (9,3,5,0,0)$, the monomial
$\, a^2 b d e^6 \,$  encodes the feasible solution $(2,1,0,1,0,6)$
of the integer program (\ref{IP}) with right hand side vector 
$\,\binom{10}{56}$.
There are $19$ different paths from this feasible solution 
to the optimal solution $(0,3,0,3,4)$ encoded by the monomial 
$\, b^3 d^3 e^4 $. The generating function for these paths is:
\begin{eqnarray*}
& z_1^2 + 3 z_1 z_2^2 z_5 z_7 + 2 z_1 z_2 z_5 z_7^2 z_{12}
 + 2 z_1 z_2 z_5 z_8 + 2 z_1 z_2 z_{12} z_{13} + z_1 z_5 z_9 \\ &
 + z_2^3 z_4 z_5 z_7^2 + z_2^3 z_4 z_{13} + z_2^3 z_5 z_{11}
 + 2 z_2 z_3 z_5 z_7 + z_3 z_5 z_7^2 z_{12} + z_3 z_5 z_8 + z_3 z_{12} z_{13}.
\end{eqnarray*}
The difference of any two monomials in this generating function 
is a valid local equation for the toric Hilbert scheme of
(\ref{OurMatrix}). For instance, the binomial 
$\,  z_3 z_5 z_7^2 z_{12} - z_3 z_{12} z_{13} \,$ lies in $J_M$,
and, conversely, $J_M$ is generated by binomials obtained in this manner.

The scheme structure of $J_M$ encodes obstructions to making certain
reductions when solving our family of integer programs. For instance,
the variable $z_3$ is a zero-divisor modulo $J_M$. If we factor it
out from the binomial $\,  z_3 z_5 z_7^2 z_{12} - z_3 z_{12} z_{13}
\in J_M \,$, we get 
$\, z_5 z_7^2 z_{12} - z_{12} z_{13} \,$
which does not lie in $J_M$. Thus there is no monomial
$a^{i_1} b^{i_2} c^{i_3} d^{i_4} e^{i_5} \,$
for which both the paths $ z_5 z_7^2 z_{12} $ and
$ z_{12} z_{13} $ are used to reach the optimum.
It would be a worthwhile combinatorial project to
study the path generating functions and their relation
to the ideal $J_M$ in more detail.

It is instructive to note that the binomials
$ \, z_{\mu_1} z_{\mu_2} \cdots z_{\mu_m} \, - \,
 z_{\nu_1} z_{\nu_2}  \cdots z_{\nu_{m'}}  $
in Theorem~\ref{paths} do not form a vector space basis
for the ideal $J_M$. We demonstrate this for the lexicographic
Gr\"obner basis (with $a \succ b \succ c \succ d \succ e$) of 
the toric ideal defining the rational normal curve of degree $4$. 
In this case, we can take $A = \left ( \begin{array}{ccccc}
1 & 1 & 1 & 1 & 1 \\ 0 & 1 & 2 & 3 & 4 \end{array} \right )$ and the 
universal family in question is :
$$ 
\bigl\{
a c    - z_1 b^2, \,
a d    - z_2 b c,\,
a e    - z_3 c^2,\,
b d    - z_4 c^2,\,
b e    - z_5 c d,\,
c e    - z_6 d^2
\bigr\}.
$$ 
The corresponding ideal of local equations is
$J_M = 
\langle z_3 - z_2  z_5, z_2 - z_1  z_4, z_5  - z_4 z_6 \rangle $,
from which we see that $M$ is a smooth point of $Hilb_A$.
The binomial $\,z_1 z_5 - z_1 z_4 z_6 \,$ lies in $J_M$
but there is no monomial which has the reduction path $z_1 z_5$
or $z_5 z_1 $ to optimality.  Indeed, any monomial
which admits the reductions $z_1 z_5$ or $z_5 z_1$ must be
divisible by either $\, a c e \, $ or $\, a b e  $.
The path generating functions for these two monomials are
$$ abe \quad \rightarrow \quad
(z_3 \, + \, z_1 z_4 z_5 \, +\, z_2 z_5) \cdot b c^2 $$
$$ ace \quad \rightarrow \quad
(z_3 \,+ \,z_1 z_4 z_5 \,+\,
z_2 z_4 z_6) \cdot c^3 . $$
Thus every reduction to optimality using $z_1$ and $z_5$ must
also use $z_4$, and we conclude that $\,z_1 z_5 - z_1 z_4 z_6 \,$
is not in the $\k$-span of the binomials listed in
Theorem~\ref{paths}.

\section{The coherent component of the toric Hilbert scheme}

In this section we study the distinguished component of the
toric Hilbert scheme $Hilb_A$. The toric ideal $I_A$ is a point 
on this component. We will show that this component need not be
normal, and will describe how its local 
and global equations can be computed using \Mtwo.
Every term order for the toric ideal $I_A$
can be realized by a weight vector 
which is an element in the lattice
$\,N = Hom_\ZZ( ker_\ZZ(A) , \ZZ) \, \simeq \, \ZZ^{n-d}$.
Two weight vectors $w$ and $w'$ in $N$ are considered {\it equivalent}
if they define the same initial ideal  $\,in_w(I_A) = in_{w'}(I_A)$.
These equivalence classes are the relatively open cones
of a projective fan $\Sigma_A$ in the lattice $N$ called the 
{\em Gr\"obner fan} of $I_A$ \cite{MR}, \cite{ST}.

\begin{theorem}
The toric ideal $I_A$ lies on a unique irreducible component of
the toric Hilbert scheme $Hilb_A$, called the coherent component.
The normalization of the coherent
component is the projective toric variety defined by 
the Gr\"obner fan of $I_A$.
\end{theorem}

\noindent {\sl Proof: }
The {\it divisor at infinity} on the toric Hilbert scheme $Hilb_A$ 
consists of all points at which at least one of the local coordinates
(around some monomial $A$-graded ideal) is zero.  This is a proper
Zariski closed subset of $Hilb_A$, consisting of all those $A$-graded
ideals which contain at least one monomial.  
By \cite[Lemma 10.12]{St2}, its complement in 
$Hilb_A$ consists of precisely  the 
orbit of $I_A$ under the action of the torus $(\k^*)^n $.
This shows that the closure of the 
$(\k^*)^n $-orbit of $I_A$  is an irreducible
component of $Hilb_A$. We call it the coherent component.

Identifying $(\k^*)^n$ with $Hom_\ZZ(\ZZ^n, \k^*)$, we note
that the stabilizer of $I_A$ consists of those linear forms
$w$ which restrict to zero on the kernel of $A$. Therefore
the coherent component is the closure in $Hilb_A$
of the  orbit of the point $I_A$ under the action of the torus
$\, N \otimes \k^* \, = \,Hom_\ZZ( ker_\ZZ(A), \k^*)$.
The $(N \otimes \k^*)$-fixed points
on this component are precisely the coherent monomial 
$A$-graded ideals, and the same holds for the
toric variety of the Gr\"obner fan. 

Fix a maximal cone $\sigma$ in the Gr\"obner fan $\Sigma_A$,
and let $M = \langle x^{u_1}, \ldots, x^{u_r} \rangle$ 
be the corresponding (monomial) 
initial ideal of $I_A$.  As before we write
$$ \left \{ x^{u_1} -  z_1 \cdot x^{v_1} \,, \,\,
 x^{u_2} - z_2  \cdot x^{v_2} \,,\, \, \ldots \, ,
\,\, x^{u_r} -   z_r \cdot x^{v_r} \right \}$$
for the universal family arising from the corresponding
reduced Gr\"obner basis of $I_A$.  Let $J_M$ be the
ideal in $\k [z_1,z_2,\ldots,z_r]$ defining this family.

The restriction of the coherent component to the 
affine neighborhood of $M$ on $Hilb_A$ is defined
by $\, J_M :  (z_1 z_2 \cdots z_r)^\infty $.
It then follows from our combinatorial description of
the ideal $J_M$ that this ideal quotient is a binomial prime ideal.
In fact, it is the ideal of algebraic relations among the 
Laurent monomials $\, x^{u_1- v_1}, \ldots, x^{u_r-v_r}$.
We conclude that the restriction of the coherent component to the
affine neighborhood of $M$ on $Hilb_A$ equals
\begin{equation}
\label{uv-algebra}
 {\rm Spec} \,\, \k \bigl[
 x^{u_1-v_1},
 x^{u_2-v_2},  \ldots,
 x^{u_r-v_r} 
\bigr] .
\end{equation}

The abelian group generated by the vectors
$\, u_1-v_1, \ldots, u_r-v_r \,$  equals
$\ker_\ZZ(A) = Hom_\ZZ(N,\ZZ)$. The cone generated
by the vectors $\, u_1-v_1, \ldots, u_r-v_r \,$  is
precisely the polar dual $\sigma^\vee$ to
the Gr\"obner cone $\sigma$.  Therefore the
normalization of the affine variety
(\ref{uv-algebra}) is the normal affine toric variety
\begin{equation}
\label{normal-uv-algebra}
 {\rm Spec} \,\, \k \bigl[
 \ker_\ZZ(A) \,\cap\, \sigma^\vee \bigr] .
\end{equation}

The various normalization maps from (\ref{normal-uv-algebra}) to
(\ref{uv-algebra}), obtained as $\sigma$ varies over the cones of 
$\Sigma_A$ are compatible with gluing along $\Sigma_A$.
The result is the desired normalization map from the
projective toric variety associated with the Gr\"obner fan of $I_A$ 
onto the coherent component of the toric Hilbert scheme $Hilb_A$.
\qed

\vskip .1cm

We now present an example which shows that the coherent component 
of $Hilb_A$ need not be normal. This example is 
derived from the matrix that appears in Example 3.15 of \cite{HM}.
This example is also mentioned in \cite{PS1} without details.
Let $d=4$ and $n=7$ and fix the matrix

\begin{equation}
\label{non-normal}
A \quad = \quad \left( \begin{array}{ccccccc}  
1 & 1 & 1 & 1 & 1 & 1 & 1 \\
0 & 6 & 7 & 5 & 8 & 4 & 3 \\
3 & 7 & 2 & 0 & 7 & 6 & 1 \\
6 & 5 & 2 & 6 & 5 & 0 & 0 \end{array} \right).
\end{equation}

The lattice $\,N =  Hom_\ZZ ( ker_\ZZ(A), \ZZ)$
is three-dimensional. The toric ideal $I_A$ is minimally 
generated by $30$ binomials of total degree between $6$ and $93$.

\par
\vskip 5 pt
\begingroup
\tteight
\baselineskip=7.4pt
\lineskip=0pt
\obeyspaces
i57\ :\ \ A\ =\ {\char`\{}{\char`\{}1,1,1,1,1,1,1{\char`\}},{\char`\{}0,6,7,5,8,4,3{\char`\}},{\char`\{}3,7,2,0,7,6,1{\char`\}},\leavevmode\hss\endgraf
\ \ \ \ \ \ \ \ \ {\char`\{}6,5,2,6,5,0,0{\char`\}}{\char`\}};\leavevmode\hss\endgraf
\endgroup
\penalty-1000
\par
\vskip 1 pt
\noindent
\par
\vskip 5 pt
\begingroup
\tteight
\baselineskip=7.4pt
\lineskip=0pt
\obeyspaces
i58\ :\ \ IA\ =\ toricIdeal\ A\ \leavevmode\hss\endgraf
\penalty-500\leavevmode\hss\endgraf
\ \ \ \ \ \ \ \ \ \ \ \ \ \ 2\ 3\ \ \ \ \ \ \ 3\ 2\ \ \ 2\ \ \ \ \ 4\ 4\ \ \ \ 8\ 4\ \ \ 4\ 3\ 6\ \ \ \ 7\ 2\ 4\ \ \ \ \ 4\  $\cdot\cdot\cdot$\leavevmode\hss\endgraf
o58\ =\ ideal\ (a\ c\ e\ -\ b*d\ f\ ,\ a\ c*d*e\ f\ \ -\ b\ g\ ,\ d\ e\ f\ \ -\ b\ c\ g\ ,\ a*b\ c $\cdot\cdot\cdot$\leavevmode\hss\endgraf
\penalty-500\leavevmode\hss\endgraf
o58\ :\ Ideal\ of\ R\leavevmode\hss\endgraf
\endgroup
\penalty-1000
\par
\vskip 1 pt
\noindent

We fix the weight vector $w = (0,0,276,220,0,0,215)$ in $N$ and 
compute the initial ideal $M = in_w(I_A)$. This initial ideal 
has $44$ minimal generators.

\par
\vskip 5 pt
\begingroup
\tteight
\baselineskip=7.4pt
\lineskip=0pt
\obeyspaces
i59\ :\ \ Y\ =\ QQ[a..g,\ MonomialSize\ =>\ 16,\leavevmode\hss\endgraf
\ \ \ \ \ \ \ \ \ \ \ \ \ \ \ \ \ Weights\ =>\ {\char`\{}0,0,276,220,0,0,215{\char`\}},\leavevmode\hss\endgraf
\ \ \ \ \ \ \ \ \ \ \ \ \ \ \ \ \ Degrees\ =>transpose\ A];\leavevmode\hss\endgraf
\endgroup
\penalty-1000
\par
\vskip 1 pt
\noindent
\par
\vskip 5 pt
\begingroup
\tteight
\baselineskip=7.4pt
\lineskip=0pt
\obeyspaces
i60\ :\ IA\ =\ substitute(IA,Y);\leavevmode\hss\endgraf
\penalty-500\leavevmode\hss\endgraf
o60\ :\ Ideal\ of\ Y\leavevmode\hss\endgraf
\endgroup
\penalty-1000
\par
\vskip 1 pt
\noindent
\par
\vskip 5 pt
\begingroup
\tteight
\baselineskip=7.4pt
\lineskip=0pt
\obeyspaces
i61\ :\ M\ =\ ideal\ leadTerm\ IA\leavevmode\hss\endgraf
\penalty-500\leavevmode\hss\endgraf
\ \ \ \ \ \ \ \ \ \ \ \ \ \ 2\ 3\ \ \ \ 8\ 4\ \ \ 7\ 2\ 4\ \ \ \ \ 4\ 7\ 3\ \ \ 5\ 4\ 3\ 5\ \ \ 2\ 6\ 5\ 4\ \ \ 3\ 3\ 1 $\cdot\cdot\cdot$\leavevmode\hss\endgraf
o61\ =\ ideal\ (a\ c\ e,\ b\ g\ ,\ b\ c\ g\ ,\ a*b\ c\ f\ ,\ b\ c\ d\ f\ ,\ a\ b\ c\ g\ ,\ a\ b\ c\  $\cdot\cdot\cdot$\leavevmode\hss\endgraf
\penalty-500\leavevmode\hss\endgraf
o61\ :\ Ideal\ of\ Y\leavevmode\hss\endgraf
\endgroup
\penalty-1000
\par
\vskip 1 pt
\noindent

\begin{proposition} The three dimensional affine variety
  (\ref{uv-algebra}), for the initial ideal $M$ with respect to $w =
  (0,0,276,220,0,0,215)$ of the toric ideal of $A$ in
  (\ref{non-normal}), is not normal.
\end{proposition} 

\proof
The universal family for the toric Hilbert scheme $Hilb_A$ at $M$ is:
\begin{eqnarray*}
\{& a^2e^{15}g^{18}-z_1b^3c^6d^{10}f^{16}, \,\,
b^{13}d^{15}f^{16}-z_2a^8ce^{21}g^{14}, \,\,
c^{59} d^{57} f^{110} - z_3  e^{92} g^{134},  \\
& a c^{14} d^{11} f^{23} - z_4  b e^{19} g^{29}, \,\,
b^7 c^2 g^4 - z_5  d^4 e^3 f^6, \,\,
\ldots, \,\,
b c^{34} d^{32} f^{62} - z_{44}  e^{53} g^{76} \}.
\end{eqnarray*}
The semigroup algebra in (\ref{uv-algebra})
is generated by $44$ Laurent  monomials 
gotten from this family. It turns out that the
first four monomials suffice to generate the semigroup.
In other words, for all $j \in \{5,6,\ldots,44\}$
there exist
$ i_1,i_2,i_3, i_4 \in \N $ such that
$\,
z_{j} - 
z_1^{i_1}
z_2^{i_2}
z_3^{i_3}
z_4^{i_4} \in  J_M : (z_1 \cdots z_{44})^\infty $.
Hence the semigroup algebra in (\ref{uv-algebra}) is:
$$ \k \bigl[
 \frac{a^2 e^{15} g^{18}}{ b^3 c^6 d^{10} f^{16}}, 
 \frac{ b^{13} d^{15} f^{16} }{a^8 c e^{21} g^{14}},
 \frac{ c^{59} d^{57} f^{110}} {e^{92} g^{134}}, 
 \frac{ a c^{14} d^{11} f^{23}} {b e^{19} g^{29}}
\bigr] \,\,\, \simeq \,\,\,
\frac{\k[z_1,z_2,z_3,z_4]}{\langle z_1^5  z_2 z_3 - z_4^2 \rangle}.
$$
This algebra is not integrally closed, since
a toric hypersurface is normal if and only if
at least one of the two monomials in the defining equation 
is square-free. Its integral closure 
in $\k[ ker_\ZZ(A) ]$ is generated by the
Laurent monomial
\begin{equation}
\label{witness}
\frac{z_4}{z_1^2} \,\, = \,\, 
( z_1 z_2 z_3)^{\frac{1}{2}}  \,\, = \,\, 
\frac{ b^5 c^{26} d^{31} f^{55}}{a^3 e^{49} g^{65}}.
\end{equation}
Hence the affine chart (\ref{normal-uv-algebra}) of the toric variety
of the Gr\"obner fan of $I_A$ is the spectrum of the normal domain
$  \k[z_1,z_2,z_3,y]/ \langle
z_1  z_2 z_3 - y^2 \rangle$, 
where $y$ maps to (\ref{witness}).
\qed

\vskip .3cm 

We now examine the local equations of $Hilb_A$ about $M$ for this 
example.
\par
\vskip 5 pt
\begingroup
\tteight
\baselineskip=7.4pt
\lineskip=0pt
\obeyspaces
i62\ :\ JM\ =\ localCoherentEquations(IA)\leavevmode\hss\endgraf
\penalty-500\leavevmode\hss\endgraf
\ \ \ \ \ \ \ \ \ \ \ \ \ \ \ \ \ \ \ \ \ \ \ \ \ \ \ \ \ \ \ \ \ \ \ \ \ \ \ \ \ \ \ \ \ \ \ \ \ \ \ \ \ \ \ \ \ \ \ \ \ \ \ \ \ \ \ \ \ \  $\cdot\cdot\cdot$\leavevmode\hss\endgraf
o62\ =\ ideal\ (z\ z\ \ -\ z\ ,\ z\ z\ \ -\ z\ ,\ z\ z\ \ -\ z\ ,\ z\ z\ \ -\ z\ ,\ z\ z\ \ -\ z\ ,\ z\  $\cdot\cdot\cdot$\leavevmode\hss\endgraf
\ \ \ \ \ \ \ \ \ \ \ \ \ \ 1\ 2\ \ \ \ 3\ \ \ 1\ 2\ \ \ \ 3\ \ \ 1\ 5\ \ \ \ 4\ \ \ 1\ 3\ \ \ \ 6\ \ \ 1\ 3\ \ \ \ 6\ \ \ 1 $\cdot\cdot\cdot$\leavevmode\hss\endgraf
\penalty-500\leavevmode\hss\endgraf
o62\ :\ Ideal\ of\ B\leavevmode\hss\endgraf
\endgroup
\penalty-1000
\par
\vskip 1 pt
\noindent
\par
\vskip 5 pt
\begingroup
\tteight
\baselineskip=7.4pt
\lineskip=0pt
\obeyspaces
i63\ :\ G\ =\ removeRedundantVariables\ JM;\ \leavevmode\hss\endgraf
\penalty-500\leavevmode\hss\endgraf
o63\ :\ RingMap\ B\ <---\ B\leavevmode\hss\endgraf
\endgroup
\penalty-1000
\par
\vskip 1 pt
\noindent
\par
\vskip 5 pt
\begingroup
\tteight
\baselineskip=7.4pt
\lineskip=0pt
\obeyspaces
i64\ :\ toString\ ideal\ gens\ gb(G\ JM)\ \leavevmode\hss\endgraf
\penalty-500\leavevmode\hss\endgraf
o64\ =\ ideal(z{\char`\_}32*z{\char`\_}42{\char`\^}2*z{\char`\_}44-z{\char`\_}37{\char`\^}2*z{\char`\_}42,z{\char`\_}32{\char`\^}3*z{\char`\_}35*z{\char`\_}37{\char`\^}2-z{\char`\_}42{\char`\^}2*z{\char`\_}4 $\cdot\cdot\cdot$\leavevmode\hss\endgraf
\penalty-500\leavevmode\hss\endgraf
o64\ :\ String\leavevmode\hss\endgraf
\endgroup
\penalty-1000
\par
\vskip 1 pt
\noindent

This ideal has six generators and decomposing it 
we see that there are five components 
through the monomial ideal $M$ on this toric Hilbert scheme. They 
are defined by the ideals: 
\begin{itemize}
\item $\langle z_{32}z_{42}z_{44}-z_{37}^2,z_{32}^4z_{35}-z_{42},
z_{32}^3z_{35}z_{37}^2-z_{42}^2z_{44},
z_{32}^2z_{35}z_{37}^4-z_{42}^3z_{44}^2,\\
\indent \indent z_{32}z_{35}z_{37}^6-z_{42}^4z_{44}^3,
z_{35}z_{37}^8-z_{42}^5z_{44}^4 \rangle$ 
\item $\langle z_{44},z_{37} \rangle$
\item $\langle z_{37},z_{42}^2 \rangle$
\item $\langle z_{42},z_{35} \rangle$ 
\item $\langle z_{42},z_{32}^3 \rangle$.
\end{itemize}
All five components are three
dimensional. The first component is an affine patch of the coherent
component and two of the components are not reduced. Let $K$ be the 
first of these ideals.

\par
\vskip 5 pt
\begingroup
\tteight
\baselineskip=7.4pt
\lineskip=0pt
\obeyspaces
i65\ :\ \ K\ =\ ideal(z{\char`\_}32*z{\char`\_}42*z{\char`\_}44-z{\char`\_}37{\char`\^}2,z{\char`\_}32{\char`\^}4*z{\char`\_}35-z{\char`\_}42,\leavevmode\hss\endgraf
\ \ \ \ \ \ \ \ \ \ z{\char`\_}32{\char`\^}3*z{\char`\_}35*z{\char`\_}37{\char`\^}2-z{\char`\_}42{\char`\^}2*z{\char`\_}44,z{\char`\_}32{\char`\^}2*z{\char`\_}35*z{\char`\_}37{\char`\^}4-z{\char`\_}42{\char`\^}3*z{\char`\_}4 $\cdot\cdot\cdot$\leavevmode\hss\endgraf
\ \ \ \ \ \ \ \ \ \ z{\char`\_}32*z{\char`\_}35*z{\char`\_}37{\char`\^}6-z{\char`\_}42{\char`\^}4*z{\char`\_}44{\char`\^}3,z{\char`\_}35*z{\char`\_}37{\char`\^}8-z{\char`\_}42{\char`\^}5*z{\char`\_}44{\char`\^}4);\ \leavevmode\hss\endgraf
\penalty-500\leavevmode\hss\endgraf
o65\ :\ Ideal\ of\ B\leavevmode\hss\endgraf
\endgroup
\penalty-1000
\par
\vskip 1 pt
\noindent

Applying {\tt removeRedundantVariables} to $K$ we see that 
the affine patch of the coherent component is, locally at $M$,
a non-normal hypersurface singularity (agreeing with (\ref{witness})).
The labels on the variables depend on the order of elements in 
the initial ideal $M$ computed by \Mtwo in line {\tt i61}.

\par
\vskip 5 pt
\begingroup
\tteight
\baselineskip=7.4pt
\lineskip=0pt
\obeyspaces
i66\ :\ \ GG\ =\ removeRedundantVariables\ K;\ \leavevmode\hss\endgraf
\penalty-500\leavevmode\hss\endgraf
o66\ :\ RingMap\ B\ <---\ B\leavevmode\hss\endgraf
\endgroup
\penalty-1000
\par
\vskip 1 pt
\noindent
\par
\vskip 5 pt
\begingroup
\tteight
\baselineskip=7.4pt
\lineskip=0pt
\obeyspaces
i67\ :\ \ ideal\ gens\ gb\ (GG\ K)\ \leavevmode\hss\endgraf
\penalty-500\leavevmode\hss\endgraf
\ \ \ \ \ \ \ \ \ \ \ \ \ 5\ \ \ \ \ \ \ \ \ \ \ 2\leavevmode\hss\endgraf
o67\ =\ ideal(z\ \ z\ \ z\ \ \ -\ z\ \ )\leavevmode\hss\endgraf
\ \ \ \ \ \ \ \ \ \ \ \ \ 32\ 35\ 44\ \ \ \ 37\leavevmode\hss\endgraf
\penalty-500\leavevmode\hss\endgraf
o67\ :\ Ideal\ of\ B\leavevmode\hss\endgraf
\endgroup
\penalty-1000
\par
\vskip 1 pt
\noindent

There is a general algorithm due to de Jong \cite{DJ} for 
computing the normalization of any affine variety. 
In the toric case, the problem of normalization amounts to  
computing the minimal {\em Hilbert basis} of a given convex
rational polyhedral cone \cite{Sch}. An efficient implementation can be 
found in the software package {\tt Normaliz} by Bruns and
Koch \cite{BK}.

Our computational study of the toric Hilbert scheme in this
chapter was based on local equations rather than
global equations (arising from a projective embedding of  $Hilb_A$),
because the latter system of equations tends to be too large 
for most purposes. Nonetheless, they are interesting.
In the remainder of this section, we present a canonical 
projective embedding of the coherent component of $Hilb_A$.

Let $G_1, G_2, G_3, \ldots, G_s$ denote all the {\it Graver fibers} of
the matrix $A$. In Section 1 we saw how to compute them in \Mtwo. Each
set $G_i$ consists of the monomials in $\k[x]$ 
which have a fixed Graver degree.  Consider the set $\, {\mathbf G} \,
:= \, G_1 G_2 G_3 \cdots G_s \,$ which consists of all monomials which
are products of monomials, one from each of the distinct Graver
fibers. Let $t$ denote the cardinality of ${\mathbf G}$.  We introduce
an extra indeterminate $z$, and we consider the $\N$-graded semigroup
algebra $\,\k[z {\mathbf G}] $, which is a subalgebra of
$\k[x_1,\ldots,x_n,z]$.  Labeling the elements of ${\mathbf G}$ with
indeterminates $y_i$, We can write
$$ \k[z {\mathbf G}]  \quad = \quad
\k[y_1,y_2,\ldots,y_t]/P_A, $$
where $P_A$ is a homogeneous toric ideal
associated with a configuration of $t$ vectors in $\ZZ^{n+1}$.
We note that the torus $(\k^*)^n$ acts naturally on 
$\, \k[z {\mathbf G}]$.

\begin{example} \rm
Let $n=4,d=2$ and 
 $\, A \, = \, \left( \begin{array}{cccc}
3 & 2 & 1 & 0 \\
0 & 1 & 2 & 3 
 \end{array} \right) $,
so that $I_A$ is the ideal of the twisted cubic curve.
There are five Graver fibers:

\vskip .2cm

\par
\vskip 5 pt
\begingroup
\tteight
\baselineskip=7.4pt
\lineskip=0pt
\obeyspaces
i68\ :\ A\ =\ {\char`\{}{\char`\{}1,1,1,1{\char`\}},{\char`\{}0,1,2,3{\char`\}}{\char`\}};\leavevmode\hss\endgraf
\endgroup
\penalty-1000
\par
\vskip 1 pt
\noindent
\par
\vskip 5 pt
\begingroup
\tteight
\baselineskip=7.4pt
\lineskip=0pt
\obeyspaces
i69\ :\ I\ =\ toricIdeal\ A;\leavevmode\hss\endgraf
\penalty-500\leavevmode\hss\endgraf
o69\ :\ Ideal\ of\ R\leavevmode\hss\endgraf
\endgroup
\penalty-1000
\par
\vskip 1 pt
\noindent
\par
\vskip 5 pt
\begingroup
\tteight
\baselineskip=7.4pt
\lineskip=0pt
\obeyspaces
i70\ :\ Graver\ =\ graver\ I;\leavevmode\hss\endgraf
\penalty-500\leavevmode\hss\endgraf
\ \ \ \ \ \ \ \ \ \ \ \ \ \ 1\ \ \ \ \ \ \ 5\leavevmode\hss\endgraf
o70\ :\ Matrix\ R\ \ <---\ R\leavevmode\hss\endgraf
\endgroup
\penalty-1000
\par
\vskip 1 pt
\noindent
\par
\vskip 5 pt
\begingroup
\tteight
\baselineskip=7.4pt
\lineskip=0pt
\obeyspaces
i71\ :\ fibers\ =\ graverFibers\ Graver;\leavevmode\hss\endgraf
\endgroup
\penalty-1000
\par
\vskip 1 pt
\noindent
\par
\vskip 5 pt
\begingroup
\tteight
\baselineskip=7.4pt
\lineskip=0pt
\obeyspaces
i72\ :\ peek\ fibers\leavevmode\hss\endgraf
\penalty-500\leavevmode\hss\endgraf
o72\ =\ MutableHashTable{\char`\{}{\char`\{}2,\ 2{\char`\}}\ =>\ {\char`\{}0,\ 0{\char`\}}\ |\ ac\ b2\ |\ \ \ \ \ {\char`\}}\leavevmode\hss\endgraf
\ \ \ \ \ \ \ \ \ \ \ \ \ \ \ \ \ \ \ \ \ \ \ {\char`\{}2,\ 3{\char`\}}\ =>\ {\char`\{}0,\ 0{\char`\}}\ |\ ad\ bc\ |\leavevmode\hss\endgraf
\ \ \ \ \ \ \ \ \ \ \ \ \ \ \ \ \ \ \ \ \ \ \ {\char`\{}2,\ 4{\char`\}}\ =>\ {\char`\{}0,\ 0{\char`\}}\ |\ bd\ c2\ |\leavevmode\hss\endgraf
\ \ \ \ \ \ \ \ \ \ \ \ \ \ \ \ \ \ \ \ \ \ \ {\char`\{}3,\ 3{\char`\}}\ =>\ {\char`\{}0,\ 0{\char`\}}\ |\ a2d\ abc\ b3\ |\leavevmode\hss\endgraf
\ \ \ \ \ \ \ \ \ \ \ \ \ \ \ \ \ \ \ \ \ \ \ {\char`\{}3,\ 6{\char`\}}\ =>\ {\char`\{}0,\ 0{\char`\}}\ |\ ad2\ bcd\ c3\ |\leavevmode\hss\endgraf
\penalty-500\leavevmode\hss\endgraf
o72\ :\ Net\leavevmode\hss\endgraf
\endgroup
\penalty-1000
\par
\vskip 1 pt
\noindent
\vskip .2cm

The set ${\mathbf G} = G_1 G_2 G_3 G_4 G_5 \,$ consists of
$22$ monomials of degree $14$.

\vskip .2cm

\par
\vskip 5 pt
\begingroup
\tteight
\baselineskip=7.4pt
\lineskip=0pt
\obeyspaces
i73\ :\ G\ =\ trim\ product(values\ fibers,\ ideal)\leavevmode\hss\endgraf
\penalty-500\leavevmode\hss\endgraf
\ \ \ \ \ \ \ \ \ \ \ \ \ \ 5\ \ \ \ \ 5\ \ \ 4\ 3\ 5\ \ \ 5\ 3\ 4\ \ \ 4\ 2\ 2\ 4\ \ \ 3\ 4\ \ \ 4\ \ \ 2\ 6\ 4\ \ \ 4\  $\cdot\cdot\cdot$\leavevmode\hss\endgraf
o73\ =\ ideal\ (a\ b*c*d\ ,\ a\ b\ d\ ,\ a\ c\ d\ ,\ a\ b\ c\ d\ ,\ a\ b\ c*d\ ,\ a\ b\ d\ ,\ a\ b $\cdot\cdot\cdot$\leavevmode\hss\endgraf
\penalty-500\leavevmode\hss\endgraf
o73\ :\ Ideal\ of\ R\leavevmode\hss\endgraf
\endgroup
\penalty-1000
\par
\vskip 1 pt
\noindent
\par
\vskip 5 pt
\begingroup
\tteight
\baselineskip=7.4pt
\lineskip=0pt
\obeyspaces
i74\ :\ numgens\ G\leavevmode\hss\endgraf
\penalty-500\leavevmode\hss\endgraf
o74\ =\ 22\leavevmode\hss\endgraf
\endgroup
\penalty-1000
\par
\vskip 1 pt
\noindent

\vskip .2cm

We introduce a polynomial ring in $22$ variables
$y_1,y_2,\ldots,y_{22}$, and we compute the ideal $P_A$.
It is generated by $180$ binomial quadrics. 

\vskip .2cm
\par
\vskip 5 pt
\begingroup
\tteight
\baselineskip=7.4pt
\lineskip=0pt
\obeyspaces
i75\ :\ z\ =\ symbol\ z;\leavevmode\hss\endgraf
\endgroup
\penalty-1000
\par
\vskip 1 pt
\noindent
\par
\vskip 5 pt
\begingroup
\tteight
\baselineskip=7.4pt
\lineskip=0pt
\obeyspaces
i76\ :\ S\ =\ QQ[a,b,c,d,z];\leavevmode\hss\endgraf
\endgroup
\penalty-1000
\par
\vskip 1 pt
\noindent
\par
\vskip 5 pt
\begingroup
\tteight
\baselineskip=7.4pt
\lineskip=0pt
\obeyspaces
i77\ :\ zG\ =\ z\ **\ substitute(gens\ G,\ S);\leavevmode\hss\endgraf
\penalty-500\leavevmode\hss\endgraf
\ \ \ \ \ \ \ \ \ \ \ \ \ \ 1\ \ \ \ \ \ \ 22\leavevmode\hss\endgraf
o77\ :\ Matrix\ S\ \ <---\ S\leavevmode\hss\endgraf
\endgroup
\penalty-1000
\par
\vskip 1 pt
\noindent
\par
\vskip 5 pt
\begingroup
\tteight
\baselineskip=7.4pt
\lineskip=0pt
\obeyspaces
i78\ :\ R\ =\ QQ[y{\char`\_}1\ ..\ y{\char`\_}22];\leavevmode\hss\endgraf
\endgroup
\penalty-1000
\par
\vskip 1 pt
\noindent
\par
\vskip 5 pt
\begingroup
\tteight
\baselineskip=7.4pt
\lineskip=0pt
\obeyspaces
i79\ :\ F\ =\ map(S,R,zG)\leavevmode\hss\endgraf
\penalty-500\leavevmode\hss\endgraf
\ \ \ \ \ \ \ \ \ \ \ \ \ \ \ \ 5\ \ \ \ \ 5\ \ \ \ 4\ 3\ 5\ \ \ \ 5\ 3\ 4\ \ \ \ 4\ 2\ 2\ 4\ \ \ \ 3\ 4\ \ \ 4\ \ \ \ 2\ 6 $\cdot\cdot\cdot$\leavevmode\hss\endgraf
o79\ =\ map(S,R,{\char`\{}a\ b*c*d\ z,\ a\ b\ d\ z,\ a\ c\ d\ z,\ a\ b\ c\ d\ z,\ a\ b\ c*d\ z,\ a\ b\  $\cdot\cdot\cdot$\leavevmode\hss\endgraf
\penalty-500\leavevmode\hss\endgraf
o79\ :\ RingMap\ S\ <---\ R\leavevmode\hss\endgraf
\endgroup
\penalty-1000
\par
\vskip 1 pt
\noindent
\par
\vskip 5 pt
\begingroup
\tteight
\baselineskip=7.4pt
\lineskip=0pt
\obeyspaces
i80\ :\ PA\ =\ trim\ ker\ F\leavevmode\hss\endgraf
\penalty-500\leavevmode\hss\endgraf
\ \ \ \ \ \ \ \ \ \ \ \ \ \ 2\ \ \ \ \ \ \ \ \ \ \ \ \ \ \ \ \ \ \ \ \ \ \ \ \ \ \ \ \ \ \ \ \ \ \ \ \ \ \ \ \ \ \ \ \ \ \ \ \ \ \ \ \ \ \  $\cdot\cdot\cdot$\leavevmode\hss\endgraf
o80\ =\ ideal\ (y\ \ \ -\ y\ \ y\ \ ,\ y\ \ y\ \ \ -\ y\ \ y\ \ ,\ y\ \ y\ \ \ -\ y\ \ y\ \ ,\ y\ \ y\ \ \ -\  $\cdot\cdot\cdot$\leavevmode\hss\endgraf
\ \ \ \ \ \ \ \ \ \ \ \ \ \ 21\ \ \ \ 20\ 22\ \ \ 19\ 21\ \ \ \ 18\ 22\ \ \ 18\ 21\ \ \ \ 17\ 22\ \ \ 17\ 21\ \ \  $\cdot\cdot\cdot$\leavevmode\hss\endgraf
\penalty-500\leavevmode\hss\endgraf
o80\ :\ Ideal\ of\ R\leavevmode\hss\endgraf
\endgroup
\penalty-1000
\par
\vskip 1 pt
\noindent

\vskip .2cm

These equations define a toric surface
of degree $30$ in projective $21$-space.
\par
\vskip 5 pt
\begingroup
\tteight
\baselineskip=7.4pt
\lineskip=0pt
\obeyspaces
i81\ :\ codim\ PA\leavevmode\hss\endgraf
\penalty-500\leavevmode\hss\endgraf
o81\ =\ 19\leavevmode\hss\endgraf
\endgroup
\penalty-1000
\par
\vskip 1 pt
\noindent
\par
\vskip 5 pt
\begingroup
\tteight
\baselineskip=7.4pt
\lineskip=0pt
\obeyspaces
i82\ :\ degree\ PA\leavevmode\hss\endgraf
\penalty-500\leavevmode\hss\endgraf
o82\ =\ 30\leavevmode\hss\endgraf
\endgroup
\penalty-1000
\par
\vskip 1 pt
\noindent

The surface is smooth, but there are too many equations and the
codimension is too large to use the Jacobian criterion for smoothness
\cite[\S 16.6]{Eis} directly. Instead we check smoothness for each
open set $y_i \neq 0$. 

\par
\vskip 5 pt
\begingroup
\tteight
\baselineskip=7.4pt
\lineskip=0pt
\obeyspaces
i83\ :\ Aff\ =\ apply(1..22,\ v\ ->\ (\leavevmode\hss\endgraf
\ \ \ \ \ \ \ \ \ \ \ \ \ \ \ \ \ \ \ \ \ \ \ \ \ \ \ \ \ K\ =\ substitute(PA,y{\char`\_}v\ =>\ 1);\leavevmode\hss\endgraf
\ \ \ \ \ \ \ \ \ \ \ \ \ \ \ \ \ \ \ \ \ \ \ \ \ \ \ \ \ FF\ =\ removeRedundantVariables\ K;\leavevmode\hss\endgraf
\ \ \ \ \ \ \ \ \ \ \ \ \ \ \ \ \ \ \ \ \ \ \ \ \ \ \ \ \ ideal\ gens\ gb\ (FF\ K)));\ \leavevmode\hss\endgraf
\endgroup
\penalty-1000
\par
\vskip 1 pt
\noindent
\par
\vskip 5 pt
\begingroup
\tteight
\baselineskip=7.4pt
\lineskip=0pt
\obeyspaces
i84\ :\ scan(Aff,\ i\ ->\ print\ toString\ i);\leavevmode\hss\endgraf
ideal()\leavevmode\hss\endgraf
ideal()\leavevmode\hss\endgraf
ideal()\leavevmode\hss\endgraf
ideal(y{\char`\_}1{\char`\^}4*y{\char`\_}5*y{\char`\_}21-1)\leavevmode\hss\endgraf
ideal(y{\char`\_}1{\char`\^}4*y{\char`\_}6{\char`\^}6*y{\char`\_}21-1)\leavevmode\hss\endgraf
ideal()\leavevmode\hss\endgraf
ideal(y{\char`\_}1{\char`\^}2*y{\char`\_}11{\char`\^}2*y{\char`\_}17-1)\leavevmode\hss\endgraf
ideal(y{\char`\_}1{\char`\^}3*y{\char`\_}9{\char`\^}2*y{\char`\_}21{\char`\^}2-1)\leavevmode\hss\endgraf
ideal(y{\char`\_}6{\char`\^}3*y{\char`\_}21-y{\char`\_}10,y{\char`\_}1*y{\char`\_}10{\char`\^}3-y{\char`\_}6{\char`\^}2,y{\char`\_}1*y{\char`\_}6*y{\char`\_}10{\char`\^}2*y{\char`\_}21-1)\leavevmode\hss\endgraf
ideal(y{\char`\_}6*y{\char`\_}15-1,y{\char`\_}2*y{\char`\_}15{\char`\^}2-y{\char`\_}6*y{\char`\_}14,y{\char`\_}6{\char`\^}2*y{\char`\_}14-y{\char`\_}2*y{\char`\_}15)\leavevmode\hss\endgraf
ideal()\leavevmode\hss\endgraf
ideal(y{\char`\_}11*y{\char`\_}13-1,y{\char`\_}1{\char`\^}2*y{\char`\_}21{\char`\^}3-y{\char`\_}13{\char`\^}2)\leavevmode\hss\endgraf
ideal(y{\char`\_}1{\char`\^}2*y{\char`\_}14{\char`\^}3*y{\char`\_}21{\char`\^}3-1)\leavevmode\hss\endgraf
ideal(y{\char`\_}10{\char`\^}2*y{\char`\_}21-1,y{\char`\_}1*y{\char`\_}15{\char`\^}4-y{\char`\_}10{\char`\^}3)\leavevmode\hss\endgraf
ideal()\leavevmode\hss\endgraf
ideal(y{\char`\_}11*y{\char`\_}20-1,y{\char`\_}3*y{\char`\_}20{\char`\^}2-y{\char`\_}11*y{\char`\_}17,y{\char`\_}11{\char`\^}2*y{\char`\_}17-y{\char`\_}3*y{\char`\_}20)\leavevmode\hss\endgraf
ideal(y{\char`\_}11*y{\char`\_}18*y{\char`\_}21-1,y{\char`\_}1*y{\char`\_}21{\char`\^}3-y{\char`\_}11*y{\char`\_}18{\char`\^}2,y{\char`\_}11{\char`\^}2*y{\char`\_}18{\char`\^}3-y{\char`\_}1*y{\char`\_}21{\char`\^}2 $\cdot\cdot\cdot$\leavevmode\hss\endgraf
ideal(y{\char`\_}1*y{\char`\_}19{\char`\^}4*y{\char`\_}21{\char`\^}4-1)\leavevmode\hss\endgraf
ideal(y{\char`\_}15*y{\char`\_}22-1)\leavevmode\hss\endgraf
ideal()\leavevmode\hss\endgraf
ideal(y{\char`\_}20*y{\char`\_}22-1)\leavevmode\hss\endgraf
ideal()\leavevmode\hss\endgraf
\endgroup
\penalty-1000
\par
\vskip 1 pt
\noindent

By examining these local equations, we see that $Hilb_A$ is smooth, and also
that there are eight fixed points under the
action of the $2$-dimensional torus. They correspond
to the variables $y_1,y_2,y_3,y_6,y_{11},y_{15},y_{20}$ and $y_{22}$. 
By setting any of these eight variables to $1$ in the
$180$ quadrics above, we obtain an affine variety
isomorphic to the affine plane.
\end{example}

\begin{theorem} \label{isomorphism}
The coherent component of the toric Hilbert scheme $Hilb_A$ is
isomorphic to $Proj \k[z {\mathbf G}]$.
\end{theorem}

We present a sketch of the proof of this theorem. The first
step is to define a morphism from $\,Hilb_A \,$ to
$\P({\mathbf G}) = Proj \, \k[y_1,y_2,\ldots,y_t]$.
Consider any point $I$ on $Hilb_A$. The component 
$I \cap \k G_i$ of $I$ in the $i$-th Graver degree
is a linear subspace of codimension $1$ in $\k G_i$.
We can represent this hyperplane by a polynomial
$\, g_i(I) = \sum_{u \in G_i } c_u x^u \,$ which is defined
uniquely up to scaling. By taking the product of these
polynomials, we obtain a unique (up to scaling)
polynomial which is supported on ${\mathbf G} = G_1 G_2 \cdots G_t$.
The assignment $\, I \mapsto g_1(I) g_2(I) \cdots g_t(I)\,$
defines a morphism from $Hilb_A \,$ to $\P({\mathbf G})$.
Now consider the restriction of this morphism to the
coherent component. It is clearly an isomorphism on the
$(\k^*)^n$-orbit of the toric ideal $I_A$. To show that
it is an isomorphism on the entire coherent component,
we consider the affine chart around an initial monomial ideal
$M = in_w(I_A)$. It corresponds to the vertex in direction $w$
of the Minkowski sum of the Newton polytopes of the
polynomials $\,g_1(I_A),  g_2(I_A), \ldots, g_t(I_A)$.
The normal cone at that vertex coincides with the
cone $\sigma$ of the Gr\"obner fan  $\Sigma_A$
which has $w$ in its interior \cite[\S 3]{St2}.
The restriction of our morphism to the affine chart around $M$
of the coherent component,  as described in equation
(\ref{uv-algebra}), shows that 
this restriction is an isomorphism. The reason for this
is that each pair of vectors $\{u_i, v_i\}$
seen in the reduced Gr\"obner basis lies in one of the
Graver fibers $G_j$. This concludes our sketch of the
proof of Theorem~\ref{isomorphism}. \qed

\section{Appendix 1: Fourier-Motzkin Elimination}

We now give the \Mtwo code for converting the generator/inequality
representation of a rational convex polyhedron to the other. It is
based on the Fourier-Motzkin elimination procedure for eliminating a
variable from a system of inequalities \cite {Zie}. This code was
written by Greg Smith.

Given any cone $C \subset \R^d$, the polar cone of $C$ is defined to be
$$C^{\vee} = \{ x \in \R^d \mid x \cdot y \leq 0, \mbox{for all\ } y
\in C\}.$$

\noindent For a $d \times n$ matrix $Z$, define
$cone(Z) = \{ Z x \mid x \in \R_{\geq 0}^n \} \subset \R^d,$ and 
$\mathit{affine}(Z) = \{ Z x \mid x \in \R^n \} \subset \R^d.$
For two integer matrices $Z$ and $H$, both having  $d$
rows, {\tt polarCone(Z,H)} returns a list of two integer matrices
{\tt\{A,E\}} such that $$cone(Z) + \mathit{affine}(H) = \{ x \in \R^d \mid A^t
x \leq 0, E^t x = 0\}.$$ 

Equivalently, $(cone(Z) + \mathit{affine}(H))^\vee = cone(A) + \mathit{affine}(E).$

We now describe each routine in the package {\tt polarCone.m2}.  We have
simplified the code for readability, sometimes at the cost of efficiency.
We start with three simple subroutines: {\tt primitive}, {\tt toZZ}, and {\tt
rotateMatrix}. 

\medskip
The routine {\tt primitive} takes a list of integers {\tt L}, and divides
each element of this list by their greatest common denominator.

\par
\vskip 5 pt
\begingroup
\tteight
\baselineskip=7.4pt
\lineskip=0pt
\obeyspaces
i85\ :\ primitive\ =\ (L)\ ->\ (\leavevmode\hss\endgraf
\ \ \ \ \ \ \ \ \ \ \ n\ :=\ {\char`\#}L-1;\ g\ :=\ L{\char`\#}n;\leavevmode\hss\endgraf
\ \ \ \ \ \ \ \ \ \ \ while\ n\ >\ 0\ do\ (n\ =\ n-1;\ g\ =\ gcd(g,\ L{\char`\#}n);\leavevmode\hss\endgraf
\ \ \ \ \ \ \ \ \ \ \ \ \ \ \ \ if\ g\ ===\ 1\ then\ n\ =\ 0);\leavevmode\hss\endgraf
\ \ \ \ \ \ \ \ \ \ \ if\ g\ ===\ 1\ then\ L\ else\ apply(L,\ i\ ->\ i\ //\ g));\leavevmode\hss\endgraf
\endgroup
\penalty-1000
\par
\vskip 1 pt
\noindent

\medskip
The routine {\tt toZZ} converts a list of rational numbers to a list of 
integers, by multiplying by their common denominator.
\par
\vskip 5 pt
\begingroup
\tteight
\baselineskip=7.4pt
\lineskip=0pt
\obeyspaces
i86\ :\ toZZ\ =\ (L)\ ->\ (\leavevmode\hss\endgraf
\ \ \ \ \ \ \ \ \ \ \ d\ :=\ apply(L,\ e\ ->\ denominator\ e);\leavevmode\hss\endgraf
\ \ \ \ \ \ \ \ \ \ \ R\ :=\ ring\ d{\char`\#}0;\ \ \ \ \ \ \ \ \ \ \ \ \ l\ :=\ 1{\char`\_}R;\leavevmode\hss\endgraf
\ \ \ \ \ \ \ \ \ \ \ scan(d,\ i\ ->\ (l\ =\ (l*i\ //\ gcd(l,i))));\ \ \ \ \leavevmode\hss\endgraf
\ \ \ \ \ \ \ \ \ \ \ apply(L,\ e\ ->\ (numerator(l*e))));\leavevmode\hss\endgraf
\endgroup
\penalty-1000
\par
\vskip 1 pt
\noindent

\medskip
The routine {\tt rotateMatrix} is a kind of transpose.  Its input is a
matrix, and its output is a matrix of the same shape as the transpose.
It places the matrix in the form so that in the routine {\tt polarCone},
computing a Gr\"obner basis will do the Gaussian elimination that is needed.
\par
\vskip 5 pt
\begingroup
\tteight
\baselineskip=7.4pt
\lineskip=0pt
\obeyspaces
i87\ :\ rotateMatrix\ :=\ (M)\ ->\ (\leavevmode\hss\endgraf
\ \ \ \ \ \ \ \ \ \ \ r\ :=\ rank\ source\ M;\ \ \ \ \ \ \ \ c\ :=\ rank\ target\ M;\leavevmode\hss\endgraf
\ \ \ \ \ \ \ \ \ \ \ matrix\ table(r,\ c,\ (i,j)\ ->\ M{\char`\_}(c-j-1,\ r-i-1)));\leavevmode\hss\endgraf
\endgroup
\penalty-1000
\par
\vskip 1 pt
\noindent

\medskip
The procedure of Fourier-Motzkin elimination as presented by 
Ziegler in \cite{Zie} is used, together with some heuristics which he
presents as exercises.  The following, which is a kind of $S$-pair
criterion for inequalities, comes from exercise 2.15(i) in \cite{Zie}.

The routine {\tt isRedundant} determines if a row vector (inequality)
is redundant. Its input argument {\tt V} is the same input that is
used in {\tt fourierMotzkin}: it is a list of sets of integers.  Each
entry contains indices of the original rays which do {\sl not} vanish
at the corresponding row vector.  {\tt vert} is a set of integers; the
original rays for the row vector in question.  A boolean value is
returned.  

\par
\vskip 5 pt
\begingroup
\tteight
\baselineskip=7.4pt
\lineskip=0pt
\obeyspaces
i88\ :\ isRedundant\ =\ (V,\ vert)\ ->\ (\ \leavevmode\hss\endgraf
\ \ \ \ \ \ \ \ \ \ \ \ \ --\ the\ row\ vector\ is\ redundant\ iff\ 'vert'\ contains\ an\ \leavevmode\hss\endgraf
\ \ \ \ \ \ \ \ \ \ \ \ \ --\ entry\ in\ 'V'.\ \ \leavevmode\hss\endgraf
\ \ \ \ \ \ \ \ \ \ \ \ \ x\ :=\ 0;\ k\ :=\ 0;\leavevmode\hss\endgraf
\ \ \ \ \ \ \ \ \ \ \ \ \ numRow\ :=\ {\char`\#}V;\ --\ equals\ the\ number\ of\ inequalities\ \leavevmode\hss\endgraf
\ \ \ \ \ \ \ \ \ \ \ \ \ while\ x\ <\ 1\ and\ k\ <\ numRow\ do\ (\ \leavevmode\hss\endgraf
\ \ \ \ \ \ \ \ \ \ \ \ \ \ \ \ \ \ if\ isSubset(V{\char`\#}k,\ vert)\ then\ x\ =\ x+1;\ \leavevmode\hss\endgraf
\ \ \ \ \ \ \ \ \ \ \ \ \ \ \ \ \ \ k\ =\ k+1;);\ \leavevmode\hss\endgraf
\ \ \ \ \ \ \ \ \ \ \ \ \ x\ ===\ 1);\leavevmode\hss\endgraf
\endgroup
\penalty-1000
\par
\vskip 1 pt
\noindent

\medskip
The main work horse of {\tt polarCone.m2} is the subroutine 
{\tt fourierMotzkin} which eliminates the first variable in the
inequalities {\tt A} using the double description version of
Fourier-Motzkin elimination. The set {\tt A} is a list of lists of
integers, each entry corresponding to a row vector in the system of
inequalities.  The argument {\tt V} is a list of sets of integers.
Each entry contains the indices of the original rays which do {\sl
  not} vanish at the corresponding row vector in {\tt A}.  Note that
this set is the {\sl complement} of the set $V_i$ appearing in
exercise 2.15 in \cite{Zie}. The argument {\tt spot} is the integer
index of the variable being eliminated.  

The routine returns a list {\tt \{projA,projV\}} where {\tt projA} is
a list of lists of integers.  Each entry corresponds to a row vector
in the projected system of inequalities.  The list {\tt projV} is a
list of sets of integers.  Each entry contains indices of the original
rays which do {\sl not} vanish at the corresponding row vector in {\tt
  projA}. 

\par
\vskip 5 pt
\begingroup
\tteight
\baselineskip=7.4pt
\lineskip=0pt
\obeyspaces
i89\ :\ fourierMotzkin\ :=\ (A,\ V,\ spot)\ ->\ (\leavevmode\hss\endgraf
\ \ \ \ \ \ \ \ \ \ \ --\ initializing\ local\ variables\leavevmode\hss\endgraf
\ \ \ \ \ \ \ \ \ \ \ numRow\ :=\ {\char`\#}A;\ \ \ \ \ \ \ \ \ \ \ \ \ \ \ --\ equal\ to\ the\ length\ of\ V\leavevmode\hss\endgraf
\ \ \ \ \ \ \ \ \ \ \ numCol\ :=\ {\char`\#}(A{\char`\#}0);\ \ \ \ \ \ \ \ \ \ \ pos\ :=\ {\char`\{}{\char`\}};\ \ \ \ \ \ \ \leavevmode\hss\endgraf
\ \ \ \ \ \ \ \ \ \ \ neg\ :=\ {\char`\{}{\char`\}};\ \ \ \ \ \ \ \ \ \ \ \ \ \ \ \ \ \ projA\ :=\ {\char`\{}{\char`\}};\ \ \ \ \ \leavevmode\hss\endgraf
\ \ \ \ \ \ \ \ \ \ \ projV\ :=\ {\char`\{}{\char`\}};\ \ \ \ \ \ \ \ \ \ \ \ \ \ \ \ k\ :=\ 0;\leavevmode\hss\endgraf
\ \ \ \ \ \ \ \ \ \ \ --\ divide\ the\ inequalities\ into\ three\ groups.\leavevmode\hss\endgraf
\ \ \ \ \ \ \ \ \ \ \ while\ k\ <\ numRow\ do\ (\leavevmode\hss\endgraf
\ \ \ \ \ \ \ \ \ \ \ \ \ \ \ \ if\ A{\char`\#}k{\char`\#}0\ <\ 0\ then\ neg\ =\ append(neg,\ k)\leavevmode\hss\endgraf
\ \ \ \ \ \ \ \ \ \ \ \ \ \ \ \ else\ if\ A{\char`\#}k{\char`\#}0\ >\ 0\ then\ pos\ =\ append(pos,\ k)\leavevmode\hss\endgraf
\ \ \ \ \ \ \ \ \ \ \ \ \ \ \ \ else\ (projA\ =\ append(projA,\ A{\char`\#}k);\leavevmode\hss\endgraf
\ \ \ \ \ \ \ \ \ \ \ \ \ \ \ \ \ \ \ \ \ projV\ =\ append(projV,\ V{\char`\#}k););\leavevmode\hss\endgraf
\ \ \ \ \ \ \ \ \ \ \ \ \ \ \ \ k\ =\ k+1;);\ \ \ \ \ \ \leavevmode\hss\endgraf
\ \ \ \ \ \ \ \ \ \ \ --\ generate\ new\ irredundant\ inequalities.\leavevmode\hss\endgraf
\ \ \ \ \ \ \ \ \ \ \ scan(pos,\ i\ ->\ scan(neg,\ j\ ->\ (vert\ :=\ V{\char`\#}i\ +\ V{\char`\#}j;\leavevmode\hss\endgraf
\ \ \ \ \ \ \ \ \ \ \ \ \ \ \ \ \ \ \ \ \ \ \ \ \ \ if\ not\ isRedundant(projV,\ vert)\ \ \leavevmode\hss\endgraf
\ \ \ \ \ \ \ \ \ \ \ \ \ \ \ \ \ \ \ \ \ \ \ \ \ \ then\ (iRow\ :=\ A{\char`\#}i;\ \ \ \ \ jRow\ :=\ A{\char`\#}j;\leavevmode\hss\endgraf
\ \ \ \ \ \ \ \ \ \ \ \ \ \ \ \ \ \ \ \ \ \ \ \ \ \ \ \ \ \ \ iCoeff\ :=\ -\ jRow{\char`\#}0;\leavevmode\hss\endgraf
\ \ \ \ \ \ \ \ \ \ \ \ \ \ \ \ \ \ \ \ \ \ \ \ \ \ \ \ \ \ \ jCoeff\ :=\ iRow{\char`\#}0;\leavevmode\hss\endgraf
\ \ \ \ \ \ \ \ \ \ \ \ \ \ \ \ \ \ \ \ \ \ \ \ \ \ \ \ \ \ \ a\ :=\ iCoeff*iRow\ +\ jCoeff*jRow;\leavevmode\hss\endgraf
\ \ \ \ \ \ \ \ \ \ \ \ \ \ \ \ \ \ \ \ \ \ \ \ \ \ \ \ \ \ \ projA\ =\ append(projA,\ a);\leavevmode\hss\endgraf
\ \ \ \ \ \ \ \ \ \ \ \ \ \ \ \ \ \ \ \ \ \ \ \ \ \ \ \ \ \ \ projV\ =\ append(projV,\ vert););)));\leavevmode\hss\endgraf
\ \ \ \ \ \ \ \ \ \ \ --\ don't\ forget\ the\ implicit\ inequalities\ '-t\ <=\ 0'.\leavevmode\hss\endgraf
\ \ \ \ \ \ \ \ \ \ \ scan(pos,\ i\ ->\ (vert\ :=\ V{\char`\#}i\ +\ set{\char`\{}spot{\char`\}};\leavevmode\hss\endgraf
\ \ \ \ \ \ \ \ \ \ \ \ \ \ \ \ if\ not\ isRedundant(projV,\ vert)\ then\ (\leavevmode\hss\endgraf
\ \ \ \ \ \ \ \ \ \ \ \ \ \ \ \ \ \ \ \ \ projA\ =\ append(projA,\ A{\char`\#}i);\leavevmode\hss\endgraf
\ \ \ \ \ \ \ \ \ \ \ \ \ \ \ \ \ \ \ \ \ projV\ =\ append(projV,\ vert););));\leavevmode\hss\endgraf
\ \ \ \ \ \ \ \ \ \ \ --\ remove\ the\ first\ column\ \leavevmode\hss\endgraf
\ \ \ \ \ \ \ \ \ \ \ projA\ =\ apply(projA,\ e\ ->\ e{\char`\_}{\char`\{}1..(numCol-1){\char`\}});\leavevmode\hss\endgraf
\ \ \ \ \ \ \ \ \ \ \ {\char`\{}projA,\ projV{\char`\}});\ \ \ \leavevmode\hss\endgraf
\endgroup
\penalty-1000
\par
\vskip 1 pt
\noindent

\medskip
As mentioned above, {\tt polarCone} takes two matrices {\tt Z, H},
both having $d$ rows, and outputs a pair of matrices {\tt A, E} 
such that $(cone(Z) + \mathit{affine}(H))^\vee = cone(A) + \mathit{affine}(E).$\\

\par
\vskip 5 pt
\begingroup
\tteight
\baselineskip=7.4pt
\lineskip=0pt
\obeyspaces
i90\ :\ polarCone(Matrix,\ Matrix)\ :=\ (Z,\ H)\ ->\ (\leavevmode\hss\endgraf
\ \ \ \ \ \ \ \ \ \ \ R\ :=\ ring\ source\ Z;\leavevmode\hss\endgraf
\ \ \ \ \ \ \ \ \ \ \ if\ R\ =!=\ ring\ source\ H\ then\ error\ ("polarCone:\ "\ |\ \leavevmode\hss\endgraf
\ \ \ \ \ \ \ \ \ \ \ \ \ \ \ \ "expected\ matrices\ over\ the\ same\ ring");\leavevmode\hss\endgraf
\ \ \ \ \ \ \ \ \ \ \ if\ rank\ target\ Z\ =!=\ rank\ target\ H\ then\ error\ (\leavevmode\hss\endgraf
\ \ \ \ \ \ \ \ \ \ \ \ \ \ \ \ "polarCone:\ expected\ matrices\ to\ have\ the\ "\ |\leavevmode\hss\endgraf
\ \ \ \ \ \ \ \ \ \ \ \ \ \ \ \ "same\ number\ of\ rows");\ \ \ \ \ \leavevmode\hss\endgraf
\ \ \ \ \ \ \ \ \ \ \ if\ (R\ =!=\ ZZ)\ then\ error\ ("polarCone:\ expected\ "\ |\ \leavevmode\hss\endgraf
\ \ \ \ \ \ \ \ \ \ \ \ \ \ \ \ "matrices\ over\ 'ZZ'");\leavevmode\hss\endgraf
\ \ \ \ \ \ \ \ \ \ \ --\ expressing\ 'cone(Y)+affine(B)'\ as\ '{\char`\{}x\ :\ Ax\ <=\ 0{\char`\}}'\leavevmode\hss\endgraf
\ \ \ \ \ \ \ \ \ \ \ Y\ :=\ substitute(Z,\ QQ);\ \ \ \ \ B\ :=\ substitute(H,\ QQ);\ \ \ \leavevmode\hss\endgraf
\ \ \ \ \ \ \ \ \ \ \ if\ rank\ source\ B\ >\ 0\ then\ Y\ =\ Y\ |\ B\ |\ -B;\leavevmode\hss\endgraf
\ \ \ \ \ \ \ \ \ \ \ n\ :=\ rank\ source\ Y;\ \ \ \ \ \ \ \ \ d\ :=\ rank\ target\ Y;\ \ \ \ \ \leavevmode\hss\endgraf
\ \ \ \ \ \ \ \ \ \ \ A\ :=\ Y\ |\ -id{\char`\_}(QQ{\char`\^}d);\leavevmode\hss\endgraf
\ \ \ \ \ \ \ \ \ \ \ --\ computing\ the\ row\ echelon\ form\ of\ 'A'\leavevmode\hss\endgraf
\ \ \ \ \ \ \ \ \ \ \ A\ =\ gens\ gb\ rotateMatrix\ A;\leavevmode\hss\endgraf
\ \ \ \ \ \ \ \ \ \ \ L\ :=\ rotateMatrix\ leadTerm\ A;\leavevmode\hss\endgraf
\ \ \ \ \ \ \ \ \ \ \ A\ =\ rotateMatrix\ A;\leavevmode\hss\endgraf
\ \ \ \ \ \ \ \ \ \ \ --\ find\ pivots\leavevmode\hss\endgraf
\ \ \ \ \ \ \ \ \ \ \ numRow\ =\ rank\ target\ A;\ \ \ \ \ \ \ \ \ \ \ \ \ \ \ \ \ \ --\ numRow\ <=\ d\leavevmode\hss\endgraf
\ \ \ \ \ \ \ \ \ \ \ i\ :=\ 0;\ \ \ \ \ \ \ \ \ \ \ \ \ \ \ \ \ \ \ \ \ pivotCol\ :=\ {\char`\{}{\char`\}};\leavevmode\hss\endgraf
\ \ \ \ \ \ \ \ \ \ \ while\ i\ <\ numRow\ do\ (j\ :=\ 0;\leavevmode\hss\endgraf
\ \ \ \ \ \ \ \ \ \ \ \ \ \ \ \ while\ j\ <\ n+d\ and\ L{\char`\_}(i,j)\ =!=\ 1{\char`\_}QQ\ do\ j\ =\ j+1;\leavevmode\hss\endgraf
\ \ \ \ \ \ \ \ \ \ \ \ \ \ \ \ pivotCol\ =\ append(pivotCol,\ j);\leavevmode\hss\endgraf
\ \ \ \ \ \ \ \ \ \ \ \ \ \ \ \ i\ =\ i+1;);\leavevmode\hss\endgraf
\ \ \ \ \ \ \ \ \ \ \ --\ computing\ the\ row-reduced\ echelon\ form\ of\ 'A'\leavevmode\hss\endgraf
\ \ \ \ \ \ \ \ \ \ \ A\ =\ ((submatrix(A,\ pivotCol)){\char`\^}(-1))\ *\ A;\leavevmode\hss\endgraf
\ \ \ \ \ \ \ \ \ \ \ --\ converting\ 'A'\ into\ a\ list\ of\ integer\ row\ vectors\ \leavevmode\hss\endgraf
\ \ \ \ \ \ \ \ \ \ \ A\ =\ entries\ A;\leavevmode\hss\endgraf
\ \ \ \ \ \ \ \ \ \ \ A\ =\ apply(A,\ e\ ->\ primitive\ toZZ\ e);\leavevmode\hss\endgraf
\ \ \ \ \ \ \ \ \ \ \ --\ creating\ the\ vertex\ list\ 'V'\ for\ double\ description\leavevmode\hss\endgraf
\ \ \ \ \ \ \ \ \ \ \ --\ and\ listing\ the\ variables\ 'T'\ which\ remain\ to\ be\leavevmode\hss\endgraf
\ \ \ \ \ \ \ \ \ \ \ --\ eliminated\leavevmode\hss\endgraf
\ \ \ \ \ \ \ \ \ \ \ V\ :=\ {\char`\{}{\char`\}};\ \ \ \ \ \ \ \ \ \ \ \ \ \ \ \ \ \ \ \ T\ :=\ toList(0..(n-1));\leavevmode\hss\endgraf
\ \ \ \ \ \ \ \ \ \ \ scan(pivotCol,\ e\ ->\ (if\ e\ <\ n\ then\ (T\ =\ delete(e,\ T);\leavevmode\hss\endgraf
\ \ \ \ \ \ \ \ \ \ \ \ \ \ \ \ \ \ \ \ \ \ \ \ \ \ V\ =\ append(V,\ set{\char`\{}e{\char`\}});)));\leavevmode\hss\endgraf
\ \ \ \ \ \ \ \ \ \ \ --\ separating\ inequalities\ 'A'\ and\ equalities\ 'E'\leavevmode\hss\endgraf
\ \ \ \ \ \ \ \ \ \ \ eqnRow\ :=\ {\char`\{}{\char`\}};\ \ \ \ \ \ \ \ \ \ \ \ \ \ \ ineqnRow\ :=\ {\char`\{}{\char`\}};\leavevmode\hss\endgraf
\ \ \ \ \ \ \ \ \ \ \ scan(numRow,\ i\ ->\ (if\ pivotCol{\char`\#}i\ >=\ n\ then\ \leavevmode\hss\endgraf
\ \ \ \ \ \ \ \ \ \ \ \ \ \ \ \ \ \ \ \ \ eqnRow\ =\ append(eqnRow,\ i)\leavevmode\hss\endgraf
\ \ \ \ \ \ \ \ \ \ \ \ \ \ \ \ \ \ \ \ \ else\ ineqnRow\ =\ append(ineqnRow,\ i);));\ \ \ \ \leavevmode\hss\endgraf
\ \ \ \ \ \ \ \ \ \ \ E\ :=\ apply(eqnRow,\ i\ ->\ A{\char`\#}i);\leavevmode\hss\endgraf
\ \ \ \ \ \ \ \ \ \ \ E\ =\ apply(E,\ e\ ->\ e{\char`\_}{\char`\{}n..(n+d-1){\char`\}});\leavevmode\hss\endgraf
\ \ \ \ \ \ \ \ \ \ \ A\ =\ apply(ineqnRow,\ i\ ->\ A{\char`\#}i);\leavevmode\hss\endgraf
\ \ \ \ \ \ \ \ \ \ \ A\ =\ apply(A,\ e\ ->\ e{\char`\_}(T\ |\ toList(n..(n+d-1))));\ \leavevmode\hss\endgraf
\ \ \ \ \ \ \ \ \ \ \ --\ successive\ projections\ eliminate\ the\ variables\ 'T'.\leavevmode\hss\endgraf
\ \ \ \ \ \ \ \ \ \ \ if\ A\ =!=\ {\char`\{}{\char`\}}\ then\ scan(T,\ t\ ->\ (\leavevmode\hss\endgraf
\ \ \ \ \ \ \ \ \ \ \ \ \ \ \ \ \ \ \ \ \ D\ :=\ fourierMotzkin(A,\ V,\ t);\leavevmode\hss\endgraf
\ \ \ \ \ \ \ \ \ \ \ \ \ \ \ \ \ \ \ \ \ A\ =\ D{\char`\#}0;\ \ \ \ \ \ \ \ \ \ V\ =\ D{\char`\#}1;));\leavevmode\hss\endgraf
\ \ \ \ \ \ \ \ \ \ \ --\ output\ formating\leavevmode\hss\endgraf
\ \ \ \ \ \ \ \ \ \ \ A\ =\ apply(A,\ e\ ->\ primitive\ e);\leavevmode\hss\endgraf
\ \ \ \ \ \ \ \ \ \ \ if\ A\ ===\ {\char`\{}{\char`\}}\ then\ A\ =\ map(ZZ{\char`\^}d,\ ZZ{\char`\^}0,\ 0)\leavevmode\hss\endgraf
\ \ \ \ \ \ \ \ \ \ \ else\ A\ =\ transpose\ matrix\ A;\leavevmode\hss\endgraf
\ \ \ \ \ \ \ \ \ \ \ if\ E\ ===\ {\char`\{}{\char`\}}\ then\ E\ =\ map(ZZ{\char`\^}d,\ ZZ{\char`\^}0,\ 0)\leavevmode\hss\endgraf
\ \ \ \ \ \ \ \ \ \ \ else\ E\ =\ transpose\ matrix\ E;\leavevmode\hss\endgraf
\ \ \ \ \ \ \ \ \ \ \ (A,\ E));\ \leavevmode\hss\endgraf
\endgroup
\penalty-1000
\par
\vskip 1 pt
\noindent

If the input matrix $H$ has no columns, it can be omitted.  A sequence of two
matrices is returned, as above.
\par
\vskip 5 pt
\begingroup
\tteight
\baselineskip=7.4pt
\lineskip=0pt
\obeyspaces
i91\ :\ polarCone(Matrix)\ :=\ (Z)\ ->\ (\leavevmode\hss\endgraf
\ \ \ \ \ \ \ \ \ \ \ polarCone(Z,\ map(ZZ{\char`\^}(rank\ target\ Z),\ ZZ{\char`\^}0,\ 0)));\leavevmode\hss\endgraf
\endgroup
\penalty-1000
\par
\vskip 1 pt
\noindent

As a simple example, consider the permutahedron in $\R^3$ 
whose vertices are the following six points. 

\par
\vskip 5 pt
\begingroup
\tteight
\baselineskip=7.4pt
\lineskip=0pt
\obeyspaces
i92\ :\ H\ =\ transpose\ matrix{\char`\{}\leavevmode\hss\endgraf
\ \ \ \ \ \ {\char`\{}1,2,3{\char`\}},\leavevmode\hss\endgraf
\ \ \ \ \ \ {\char`\{}1,3,2{\char`\}},\leavevmode\hss\endgraf
\ \ \ \ \ \ {\char`\{}2,1,3{\char`\}},\leavevmode\hss\endgraf
\ \ \ \ \ \ {\char`\{}2,3,1{\char`\}},\leavevmode\hss\endgraf
\ \ \ \ \ \ {\char`\{}3,1,2{\char`\}},\leavevmode\hss\endgraf
\ \ \ \ \ \ {\char`\{}3,2,1{\char`\}}{\char`\}};\leavevmode\hss\endgraf
\penalty-500\leavevmode\hss\endgraf
\ \ \ \ \ \ \ \ \ \ \ \ \ \ \ 3\ \ \ \ \ \ \ \ 6\leavevmode\hss\endgraf
o92\ :\ Matrix\ ZZ\ \ <---\ ZZ\leavevmode\hss\endgraf
\endgroup
\penalty-1000
\par
\vskip 1 pt
\noindent

The inequality representation of the permutahedron is obtained 
by calling {\tt polarCone} on $H$: the facet normals of the 
polytope are the columns of the matrix in the first argument of the 
output. The second argument is trivial since our input is a polytope 
and hence there are is no non-trivial affine space contained in it.
If we call {\tt polarCone} on the output, we will get back H as
expected. 

\par
\vskip 5 pt
\begingroup
\tteight
\baselineskip=7.4pt
\lineskip=0pt
\obeyspaces
i93\ :\ P\ =\ polarCone(H)\leavevmode\hss\endgraf
\penalty-500\leavevmode\hss\endgraf
o93\ =\ (|\ 1\ \ 1\ \ 1\ \ -1\ -1\ -5\ |,\ 0)\leavevmode\hss\endgraf
\ \ \ \ \ \ \ |\ -1\ 1\ \ -5\ 1\ \ -1\ 1\ \ |\leavevmode\hss\endgraf
\ \ \ \ \ \ \ |\ -1\ -5\ 1\ \ -1\ 1\ \ 1\ \ |\leavevmode\hss\endgraf
\penalty-500\leavevmode\hss\endgraf
o93\ :\ Sequence\leavevmode\hss\endgraf
\endgroup
\penalty-1000
\par
\vskip 1 pt
\noindent
\par
\vskip 5 pt
\begingroup
\tteight
\baselineskip=7.4pt
\lineskip=0pt
\obeyspaces
i94\ :\ Q\ =\ polarCone(P{\char`\_}0)\leavevmode\hss\endgraf
\penalty-500\leavevmode\hss\endgraf
o94\ =\ (|\ 1\ 1\ 2\ 2\ 3\ 3\ |,\ 0)\leavevmode\hss\endgraf
\ \ \ \ \ \ \ |\ 2\ 3\ 1\ 3\ 1\ 2\ |\leavevmode\hss\endgraf
\ \ \ \ \ \ \ |\ 3\ 2\ 3\ 1\ 2\ 1\ |\leavevmode\hss\endgraf
\penalty-500\leavevmode\hss\endgraf
o94\ :\ Sequence\leavevmode\hss\endgraf
\endgroup
\penalty-1000
\par
\vskip 1 pt
\noindent

\section{Appendix 2: Minimal presentation of rings}

Throughout this chapter, we have used on several occasions the simple, yet
useful subroutine {\tt removeRedundantVariables}.
In this appendix, we present \Mtwo code for this routine,
which is the main ingredient for finding minimal
presentations of quotients of polynomial rings.
Our code for this routine is a somewhat simplified, but less
efficient version of a routine in the \Mtwo package, {\tt minPres.m2},
written by Amelia Taylor.

The routine {\tt removeRedundantVariables} takes as input an ideal {\tt I} in
a polynomial ring {\tt A}.  It returns a ring map {\tt F} from {\tt A} to
itself, which sends redundant variables to polynomials in the non-redundant
variables, and sends non-redundant variables to themselves.  For example,
\medskip
  \par
\vskip 5 pt
\begingroup
\tteight
\baselineskip=7.4pt
\lineskip=0pt
\obeyspaces
i95\ :\ A\ =\ QQ[a..e];\leavevmode\hss\endgraf
\endgroup
\penalty-1000
\par
\vskip 1 pt
\noindent
  \par
\vskip 5 pt
\begingroup
\tteight
\baselineskip=7.4pt
\lineskip=0pt
\obeyspaces
i96\ :\ I\ =\ ideal(a-b{\char`\^}2-1,\ b-c{\char`\^}2,\ c-d{\char`\^}2,\ a{\char`\^}2-e{\char`\^}2)\leavevmode\hss\endgraf
\penalty-500\leavevmode\hss\endgraf
\ \ \ \ \ \ \ \ \ \ \ \ \ \ \ \ 2\ \ \ \ \ \ \ \ \ \ \ \ \ 2\ \ \ \ \ \ \ \ \ 2\ \ \ \ \ \ \ 2\ \ \ \ 2\leavevmode\hss\endgraf
o96\ =\ ideal\ (-\ b\ \ +\ a\ -\ 1,\ -\ c\ \ +\ b,\ -\ d\ \ +\ c,\ a\ \ -\ e\ )\leavevmode\hss\endgraf
\penalty-500\leavevmode\hss\endgraf
o96\ :\ Ideal\ of\ A\leavevmode\hss\endgraf
\endgroup
\penalty-1000
\par
\vskip 1 pt
\noindent
  \par
\vskip 5 pt
\begingroup
\tteight
\baselineskip=7.4pt
\lineskip=0pt
\obeyspaces
i97\ :\ F\ =\ removeRedundantVariables\ I\leavevmode\hss\endgraf
\penalty-500\leavevmode\hss\endgraf
\ \ \ \ \ \ \ \ \ \ \ \ \ \ \ \ 8\ \ \ \ \ \ \ 4\ \ \ 2\leavevmode\hss\endgraf
o97\ =\ map(A,A,{\char`\{}d\ \ +\ 1,\ d\ ,\ d\ ,\ d,\ e{\char`\}})\leavevmode\hss\endgraf
\penalty-500\leavevmode\hss\endgraf
o97\ :\ RingMap\ A\ <---\ A\leavevmode\hss\endgraf
\endgroup
\penalty-1000
\par
\vskip 1 pt
\noindent
The non-redundant variables are $d$ and $e$.  The image of $I$ under $F$
gives the elements in this smaller set of variables.  We take the ideal of a 
Gr\"obner basis of the image:
  \par
\vskip 5 pt
\begingroup
\tteight
\baselineskip=7.4pt
\lineskip=0pt
\obeyspaces
i98\ :\ I1\ =\ ideal\ gens\ gb(F\ I)\leavevmode\hss\endgraf
\penalty-500\leavevmode\hss\endgraf
\ \ \ \ \ \ \ \ \ \ \ \ \ 16\ \ \ \ \ 8\ \ \ \ 2\leavevmode\hss\endgraf
o98\ =\ ideal(d\ \ \ +\ 2d\ \ -\ e\ \ +\ 1)\leavevmode\hss\endgraf
\penalty-500\leavevmode\hss\endgraf
o98\ :\ Ideal\ of\ A\leavevmode\hss\endgraf
\endgroup
\penalty-1000
\par
\vskip 1 pt
\noindent
The original ideal can be written in a cleaner way as
  \par
\vskip 5 pt
\begingroup
\tteight
\baselineskip=7.4pt
\lineskip=0pt
\obeyspaces
i99\ :\ ideal\ compress\ (F.matrix\ -\ vars\ A)\ +\ I1\leavevmode\hss\endgraf
\penalty-500\leavevmode\hss\endgraf
\ \ \ \ \ \ \ \ \ \ \ \ \ \ 8\ \ \ \ \ \ \ \ \ \ \ 4\ \ \ \ \ \ \ 2\ \ \ \ \ \ \ 16\ \ \ \ \ 8\ \ \ \ 2\leavevmode\hss\endgraf
o99\ =\ ideal\ (d\ \ -\ a\ +\ 1,\ d\ \ -\ b,\ d\ \ -\ c,\ d\ \ \ +\ 2d\ \ -\ e\ \ +\ 1)\leavevmode\hss\endgraf
\penalty-500\leavevmode\hss\endgraf
o99\ :\ Ideal\ of\ A\leavevmode\hss\endgraf
\endgroup
\penalty-1000
\par
\vskip 1 pt
\noindent
  
  Let us now describe the \Mtwo code.  The subroutine {\tt
    findRedundant} takes a polynomial $f$, and finds a variable $x_i$
  in the ring of $f$, such that $f = c x_i + g$, for a non-zero
  constant $c$, and a polynomial $g$ which does not involve the
  variable $x_i$.  If there is no such variable, {\tt null} is
  returned.  Otherwise, if $x_i$ is the first such variable , the list
  $\{i, c^{-1} g\}$ is returned.

\par
\vskip 5 pt
\begingroup
\tteight
\baselineskip=7.4pt
\lineskip=0pt
\obeyspaces
i100\ :\ findRedundant=(f)->(\leavevmode\hss\endgraf
\ \ \ \ \ \ \ \ \ \ \ \ A\ :=\ ring(f);\leavevmode\hss\endgraf
\ \ \ \ \ \ \ \ \ \ \ \ p\ :=\ first\ entries\ contract(vars\ A,f);\leavevmode\hss\endgraf
\ \ \ \ \ \ \ \ \ \ \ \ i\ :=\ position(p,\ g\ ->\ g\ !=\ 0\ and\ first\ degree\ g\ ===\ 0);\leavevmode\hss\endgraf
\ \ \ \ \ \ \ \ \ \ \ \ if\ i\ ===\ null\ then\leavevmode\hss\endgraf
\ \ \ \ \ \ \ \ \ \ \ \ \ \ \ \ null\leavevmode\hss\endgraf
\ \ \ \ \ \ \ \ \ \ \ \ else\ (\leavevmode\hss\endgraf
\ \ \ \ \ \ \ \ \ \ \ \ \ \ \ \ \ v\ :=\ A{\char`\_}i;\leavevmode\hss\endgraf
\ \ \ \ \ \ \ \ \ \ \ \ \ \ \ \ \ c\ :=\ f{\char`\_}v;\leavevmode\hss\endgraf
\ \ \ \ \ \ \ \ \ \ \ \ \ \ \ \ \ {\char`\{}i,(-1)*(c{\char`\^}(-1)*(f-c*v)){\char`\}}\leavevmode\hss\endgraf
\ \ \ \ \ \ \ \ \ \ \ \ \ \ \ \ \ )\leavevmode\hss\endgraf
\ \ \ \ \ \ \ \ \ \ \ \ );\leavevmode\hss\endgraf
\endgroup
\penalty-1000
\par
\vskip 1 pt
\noindent

The main function {\tt removeRedundantVariables} requires an ideal in a
polynomial ring (not a quotient ring) as input.  The internal
routine {\tt findnext} finds the first entry of the (one row) matrix {\tt M}
which contains a redundancy.  This redundancy is used to modify the list {\tt
xmap}, which contains the images of the redundant variables.
The matrix {\tt M}, and the list {\tt xmap} are both updated, and 
then we continue to look for more redundancies.

\par
\vskip 5 pt
\begingroup
\tteight
\baselineskip=7.4pt
\lineskip=0pt
\obeyspaces
i101\ :\ removeRedundantVariables\ =\ (I)\ ->\ (\leavevmode\hss\endgraf
\ \ \ \ \ \ \ \ \ \ \ \ A\ :=\ ring\ I;\leavevmode\hss\endgraf
\ \ \ \ \ \ \ \ \ \ \ \ xmap\ :=\ new\ MutableList\ from\ gens\ A;\ \ \ \ \ \ \ \leavevmode\hss\endgraf
\ \ \ \ \ \ \ \ \ \ \ \ M\ :=\ gens\ I;\leavevmode\hss\endgraf
\ \ \ \ \ \ \ \ \ \ \ \ findnext\ :=\ ()\ ->\ (\leavevmode\hss\endgraf
\ \ \ \ \ \ \ \ \ \ \ \ \ \ \ \ \ p\ :=\ null;\leavevmode\hss\endgraf
\ \ \ \ \ \ \ \ \ \ \ \ \ \ \ \ \ next\ :=\ 0;\leavevmode\hss\endgraf
\ \ \ \ \ \ \ \ \ \ \ \ \ \ \ \ \ done\ :=\ false;\leavevmode\hss\endgraf
\ \ \ \ \ \ \ \ \ \ \ \ \ \ \ \ \ ngens\ :=\ numgens\ source\ M;\leavevmode\hss\endgraf
\ \ \ \ \ \ \ \ \ \ \ \ \ \ \ \ \ while\ next\ <\ ngens\ and\ not\ done\ do\ (\leavevmode\hss\endgraf
\ \ \ \ \ \ \ \ \ \ \ \ \ \ \ \ \ \ \ p\ =\ findRedundant(M{\char`\_}(0,next));\leavevmode\hss\endgraf
\ \ \ \ \ \ \ \ \ \ \ \ \ \ \ \ \ \ \ if\ p\ =!=\ null\ then\leavevmode\hss\endgraf
\ \ \ \ \ \ \ \ \ \ \ \ \ \ \ \ \ \ \ \ \ \ \ \ done\ =\ true\leavevmode\hss\endgraf
\ \ \ \ \ \ \ \ \ \ \ \ \ \ \ \ \ \ \ else\ next=next+1;\leavevmode\hss\endgraf
\ \ \ \ \ \ \ \ \ \ \ \ \ \ \ \ \ );\leavevmode\hss\endgraf
\ \ \ \ \ \ \ \ \ \ \ \ \ \ \ \ \ p);\leavevmode\hss\endgraf
\ \ \ \ \ \ \ \ \ \ \ \ p\ :=\ findnext();\leavevmode\hss\endgraf
\ \ \ \ \ \ \ \ \ \ \ \ while\ p\ =!=\ null\ do\ (\leavevmode\hss\endgraf
\ \ \ \ \ \ \ \ \ \ \ \ \ \ \ \ \ xmap{\char`\#}(p{\char`\#}0)\ =\ p{\char`\#}1;\leavevmode\hss\endgraf
\ \ \ \ \ \ \ \ \ \ \ \ \ \ \ \ \ F1\ :=\ map(A,A,toList\ xmap);\leavevmode\hss\endgraf
\ \ \ \ \ \ \ \ \ \ \ \ \ \ \ \ \ F2\ :=\ map(A,A,\ F1\ (F1.matrix));\leavevmode\hss\endgraf
\ \ \ \ \ \ \ \ \ \ \ \ \ \ \ \ \ xmap\ =\ new\ MutableList\ from\ first\ entries\ F2.matrix;\leavevmode\hss\endgraf
\ \ \ \ \ \ \ \ \ \ \ \ \ \ \ \ \ M\ =\ compress(F2\ M);\leavevmode\hss\endgraf
\ \ \ \ \ \ \ \ \ \ \ \ \ \ \ \ \ p\ =\ findnext();\leavevmode\hss\endgraf
\ \ \ \ \ \ \ \ \ \ \ \ \ \ \ \ \ );\leavevmode\hss\endgraf
\ \ \ \ \ \ \ \ \ \ \ \ map(A,A,toList\ xmap));\leavevmode\hss\endgraf
\endgroup
\penalty-1000
\par
\vskip 1 pt
\noindent
\bibliography{papers}
\end{document}